%%%%%%%%%%%%%%%%%AMSTeXfile%%%%%%%%%%%%%%%%%%%%%%%%%%%%%%%%%%%%%%%%%
\input amstex.tex
\input epsf
\documentstyle{amsppt}
\magnification=1200
\baselineskip=13pt
\hsize=6.5truein
\vsize=8.9truein
\parindent=20pt
%\NoBlackBoxes
%%%%%%%%%%%%%%%%% M a c r o s %%%%%%%%%%%
%Section, Theorem and Formula Numbering%%
%%%%%%%%%%%%%%%%%%%%%%%%%%%%%%%%%%%%%%%%%
\newcount\sectionnumber
\newcount\equationnumber
\newcount\thnumber
\newcount\countrefno

\def\ifundefined#1{\expandafter\ifx\csname#1\endcsname\relax}
\def\assignnumber#1#2{%
	\ifundefined{#1}\relax\else\message{#1 already defined}\fi
	\expandafter\xdef\csname#1\endcsname
        {\the\sectionnumber.\the#2}}%
%
% macros on section numbers
%
\def\newsec{
  \global\advance\sectionnumber by 1
  \global\equationnumber=0
  \global\thnumber=0
  \the\sectionnumber .\ }
\def\newappA{
  \global\advance\sectionnumber by 1
  \global\equationnumber=0
  \global\thnumber=0
  Appendix A.\ }
\def\newappB{
  \global\advance\sectionnumber by 1
  \global\equationnumber=0
  \global\thnumber=0
  Appendix B.\ }
%  
%\comment
\def\eq#1{\relax
  \global\advance\equationnumber by 1
  \assignnumber{EN#1}\equationnumber
  {\rm \csname EN#1\endcsname}}
\def\eqtag#1{\ifundefined{EN#1}\message{EN#1 undefined}{\sl (#1)}%
  \else\thetag{\csname EN#1\endcsname}\fi}
%
% macros on theorem numbers
%
\def\thname#1{\relax
  \global\advance\thnumber by 1
  \assignnumber{TH#1}\thnumber
  \csname TH#1\endcsname}
\def\thtag#1{\ifundefined{TH#1}\message{TH#1 undefined}{\sl #1}%
  \else\csname TH#1\endcsname\fi}
%
%same macros for Appendix A
%  
\def\Aassignnumber#1#2{%
	\ifundefined{#1}\relax\else\message{#1 already defined}\fi
	\expandafter\xdef\csname#1\endcsname
        {A.\the#2}}
\def\Aeq#1{\relax
  \global\advance\equationnumber by 1
  \Aassignnumber{EN#1}\equationnumber
  {\rm \csname EN#1\endcsname}}
\def\Athname#1{\relax
  \global\advance\thnumber by 1
  \Aassignnumber{TH#1}\thnumber
  \csname TH#1\endcsname} 
%
%same macros for Appendix B
%  
\def\Bassignnumber#1#2{%
	\ifundefined{#1}\relax\else\message{#1 already defined}\fi
	\expandafter\xdef\csname#1\endcsname
        {B.\the#2}}
\def\Beq#1{\relax
  \global\advance\equationnumber by 1
  \Bassignnumber{EN#1}\equationnumber
  {\rm \csname EN#1\endcsname}}
\def\Bthname#1{\relax
  \global\advance\thnumber by 1
  \Bassignnumber{TH#1}\thnumber
  \csname TH#1\endcsname} 
%
%\endcomment
%not using automatic formula numbering; comment this otherwise
\comment
\def\eq{}
\def\eqtag#1{(#1)}
\def\thname{}
\def\thtag{}
\def\Aeq{}
\def\Athname{}
\def\Beq{}
\def\Bthname{}
\endcomment
%not using automatic formula numbering; comment this otherwise
%
% macros on reference numbers
%
\global\countrefno=1

\def\refno#1{\xdef#1{{\the\countrefno}}
\global\advance\countrefno by 1}
%%%%Abbreviations%%%%%%%%%%%%%%%%%%%%%%%%%%%
\def\R{{\Bbb R}}
\def\N{{\Bbb N}}
\def\C{{\Bbb C}}

\def\Z{{\Bbb Z}}
\def\T{{\Bbb T}}

\def\Zp{{\Bbb Z}_{\geq 0}}
\def\hZp{{1\over 2}\Zp}
\def\hf{{1\over 2}}

\def\al{\alpha}
\def\be{\beta}
\def\ga{\gamma}
\def\Ga{\Gamma}
\def\de{\delta}
\def\De{\Delta}

\def\ep{\varepsilon}
\def\si{\sigma}

\def\th{\theta}
\def\la{\lambda}
\def\vp{\varphi}
\def\om{\omega}
\def\Om{\Omega}
\def\H{\ell^2(\Z)}
\def\id{\text{\rm id}}
\def\pt{\pi_{e^{i\theta}}}
\def\pp{\pi_{e^{i\psi}}}
\def\pf{\pi_{e^{i\phi}}}
\def\rst{\rho_{s,t}}
\def\rit{\rho_{\infty,t}}
\def\U{U_q({\frak{su}}(1,1))}
\def\Uc{U_q({\frak{su}}(2))}
\def\UC{U_q({\frak{sl}}(2,\C))}
%%%Reference Numbering%%%%%%%%%%%%%%%%%%%%%%%%
\refno{\Akhi}
\refno{\AlSaC}
\refno{\AskeI}
\refno{\AskeRS}
\refno{\AskeW}
\refno{\Baaj}
\refno{\Bere}
\refno{\BurbK}
\refno{\BustS}
\refno{\CharP}
\refno{\CiccKK}
\refno{\Comb}
\refno{\DijkN}
\refno{\DijkS}
\refno{\DixmvNA}
\refno{\DunfS}
\refno{\GaspR}
\refno{\GuptIM}
\refno{\IsmaR}
\refno{\IsmaS}
\refno{\Kake}
\refno{\KakeMU}
\refno{\KoelPhD}
\refno{\KoelITSF}
\refno{\KoelSIAM}
\refno{\KoelAAM}
\refno{\KoelFIC}
\refno{\KoelSbig}
\refno{\KoelSAW}
\refno{\KoelVdJ}
\refno{\KoelV}
\refno{\KoorJF}
\refno{\KoorAF}
\refno{\KoorCM}
\refno{\KoorZSE}
\refno{\Kust}
\refno{\KustV}
\refno{\KustVCR}
\refno{\Lanc}
\refno{\MassR}
\refno{\MasuMNNSU}
\refno{\MasuW}
\refno{\Noum} 
\refno{\NoumAM}
\refno{\NoumMPJA}
\refno{\NoumMLNM}
\refno{\Pede}
\refno{\QuaeV}
\refno{\Rahm}
\refno{\Rose}
\refno{\Rudi}
\refno{\Schm}
\refno{\Sheu}
\refno{\Simo}
\refno{\SuslJPA}
\refno{\SuslPP}
\refno{\VaksK}
\refno{\VandJ}
\refno{\Verd}
\refno{\Vile}
\refno{\VileK}
\refno{\WoroEtwee}
%%%%%%%%%%%%%%%%%%%%%%%%%%%%%%%%%%%%%%%%%%%%%%%%%%%%%%%%%%%%%%%%%%%%
%Beginning of Text%
%%%%%%%%%%%%%%%%%%%%%%%%%%%%%%%%%%%%%%%%%%%%%%%%%%%%%%%%%%%%%%%%%%%%
\topmatter
\title Fourier transforms on the quantum $SU(1,1)$ group\endtitle
\author Erik Koelink and Jasper V.~Stokman\\ 
with an appendix by Mizan Rahman
\endauthor
\leftheadtext{Erik Koelink and Jasper Stokman}
\rightheadtext{Fourier transforms on the quantum $SU(1,1)$ group}
\address E.K.: Technische Universiteit Delft, Faculteit
Informatietechnologie en Systemen, Afd. Toegepaste Wiskundige
Analyse, Postbus 5031, 2600 GA Delft, the Netherlands\endaddress
\email koelink\@twi.tudelft.nl\endemail
\address J.S.: Centre de Math\'ematiques de Jussieu, Universit\'e 
Paris 6 Pierre et Marie Curie, 4, place Jussieu, 
F-75252 Paris Cedex 05, France\endaddress
\email stokman\@math.jussieu.fr\endemail
\address M.R.: Carleton University, Department of Mathematics and
Statistics, 1125 Colonel By Drive, Ottawa, Ontario K1S\, 5B6,
Canada\endaddress
\email mrahman\@math.carleton.ca\endemail
\date November 2, 1999\enddate
\abstract 
The main goal is to interpret the Askey-Wilson function and the
corresponding transform pair on the quantum $SU(1,1)$ group. 
A weight on the $\text{C}^\ast$-algebra of continuous functions
vanishing at infinity on the quantum $SU(1,1)$ group is
studied, which is left and right invariant in a weak sense 
with respect to a product defined using Wall functions.  
The Haar
weight restricted to certain subalgebras are explicitly
determined in terms of an infinitely supported Jackson integral and
in terms of an infinitely supported Askey-Wilson type measure. 
For the evaluation the spectral analysis of explicit 
unbounded doubly infinite Jacobi matrices and
some new summation formulas for basic hypergeometric series
are needed.  
The spherical functions are calculated in terms of Askey-Wilson
functions and big $q$-Jacobi functions. The 
corresponding spherical Fourier transforms are identified with
special cases of the big $q$-Jacobi function transform
and of the Askey-Wilson function transform. 
\endabstract
\endtopmatter
\document

%%%%%%%%%%%%%%%%%%%%%%%%%%%%%%%%%%%%%%%%%%%%%%%%%%%%%%%%%%%%%%%%%%%%
%N E W   S E C T I O N%
%%%%%%%%%%%%%%%%%%%%%%%%%%%%%%%%%%%%%%%%%%%%%%%%%%%%%%%%%%%%%%%%%%%%
\head \newsec Introduction\endhead

The motivation for the study in this paper is twofold. On the one
hand we are interested in the study of the simplest non-compact
semisimple quantum group, namely the quantum $SU(1,1)$ group, and in
particular in its corresponding Haar functional. This quantum group
is resisting any of the theories on locally compact quantum groups
like e.g. \cite{\KustVCR}.
On the other hand
we are interested in special functions associated to quantum groups,
and in particular in the so-called Askey-Wilson functions. 
Let us first say something on the second subject, which is our main
concern. 

A very general set of orthogonal polynomials in 
one variable is the
set of Askey-Wilson polynomials introduced in 1985 in \cite{\AskeW}.
As the title of the memoir indicates, Askey-Wilson polynomials can
be considered as $q$-analogues of the Jacobi 
polynomials which are
orthogonal on $[-1,1]$ with respect to the beta integral
$(1-x)^\al(1+x)^\be$. The Jacobi polynomials are the polynomial
solutions of the hypergeometric differential operator, whereas 
the Askey-Wilson polynomials are the polynomial
solutions of a certain second-order difference operator. 
The Jacobi polynomials naturally arise as spherical functions 
on rank one compact Riemannian symmetric spaces. On the other hand, 
the spherical functions on non-compact rank one Riemannian symmetric
spaces can be expressed in terms of Jacobi functions, which
are non-polynomial eigenfunctions of the hypergeometric differential
operator. The corresponding Fourier transforms are special cases of
the Jacobi function transform in which the
kernel is a Jacobi function. By now 
these Jacobi function transforms, containing as special cases
the Fourier-cosine and Mehler-Fock transforms, are very well
understood, see e.g. the survey paper \cite{\KoorJF} by Koornwinder
and references therein. There are inversion formulas, as well as
the appropriate
analogues of the theorems of Plancherel, Parseval and Paley-Wiener. 
Furthermore, there are several different approaches to the study of
the $L^2$-theory of the Jacobi function transform. One particular
approach is by  spectral analysis of the hypergeometric 
differential operator on a
weighted $L^2$-space. 

Although the Askey-Wilson functions are known, see \cite{\IsmaR},
\cite{\Rahm}, \cite{\SuslJPA}, \cite{\SuslPP}, as are all the
solutions to the Askey-Wilson second order difference equation 
and their interrelations, it was not yet known what the
appropriate Askey-Wilson function transform should be. 
The reason for this is our lack of understanding of the 
Hilbert space on which the Askey-Wilson difference operator  
has to be diagonalised. In case the Jacobi functions have an
interpretation as spherical functions, the weighted $L^2$-space 
can be obtained by restricting the Haar measure 
to functions which behave as a character under the left and right
action of a maximal compact subgroup. 
In this paper we show how the
study of the Haar functional on the quantum $SU(1,1)$ group can be
used to find the right Hilbert spaces 
as weighted $L^2$-space and this is one of 
the main results of this paper. The actual analytic study of the
Askey-Wilson function transform is done in another paper 
\cite{\KoelSAW}, and
that of the appropriate limit case to the 
big $q$-Jacobi transform is done in \cite{\KoelSbig}. 
This is for two reasons. The quantum group theoretic approach 
does not lead to a rigorous proof, and secondly the interpretation
on the quantum group only holds for a restricted set of the
parameters involved. 
It has to be noted that the Askey-Wilson function transform that
occurs in this paper is different from the orthogonality 
relations introduced by Suslov \cite{\SuslJPA}, \cite{\SuslPP}, 
see also \cite{\BustS} for the more extensively studied little 
$q$-Jacobi case, which is analogous to Fourier(-Bessel) series. 

The motivation for the method we employ is the relation between
special functions and the theory of group (and quantum group) 
representations. The Jacobi polynomials occur as matrix coefficients
of irreducible unitary representations 
of the compact Lie group $SU(2)$ and the Jacobi functions arise as
matrix coefficients of irreducible unitary representations
of the non-compact Lie group $SU(1,1)\cong
SL(2,\R)$, see \cite{\Vile}, \cite{\VileK}, \cite{\KoorJF}. These
groups are both real forms of the same complex Lie group $SL(2,\C)$.
In the theory of quantum groups, the quantum analogue of the complex
case $SL(2,\C)$ is much studied, as is the quantum analogue of the
compact $SU(2)$, see e.g. \cite{\CharP}. One of the first
indications that the relation between quantum groups and special
functions is very strong, is the interpretation of the
little $q$-Jacobi polynomials on the quantum $SU(2)$ group as matrix
elements on which the subgroup $K=S(U(1)\times U(1))\cong U(1)$ 
acts by a character. 
Since we can view little $q$-Jacobi polynomials as
limiting cases of the Askey-Wilson polynomials, this is a first
step. The breakthrough has come with Koornwinder's paper 
\cite{\KoorZSE} in which he gives an 
infinitesimal characterisation of quantum subgroups. 
This gives a one-parameter family of quantum subgroups, 
which we denote by $K_t$. The subgroups $K_t$ and $K_s$
are formally conjugated, see \cite{\Rose, \S 4}. 
Then the matrix elements on which $K_s$, respectively $K_t$,
acts by a character from the left, respectively right, can be
expressed in terms of Askey-Wilson polynomials. 
The in-between case of the big $q$-Jacobi polynomials can 
be obtained in a similar way.
As a corollary to these results we get from the Schur  
orthogonality relations an explicit expression for
the Haar functional on certain commutative subalgebras 
in terms of the Askey-Wilson orthogonality measure. 
For the spherical case, i.e. the matrix elements that are left
and right invariant under $K$, we state this in
the following table for the quantum $SU(2)$ group. The spherical
case is the important case to calculate. 

\vskip10pt\noindent
\vbox{\offinterlineskip
\hrule
\halign{&\vrule#&\strut\quad\hfil#\hfil\quad\cr
height2pt&\omit&&\omit&&\omit&\cr
&{\sl subgroups}&&{\sl Haar functional}&&{\sl spherical 
functions}&\cr
height2pt&\omit&&\omit&&\omit&\cr
\noalign{\hrule}
height2pt&\omit&&\omit&&\omit&\cr
&$(K,K)$&&Jackson integral on $[0,1]$&&little 
$q$-Legendre polynomials&\cr
&$(K,K_t)$&&Jackson integral on $[-t,1]$&&big  
$q$-Legendre polynomials&\cr
&$(K_s,K_t)$&&Askey-Wilson integral&&2-parameter   
Askey-Wilson polynomials&\cr
height2pt&\omit&&\omit&&\omit&\cr}
\hrule}
\vskip-10pt
\botcaption{Table 1.1} Spherical functions for $SU_q(2)$
\endcaption
\vskip10pt

For the quantum $SU(1,1)$ group the matrix elements that behave
nicely under the action of the subgroup 
$K=S(U(1)\times U(1))\cong U(1)$ have been
calculated explicitly by Masuda et al. \cite{\MasuMNNSU} and 
Vaksman and Korogodski\u\i\ \cite{\VaksK}. These can be expressed in
terms of little $q$-Jacobi functions, and the Haar functional 
is also known in terms of a Jackson integral on $[0,\infty)$, 
see \cite{\Kake}, \cite{\KakeMU}, \cite{\MasuW}, \cite{\VaksK}. 
This gives rise to the first line in the following table. 

\vskip10pt\noindent
\vbox{\offinterlineskip
\hrule
\halign{&\vrule#&\strut\quad\hfil#\hfil\quad\cr
height2pt&\omit&&\omit&&\omit&\cr
&{\sl subgroups}&&{\sl Haar functional}&&{\sl spherical 
functions}&\cr
height2pt&\omit&&\omit&&\omit&\cr
\noalign{\hrule}
height2pt&\omit&&\omit&&\omit&\cr
&$(K,K)$&&Jackson integral on $[0,\infty)$&&little 
$q$-Legendre functions&\cr
&$(K,K_t)$&&Jackson integral on $[-t,\infty)$&&big  
$q$-Legendre functions&\cr
&$(K_s,K_t)$&&Askey-Wilson type integral&&2-parameter  
Askey-Wilson functions&\cr
height2pt&\omit&&\omit&&\omit&\cr}
\hrule}
\vskip-10pt
\botcaption{Table 1.2} Spherical functions for $SU_q(1,1)$
\endcaption
\vskip10pt

The purpose of this paper is to prove the last two lines of 
Table~1.2. The
proof of the explicit expression for the Haar functional in the last
two cases of Table~1.2 is the main result of this paper.
Koornwinder's proof for the cases in the compact setting listed in
Table~1.1 cannot be used here, but the alternative proof using
spectral theory and bilinear generating functions given in
\cite{\KoelV} can be generalised to the quantum $SU(1,1)$ group. 
For this we give an expression for a Haar functional 
on the quantum $SU(1,1)$ group 
in terms of representations of the quantised function algebra. 
In the last section we then formally show how the big
$q$-Jacobi function transform and the Askey-Wilson function 
transform can be interpreted as Fourier transforms on the quantum 
$SU(1,1)$ group. Because the Fourier transforms 
associated with $SU(1,1)$ are special case of the Jacobi
function transforms, see \cite{\KoorJF}, we
view the big $q$-Jacobi and Askey-Wilson function transform as
$q$-analogues of the Jacobi function transform. 
The complete analytic study of the big
$q$-Jacobi and Askey-Wilson function transform is developed in
\cite{\KoelSbig} and \cite{\KoelSAW}.

We expect that the Askey-Wilson function transform will play a
central role in the theory of integral transforms with basic
hypergeometric kernels. Indeed, in the polynomial setting,
the Askey-Wilson polynomials have had a tremendous impact in the
theory of basic hypergeometric orthogonal polynomials. 
Furthermore, the Jacobi function transform, which we consider
as the classical counterpart of the Askey-Wilson function 
transform, has turned out to be an important integral transform
in the theory of special functions and its applications.
In particular, due 
to the quantum group theoretic interpretation of Askey-Wilson
functions in this paper, we may expect that 
appropriate non-polynomial analogues
of the results on Askey-Wilson polynomials 
in e.g. \cite{\DijkN}, \cite{\DijkS}, 
\cite{\KoelFIC}, \cite{\KoorCM}, \cite{\NoumAM} exist. 

The theory of locally compact quantum groups has
not yet reached the state of maturity, but see Kustermans and Vaes
\cite{\KustV}, \cite{\KustVCR}. The quantum $SU(1,1)$ group does not
fit into these theories because it lacks a good definition of the
comultiplication defined on the ${\text C}^\ast$-algebra level, see
\cite{\WoroEtwee}, and without a comultiplication it is not possible
to speak of left- and right invariance of a functional. In \S 2 we
propose a weak version of the comultiplication, in the sense that we
define a product for linear functionals in terms of Wall
functions. Then we can show that our definition of
the  Haar functional is indeed left- and
right invariant with respect to this weak version of the
comultiplication.  

Let us now turn to the contents of this paper. In \S 2 we introduce
the quantum $SU(1,1)$ group and a corresponding 
${\text C}^\ast$-algebra, which can be seen as the algebra of
continuous functions vanishing at infinity on the quantum $SU(1,1)$
group.   
We work here with a faithful representation of the 
${\text C}^\ast$-algebra, and
we can introduce a weight, i.e. an unbounded functional, that is
left and right invariant in the weak sense. This analogue of the
Haar functional is an integral of weighted traces in irreducible
representations of the ${\text C}^\ast$-algebra. 
In \S 3 we recall some
facts on the algebraic level, both for the quantised algebra of
polynomials on $SU(1,1)$ and for the quantised universal enveloping
algebra. This part is mainly intended for notational purposes and
for stating properties that are needed in the sequel. In \S 4 we 
prove the statement for the Haar functional in the second line of
Table~1.2. This is done by a spectral analysis of a three-term
recurrence operator in $\H$ previously studied in \cite{\CiccKK}. 
In \S 5  we prove the statement for the Haar functional in the 
third line of Table~1.2. A spectral analysis of a five-term
recurrence operator in $\H$ 
is needed, and we can do it by factorising it as
the product of two three-term recurrence operators. The
factorisation is motivated by factorisation results on the quantum
group level. At a certain point,
Lemma~5.5, %harde referentie
we require a highly
non-trivial summation formula for basic hypergeometric series, and
the derivation by Mizan Rahman is given in Appendix B. 
The result is an
Askey-Wilson type measure with absolutely continuous part supported
on $[-1,1]$ plus an infinite set of discrete mass points tending to
infinity. In \S 6 we mainly study the spherical Fourier transforms
on the quantum $SU(1,1)$ group. In this section we have to take a
number of formal steps. We show that the radial part of the Casimir
operator corresponds to a 2-parameter Askey-Wilson difference
operator and we calculate the spherical functions in terms of
very-well-poised ${}_8\vp_7$-series. By the results of
\cite{\KoelSAW} we can invert the spherical Fourier transform and we
see that the Plancherel measure is supported on the principal
unitary series representations and an infinite discrete subset of
the strange series representations.  
Finally, Appendix A contains the spectral analysis of a three-term
operator on $\H$ extending the results of Kakehi \cite{\Kake} 
and Appendix B, by Mizan Rahman, contains a number of summation
formulas needed in  the paper. 

\demo{Notation} We use $\N=\{1,2,\ldots\}$, $\Zp =\{0,1,\ldots\}$,
and $q$ is a fixed number with $0<q<1$. For basic hypergeometric
series we use Gasper and Rahman \cite{\GaspR}. So for $k\in \Zp\cup
\{\infty\}$ we use the notation 
$(a;q)_k=\prod_{i=0}^{k-1}(1-aq^i)$ for $q$-shifted factorials, and
also $(a_1,\ldots,a_r;q)_k=\prod_{i=1}^r (a_i;q)_k$. The basic
hypergeometric series is defined by 
$$
{}_r\vp_s\left( {{a_1,\ldots,a_r}\atop{b_1,\ldots,b_s}};q,z\right) 
= \sum_{k=0}^\infty {{(a_1,\ldots,a_r;q)_k
z^k}\over{(q,b_1,\ldots,b_s;q)_k}} \Bigl( (-1)^k q^{\hf k(k-1)}
\Bigr)^{s+1-r},
$$
whenever it is well-defined. The series is balanced if $r=s+1$, 
$b_1\ldots b_s=qa_1\ldots a_{s+1}$ and $z=q$. The series is called
very-well-poised if $r=s+1$, $qa_1=a_2b_1=a_3b_2=\ldots =a_{s+1}b_s$
and $a_2=q\sqrt{a_1}$, $a_3=-q\sqrt{a_1}$. For the very-well-poised
series we use the notation
$$
\align
{}_{s+1}W_s(a_1;a_4,\ldots,a_{s+1};q,z) &= 
{}_{s+1}\vp_s\left( {{a_1,qa_1^\hf,-qa_1^\hf,a_4,\ldots,a_{s+1}}
\atop{a_1^\hf,-a_1^\hf,qa_1/a_4,\ldots,qa_1/a_{s+1}}};q,z\right) \\
&= \sum_{k=0}^\infty {{1-a_1q^{2k}}\over{1-a_1}}
{{(a_1,a_4,\ldots, a_{s+1};q)_k z^k}\over{(q,
qa_1/a_4,\ldots,qa_1/a_{s+1};q)_k}}.
\endalign
$$
\enddemo

\demo{Acknowledgement} We are grateful to Mizan Rahman for his 
crucial help, which resulted in Appendix B. We thank Tom Koornwinder
for giving us access to unpublished notes, that have served as the
basis for Appendix A, and for creating the stimulating
environment at the Universiteit van Amsterdam, where both authors
were based when this project began. 
The second author is supported by a NWO-Talent stipendium of the 
Netherlands Organization for Scientific Research (NWO).
Part of the research was done while the second author
was supported by the EC TMR network ``Algebraic Lie
Representations'', grant no. ERB FMRX-CT97-0100. 
We thank Alfons Van Daele,
Johan Kustermans, Erik Thomas and Hjalmar Rosengren 
for useful discussions. 
\enddemo 

%%%%%%%%%%%%%%%%%%%%%%%%%%%%%%%%%%%%%%%%%%%%%%%%%%%%%%%%%%%%%%%%%%%%
%N E W   S E C T I O N%
%%%%%%%%%%%%%%%%%%%%%%%%%%%%%%%%%%%%%%%%%%%%%%%%%%%%%%%%%%%%%%%%%%%%
\head \newsec The quantum $SU(1,1)$ group\endhead

The quantum $SU(1,1)$ group is introduced as a Hopf $\ast$-algebra.
In \S \the\sectionnumber.1 we describe its irreducible
$\ast$-representations in terms
of unbounded operators. In \S \the\sectionnumber.2 we use these
representations to define a $\text{C}^\ast$-algebra,
which we regard as the algebra of
continuous functions on the quantum $SU(1,1)$ group 
which tend to zero 
at infinity, and we define
a Haar weight on the $\text{C}^\ast$-algebra. The Haar weight and
the $\text{C}^\ast$-algebra are the same as previously introduced
for the quantum group of plane motions by Woronowicz
\cite{\WoroEtwee} and also studied by Baaj \cite{\Baaj},
Quaegebeur and Verding \cite{\QuaeV}, Verding \cite{\Verd}.
We define a new product on certain linear functionals
in terms of Wall functions, which reflect the comultiplication
of the quantum $SU(1,1)$ group. For this product we show that
the Haar weight is right and left invariant.

%%%%%%%%%%%%%%%%%%%%%%%%%%%%%%%%%%%%%%%%%%%%%%%%%%%%%%%%%%%%%%%%%%%%
\subhead \the\sectionnumber.1\ Representations
of $A_q(SU(1,1))$\endsubhead
We first recall some generalities on the quantum $SL(2,\C)$
group and a non-compact real form, the quantum $SU(1,1)$ group,
see e.g. Chari and Pressley \cite{\CharP} or any other
textbook on quantum groups.
Let $A_q(SL(2,\C))$ be the unital algebra
over $\C$ generated by $\al$, $\be$,
$\ga$ and $\de$ satisfying
$$
\gathered
\al\be =q\be\al ,\quad \al\ga = q\ga\al ,\quad \be\de = q\de\be ,
\quad \ga\de = q\de\ga ,\\
\be\ga =\ga\be ,\quad \al\de -q\be\ga = \de\al - q^{-1}\be\ga =1,
\endgathered
\tag\eq{201}
$$
where $1$ denotes the unit of $A_q(SL(2,\C))$ and $0<q<1$.
A linear basis for this algebra is given by
$\{ \al^k\be^l\ga^m\mid k,l,m\in\Zp\}\cup
\{\de^k\be^l\ga^m\mid k\in\N,l,m\in\Zp\}$.
This is a Hopf-algebra with comultiplication $\De\colon
A_q(SL(2,\C))\to A_q(SL(2,\C))\otimes A_q(SL(2,\C))$,
which is an algebra homomorphism, given by
$$
\gather
\De(\al) = \al\otimes\al + \be\otimes\ga,  \qquad
\De(\be) = \al\otimes\be +\be\otimes\de, \\
\De(\ga) = \ga\otimes\al + \de\otimes\ga, \qquad
\De(\de) = \de\otimes\de +\ga\otimes\be, 
\endgather
$$
and counit $\ep\colon A_q(SL(2,\C))\to\C$, which is
an algebra homomorphism, given by
$\ep(\al)=\ep(\de)=1$, $\ep(\be)=\ep(\ga)=0$. There is
also an antipode $S\colon A_q(SL(2,\C))\to A_q(SL(2,\C))$,
which is an antimultiplicative
linear mapping given on the
generators by $S(\al)=\de$, $S(\be)=-q^{-1}\be$,
$S(\ga)=-q\ga$ and $S(\de)=\al$.
We say that
a linear functional $h\colon A_q(SL(2,\C))\to\C$ is
right invariant, respectively left invariant,
if $(h\otimes \id)\De(a) = h(a)1$, respectively
$(\id\otimes h)\De(a) = h(a)1$, in
$A_q(SL(2,\C))$. So $h$ is a right invariant Haar
functional if and only if $h\star \om=\om(1)h$ for any linear
functional $\om\colon A_q(SL(2,\C))\to\C$, where the product
of two linear functionals $\om, \om'$ is defined by
$\om\star \om'=(\om\otimes\om')\circ\De$.

With $A_q(SU(1,1))$ we denote the $\ast$-algebra
which is $A_q(SL(2,\C))$ as an algebra with $\ast$ given by
$\al^\ast=\de$, $\be^\ast=q\ga$,
$\ga^\ast=q^{-1}\be$, $\de^\ast=\al$. So $A_q(SU(1,1))$ is the
$\ast$-algebra generated by $\al$ and $\ga$ subject to
the relations
$$
\al\ga=q\ga\al,\quad \al\ga^\ast=q\ga^\ast\al, \quad
\ga\ga^\ast=\ga^\ast\ga, \quad
\al\al^\ast-q^2\ga^\ast\ga=1=\al^\ast\al-\ga\ga^\ast.
\tag\eq{205}
$$
This is in fact a Hopf $\ast$-algebra, implying that
$\De$ and $\ep$ are $\ast$-homomorphisms and
$S\circ\ast$ is an involution. In particular,
$$
\De(\al) =\al\otimes\al + q\ga^\ast\otimes\ga,\qquad
\De(\ga) = \ga\otimes\al + \al^\ast\otimes \ga.
\tag\eq{210}
$$

We can represent the $\ast$-algebra $A_q(SU(1,1))$ by unbounded
operators in the Hilbert space $\H$
with standard orthonormal basis $\{e_k\mid k\in \Z\}$.
Since the representation involves unbounded operators
we have to be cautious. We stick to the conventions
of Schm\"udgen \cite{\Schm, Def. 8.1.9}: given a dense linear
subspace ${\Cal D}$ of a Hilbert space ${\Cal H}$, a mapping
$\pi$ of a unital $\ast$-algebra $A$ 
into the set of linear operators
defined on ${\Cal D}$ is a {\it $\ast$-representation of $A$}
if
\roster
\item"{(i)}" $\pi(c_1a_1+c_2a_2)v=c_1\pi(a_1)v+c_2\pi(a_2)v$
and $\pi(1)v=v$
for all $a_i\in A$, $c_i\in\C$ ($i=1,2$) and all $v\in{\Cal D}$,
\item"{(ii)}" $\pi(b){\Cal D}\subseteq {\Cal D}$ and
$\pi(ab)v=\pi(a)\pi(b)v$ for all $a,b\in A$ and all $v\in{\Cal D}$,
\item"{(iii)}" $\langle\pi(a)v,w\rangle =
\langle v,\pi(a^\ast)w\rangle$
for all $a\in A$ and all $v,w\in{\Cal D}$.
\endroster
Note that (iii) states that the domain of the adjoint $\pi(a)^\ast$
contains ${\Cal D}$ for any $a\in A$ and 
$\pi(a)^\ast\vert_{\Cal D}=\pi(a^\ast)$, so that $\pi(a)$ is
closeable. 
See also Woronowicz
\cite{\WoroEtwee, \S 4}.

By ${\Cal D}(\Z)$ we denote the dense subspace of $\H$
consisting of finite linear combinations of the standard
basis vectors $e_k$, $k\in\Z$. 

\proclaim{Proposition \thname{220}}
{\rm (i)}
Let $\lambda\in \C\backslash\{0\}$. There exists a unique
$\ast$-representation $\pi_{\lambda}$ of $A_q(SU(1,1))$ acting on
$\H$ with common domain
${\Cal D}(\Z)$,
such that
$$
\gather
\pi_\la(\al)\, e_k = \sqrt{1+|\la|^2q^{-2k}} \, e_{k+1}, \quad
\pi_\la(\ga)\, e_k = \la q^{-k}\, e_k, \\
\pi_\la(\al^\ast)\, e_k = \sqrt{1+|\la|^2q^{2-2k}} \, e_{k-1},
\quad
\pi_\la(\ga^\ast)\, e_k = \bar \la q^{-k}\, e_k.
\endgather
$$
\par
\noindent
{\rm (ii)} The $\ast$-representation $\pi_\la$ is irreducible and
for $\la,\mu\in R=
\{ z\in \C\mid q<|z|\leq 1\}$ the $\ast$-representations
$\pi_\la$ and $\pi_\mu$ are inequivalent for $\la\not=\mu$. This
means that the space of intertwiners
$I_{\la,\mu} = \{ T\in{\Cal B}(\H)\mid T({\Cal D}(\Z))\subseteq
{\Cal D}(\Z),\,
T\pi_\mu(a)v=\pi_\la(a)Tv, \, \forall a\in A_q(SU(1,1)),
\forall v\in {\Cal D}(\Z)\}$ equals $\{0\}$ for $\la\not=\mu$,
$\la,\mu\in R$ and equals $\C\cdot 1$ for $\la=\mu$.
\par\noindent
\endproclaim

\demo{Remark} These are precisely the representations described
by Woronowicz \cite{\WoroEtwee, \S 4}.
\enddemo

\demo{Proof} To prove (i) we observe that
$\pi_\la(a)$ preserves ${\Cal D}(\Z)$ for
$a\in\{\al,\al^\ast,\ga,\ga^\ast\}$ so that we have
compositions of these operators. It is a straightforward
calculation to see that 
these operators satisfy the commutation relations
\eqtag{205}. It follows 
that $\pi_{\la}$ uniquely extends to an algebra homomorphism
$\pi_{\la}$ of $A_q(SU(1,1))$ into the algebra of linear operators
on ${\Cal D}(\Z)$. It remains to prove
that $\langle \pi_{\la}(a)v,w\rangle=\langle v,
\pi_{\la}(a^\ast)w\rangle$
for all $a\in A_q(SU(1,1))$ and $v,w\in{\Cal D}(\Z)$, 
which follows from checking it for
the generators $a=\al$ and $a=\ga$.

For (ii) we fix an intertwiner $T\in I_{\la,\mu}$. Then
$$
\la q^{-k}Te_k=T\bigl(\pi_{\la}(\ga)e_k\bigr)
=\pi_{\mu}(\ga)\bigl(Te_k\bigr),\qquad k\in \Z,
\tag\eq{221}
$$
so $Te_k=0$ for all $k\in \Z$ if $\la\mu^{-1}\not\in q^{\Z}$.
If $\la,\mu\in R$, then $\la\mu^{-1}\not\in q^{\Z} \Leftrightarrow
\la\not=\mu$. Hence,
$I_{\la,\mu}=\{0\}$ for $\la,\mu\in R$ with $\la\not=\mu$.

If $\la=\mu$, then it follows from \eqtag{221} that $Te_k=c_ke_k$
for some $c_k\in \C$. Since $T$ commutes with
$\pi_{\la}(\al)$, it follows that $c_k$ is independent of
$k$, proving that $I_{\la,\la}=\C \cdot 1$. 
\qed\enddemo

\demo{Remark} There is a canonical way to associate an
adjoint representation $\pi_\la^\ast$ to the representation
$\pi_\la$, see \cite{\Schm, p.~202}, by defining its
common domain as the intersection of the domains of all adjoints,
${\Cal D}^\ast = \cap_{a\in A_q(SU(1,1))} 
{\Cal D}(\pi_\la(a)^\ast)$,
and the action by $\pi_\la^\ast(a) =
\pi_\la(a^\ast)^\ast|_{{\Cal D}^\ast}$.
For the domain ${\Cal D}(\Z)$ we do not have 
self-adjointness of the
representation $\pi_\la$. However, if we replace the common
domain of $\pi_{\la}$  by
$$
{\Cal S}(\Z) = \{ \sum_{k=-\infty}^\infty c_k\, e_k\in\H\mid
\sum_{k=-\infty}^\infty q^{-2nk} |c_k|^2<\infty, \, \forall\,
n\in \Zp\}
$$
we do have $\pi_\la^\ast = \pi_\la$, i.e.
the domains and operators are all the same.
Indeed, observe that
${\Cal D}(\pi_\la(\ga^n)^\ast) = \{
\sum_{k=-\infty}^\infty c_k\, e_k\in\H\mid
\sum_{k=-\infty}^\infty q^{-nk} c_k\, e_k \in\H\}$
and hence ${\Cal D}^\ast\subseteq {\Cal S}(\Z)$. By
\cite{\Schm, Prop.~8.1.2} this implies that
the representation by unbounded operators is self-adjoint.
\enddemo

Having the irreducible $\ast$-representations of
Proposition \thtag{220} we can form the direct integral
$\ast$-representation
$\pi=(2\pi)^{-1}\int_0^{2\pi}\pf\, d\phi$, see
\cite{\Schm, Def.~12.3.1},  with its
representation space $(2\pi)^{-1}\int_0^{2\pi} \H d\phi
\cong L^2(\T;\H)$ equipped with the orthonormal
basis $e^{ix\phi}\otimes e_m$ for $x,m\in\Z$.
The common domain is by definition 
$$
\multline
{\Cal D}(L^2(\T;\H))  = \{ f\in L^2(\T;\H) \mid
f(e^{i\phi})\in{\Cal D}(\Z)\, \text{a.e. and}\\
e^{i\phi}\mapsto \pf(a)f(e^{i\phi})\in L^2(\T;\H)\, \forall
a\in A_q(SU(1,1))\}.
\endmultline
\tag\eq{231}
$$
The last condition means in particular that 
$(2\pi)^{-1}\int_0^{2\pi}
\| \pf(a)f(e^{i\phi})\|^2d\phi<\infty$ for all $a\in A_q(SU(1,1))$.
In this case ${\Cal D}(L^2(\T;\H))$ is dense in $L^2(\T;\H)$
since it contains finite linear combinations of the
basis elements $e^{ix\phi}\otimes e_m$.
The action of the generators of $A_q(SU(1,1))$
on the basis of $L^2(\T;\H)$
can be calculated explicitly
from Proposition \thtag{220};
$$
\gather
\pi(\ga)\, e^{ix\phi}\otimes e_m = q^{-m}\,
e^{i(x+1)\phi}\otimes e_m, \qquad
\pi(\al)\, e^{ix\phi}\otimes e_m = \sqrt{1+q^{-2m}}\,
e^{ix\phi}\otimes e_{m+1}, \\
\pi(\ga^\ast)\, e^{ix\phi}\otimes e_m = q^{-m}\,
e^{i(x-1)\phi}\otimes e_m, \qquad
\pi(\al^\ast)\, e^{ix\phi}\otimes e_m = \sqrt{1+q^{2-2m}}\,
e^{ix\phi}\otimes e_{m-1}.
\endgather
$$

\proclaim{Lemma \thname{235}} The direct integral representation
$\pi=(2\pi)^{-1}\int_0^{2\pi}\pf\, d\phi$ is a
faithful representation of the $\ast$-algebra
$A_q(SU(1,1))$, i.e. $\pi(\xi)f=0$ for
all $f\in{\Cal D}(L^2(\T;\H))$ implies $\xi=0$
in $A_q(SU(1,1))$.
\endproclaim

\demo{Proof}
The action of the monomial basis of $A_q(SU(1,1))$ under the
representation $\pi$ is given by
$$
\gather
\pi\bigl(\al^r(\gamma^\ast)^s\ga^t\bigr)e^{ip\th}\otimes e_l=\,
q^{-l(s+t)}(-q^{-2l};q^{-2})_r^\hf e^{i(p+t-s)\th}\otimes
e_{l+r},\\
\pi\bigl((\al^\ast)^r(\ga^\ast)^s\ga^t\bigr)e^{ip\th}\otimes e_l=\,
q^{-l(s+t)\th}(-q^{2-2l};q^2)_r^\hf e^{i(p+t-s)\th}\otimes
e_{l-r}.
\endgather
$$
It easily follows that $\pi$ is a faithful
representation of $A_q(SU(1,1))$. 
\qed\enddemo

In the next subsection we give a 
coordinate-free realisation of the
representation $\pi$.

%%%%%%%%%%%%%%%%%%%%%%%%%%%%%%%%%%%%%%%%%%%%%%%%%%%%%%%%%%%%%%%%%%%%
\subhead \the\sectionnumber.2\ The Haar functional\endsubhead
Let $L^2(X,\mu)$ be the Hilbert space of 
square integrable functions
on $X= \T\times q^\Z \cup \{0\}$ with respect to the
measure
$$
\int f \, d\mu = \sum_{k=-\infty}^\infty
{1\over{2\pi}} \int_0^{2\pi} f(q^k e^{i\th})\, d\th.
$$
Then the map $\psi\colon L^2(\T;\H)\to L^2(X,\mu)$ given by
$$
\psi\colon e^{ix\phi}\otimes e_m \mapsto 
\Bigl( f_{x,m}\colon z\mapsto
\de_{|z|, q^{-m}} \bigl( {z\over{|z|}}\bigr)^x\Bigr)
\tag\eq{249}
$$
is a unitary isomorphism. Observe that
$\psi\pi(\ga)\psi^{-1}=M_z$, $\psi\pi(\ga^\ast)
\psi^{-1}=M_{\bar z}$,
$\psi\pi(\al)\psi^{-1}=M_{\sqrt{1+q^2|z|^2}}T_q$ and
$\psi\pi(\al^\ast)\psi^{-1}=M_{\sqrt{1+|z|^2}}T_q^{-1}$, where
$M_g$ denotes the operator of multiplication by the 
function $g$ and
$T_q$ is the $q$-shift operator defined by $(T_qf)(z)=f(qz)$.

These formulas for the action of the generators of $A_q(SU(1,1))$
under the faithful representation $\psi\pi(\cdot)\psi^{-1}$
suggest the following formal definition for the 
$\text{C}^\ast$-algebra
of continuous functions on the quantum $SU(1,1)$ group 
which vanish at
infinity: it is the $C^\ast$-subalgebra of ${\Cal B}(L^2(X,\mu))$
generated by $M_g$, $g\in C_0(X)$, and the $q$-shifts $T_{q^{\pm
1}}$, where $C_0(X)$ is the $C^\ast$-algebra consisting of
continuous functions on $X$ which vanish at infinity 
(here $X$ inherits its topology from $\C$, and 
the $C^\ast$-norm is given by the supremum norm
$\|\cdot \|_{\infty}$). In
other words, one replaces the unbounded action of the subalgebra
$\C[\ga, \ga^\ast]\subset A_q(SU(1,1))$ on $L^2(X,\mu)$
by the bounded, regular action of $C_0(X)$ on $L^2(X,\mu)$.

To make the construction rigorous we use the notion of a crossed
product $\text{C}^\ast$-algebra, see \cite{\Pede, Ch.~7}. 
The crossed product needed here is the same as for the
quantum group of plane motions, see Baaj \cite{\Baaj},
Woronowicz \cite{\WoroEtwee}. Let us recall the construction
in this specific case. For $k\in\Z$ we define
the automorphisms $\tau_k$ of $C_0(X)$ by
$\tau_k(f)=(T_q)^kf$. Then $(C_0(X),\Z,\tau)$
is a $\text{C}^\ast$-dynamical system, see \cite{\Pede, \S 7.4}.

Let $\ell^1(\Z;C_0(X))$ be the $\ell^1$-functions 
$f\colon \Z\rightarrow C_0(X)$
with respect to the norm 
$\|f\|_1=\sum_{n\in\Z}\|f_n\|_{\infty}$.
The subspace $C_c(\Z;C_0(X))=
\{ f\colon \Z\rightarrow C_0(X) \, | \,
\#\hbox{supp}(f)<\infty\}$ is dense in $\ell^1(\Z;C_0(X))$.
Furthermore, $\ell^1(\Z;C_0(X))$ is a Banach $\ast$-algebra, with
$\ast$-structure and multiplication given by
$$
(f^\ast)_n = \tau_n( (f_{-n})^\ast), \quad
(fg)_n=\sum_{k=-\infty}^\infty f_k \, \tau_k(g_{n-k}).
\tag\eq{240}
$$
The crossed $C^\ast$-product $C_0(X)\times_\tau\Z$ is by
definition the strong closure of $\ell^1(\Z;C_0(X))$
under its universal representation, where the universal
representation is the direct sum of the non-degenerate
$\ast$-representations of $\ell^1(\Z;C_0(X))$, 
see \cite{\Pede, \S 7.6}.
We regard $C_0(X)\times_{\tau}\Z$ as the quantum analogue of the
$C^\ast$-algebra consisting of continuous functions on $SU(1,1)$
which vanish at infinity.

We interpret the faithful $\ast$-representation
$\pi$ of \S 2.1 as a non-degenerate
representation of $C_0(X)\times_{\tau}\Z$ in the following way.
Let $\tilde{\pi}$ be the regular representation of $C_0(X)$ on
$L^2(X,\mu)$, and let $u\colon \Z\rightarrow {\Cal B}(L^2(X,\mu))$
be the unitary representation defined by $u_nf=(T_q)^nf$. 
Then $(\tilde{\pi},u,L^2(X,\mu))$ is a covariant
representation of the $C^\ast$-dynamical system
$(C_0(X),\Z,\tau)$,
i.e. $\tilde{\pi}(\tau_nf)=u_n\tilde{\pi}(f)u_n^\ast$ 
for all $f\in C_0(X)$ and $n\in \Z$.
In the notation of \cite{\Pede}, we get a non-degenerate
representation $\pi=\tilde{\pi}\times u$ of $C_0(X)\times_\tau\Z$
on $L^2(X,\mu)$, which is defined on $C_c(\Z;C_0(X))$ by
$\pi(f)g=\sum_{n\in \Z}\tilde{\pi}(f_n)(u_ng)$. 
More explicitly, we have
$$
(\pi(f)g)(z)=\sum_{n\in\Z}f_n(z)g(q^nz),\qquad f\in
C_c(\Z;C_0(X)),\,\,\, g\in L^2(X,\mu),\,\,\, x\in X.
$$
The $\ast$-representations $\pt$ of \S 2.1 can also be considered
as a covariant representation of $C_0(X)\times_{\tau}\Z$.
Let the representation 
$\tilde\pt\colon C_0(X)\to {\Cal B}(\H)$ of
the commutative $\text{C}^\ast$-algebra be defined
by $\tilde\pt(f)e_l=f(q^{-l}e^{i\th})e_l$, and let
the unitary representation $u\colon \Z \to {\Cal B}(\H)$
be defined by $u_n =U^n$, where $U\colon\H\to\H$, 
$e_k\mapsto e_{k+1}$ is the shift operator. Then this 
gives a covariant
representation. The direct integral representation 
$(2\pi)^{-1}\int_0^{2\pi} \pt\, d\th$ in
$(2\pi)^{-1}\int_0^{2\pi} \H\, d\th \cong L^2(\T;\H)$ 
is equivalent to $\pi$ using $\psi$ as in \eqtag{249}.

\demo{Remark}
The matrix elements of $\pi$ with respect to the
orthonormal basis $f_{k,l}$, $k,l\in \Z$, 
of $L^2(X,\mu)$, see \eqtag{249}, give the linear functionals
$$
\om_{s,l}^{r-k}(f)=
\langle \pi(f)f_{k,l}, f_{r,s}\rangle_{L^2(X,\mu)}
=\frac{1}{2\pi}\int_0^{2\pi}f_{s-l}(q^{-s}e^{i\th})
e^{i(k-r)\th}d\th
$$
for $f\in C_c(\Z;C_0(X))$, which can be uniquely extended by
continuity to a linear functional on $C_0(X)\times_\tau \Z$.
Note that
$$
\om_{k,l}^x(fg) = \sum_{r,y\in\Z} \om_{k,r}^{x-y}(f)
\om_{r,l}^y(g).
$$
\enddemo

\proclaim{Proposition \thname{259}}
$\pi$ is a faithful representation of $C_0(X)\times_\tau\Z$.
\endproclaim

\demo{Proof}
Recall that $\tilde{\pi}$ is the regular 
representation of $C_0(X)$
on $L^2(X,\mu)$. The regular representation 
of $C_0(X)\times_\tau\Z$
induced by $\tilde{\pi}$, which we denote 
here by $\rho$, acts on
$\ell^2(\Z;L^2(X,\mu))$ by
$$
(\rho(f)g)_n(z)=\sum_{m\in\Z}\bigl(\tilde{\pi}
(\tau_{-n}(f_m))g_{n-m}\bigr)(z)=
\sum_{m\in\Z}f_m(q^{-n}z)g_{n-m}(z)
$$
for $f\in C_c(\Z;C_0(X))$, $g\in 
\ell^2(\Z;L^2(X,\mu))$ and $z\in
X$,
see \cite{\Pede, Ch.~7}. Explicitly, 
the action of $\rho$ in terms of
the orthonormal basis
$g_{k,l,m}=\delta_{k}(\cdot)f_{l,m}$, $k,l,m\in \Z$,  
of $\ell^2(\Z;L^2(X,\mu))$
is given by
$$
\rho(f)g_{k,l,m}=\sum_{r,s}\left(\frac{1}{2\pi}
\int_0^{2\pi}f_{r-k}(q^{-m-r}e^{i\th})
e^{i(l-s)\th}d\th\right) g_{r,s,m}.
$$
So the closure $H_m$ of $\text{span}\{ g_{k,l,m}\mid
k,l\in\Z\}$ is an invariant subspace for $\rho$ and
we have the orthogonal direct sum decomposition
$\ell^2(\Z;L^2(X,\mu))=\oplus_{m\in\Z} H_m$.
{}From the explicit formulas for $\pi$ and $\rho$ with respect to
the orthonormal basis $f_{k,l}$, $k,l\in\Z$, of $L^2(X,\mu)$,
respectively $g_{k,l,m}$, $k,l,m\in\Z$, of
$\ell^2(\Z;L^2(X,\mu))$,
it follows that
$L^2(X,\mu)\to H_m$, $f_{l,k}\mapsto
g_{k-m,l,m}$ is a unitary intertwiner between $\pi$ and
$\rho\vert_{H_m}$, so that $\rho\simeq \bigoplus_{m\in \Z}\pi$
as representations of $C_0(X)\times_\tau\Z$. By \cite{\Pede,
Cor.~7.7.8} we know that $\rho$ is a faithful representation of
$C_0(X)\times_\tau\Z$, hence so is $\pi$.
\qed\enddemo

Recall that a {\it weight} on a $\text{C}^\ast$-algebra $A$ is
a function $h\colon A^+\to [0,\infty]$ satisfying
(i) $h(\la a)=\la h(a)$ for $\la\geq 0$ and $a\in A^+$,
(ii) $h(a+b)=h(a)+h(b)$ for $a,b\in A^+$.
The weight $h$ is
said to be densely defined if
$\{ a\in A^+ \, | \, h(a)<\infty \}$ is dense in $A^+$.
Furthermore, we say that
$h$ is lower semi-continuous if
$\{ a\in A^+\mid h(a)\leq \la\}$ is closed for any $\la\geq 0$, and
that $h$ is faithful if $h(a^\ast a)=0$ implies $a=0$ in $A$.
A weight can be extended uniquely to $N_h^\ast N_h =
\{ a^\ast b\mid a,b\in N_h\}$, where $N_h=
\{ a\in A\mid h(a^\ast a)<\infty\}$. 
See Combes \cite{\Comb, \S 1}, Pedersen \cite{\Pede, Ch.~5}
for general information and
for application of weights in quantum groups
see Kustermans and Vaes \cite{\KustV},
Quaegebeur and Verding \cite{\QuaeV}, Verding \cite{\Verd}.

Let $h$ be a lower semi-continuous, densely defined weight on
$A$. 
The GNS-construction for weights gives a 
Hilbert space $H_h$ and a representation $\si_h$ of
$A$ in $H_h$ and a
linear map $\Lambda_h$ from  $N_h = \{ f\in A
\mid h(f^\ast f)<\infty\}$ onto a dense subspace
of  $H_h$ satisfying
\roster
\item"{(i)}" $h(f^\ast g) = \langle
\Lambda_h(g),\Lambda_h(f)\rangle$
for all $f,g\in N_h$,
\item"{(ii)}" $\si_h(f)\Lambda_h(g) = \Lambda_h(fg)$
for all $f\in A$ and all $g\in N_h$,
\item"{(iii)}" the representation $\si_h$ is non-degenerate, 
i.e. the closure of the linear span of elements of
the form $\si_h(f)g$ for $f\in A$ and $g\in H_h$ equals
$H_h$,
\item"{(iv)}" the map $\Lambda_h\colon N_h\to H_h$
is closed.
\endroster
Properties (i) and (ii) hold by the general GNS-construction
of weights, see \cite{\Comb, \S 2}, and properties (iii) and
(iv) follow since $h$ is lower semi-continuous, see
\cite{\Comb, \S 2}, \cite{\Verd, Prop.~2.1.11}.
If $h$ is faithful, we obtain $H_h$ as
the Hilbert space completion of $N_h$ with respect to
inner product $\langle f,g\rangle =h(g^\ast f)$. 

The following theorem has been proved by Baaj \cite{\Baaj, \S 4} 
in the setting of weights on von Neumann algebras and later
in the $\text{C}^\ast$-algebra framework by
Quaegebeur and Verding \cite{\QuaeV, \S 4}, 
Verding \cite{\Verd, \S 3.2}. 

\proclaim{Theorem \thname{260}} Let $h=\sum_{k=-\infty}^\infty
q^{-2k}\om^0_{k,k}$, then $h$ is a densely defined,
faithful,
lower semi-continuous weight on $C_0(X)\times_\tau \Z$.
\endproclaim

\demo{Remark \thname{261}} Note that $C_c(\Z;C_c(X))\subset
N_h$, since for $f\in C_c(\Z;C_c(X))$ we have
$f^\ast f\in C_c(\Z;C_c(X))$ so that $\om^0_{k,k}(f^\ast f)=0$
for $k$ sufficiently large. Also, the dense subspace
$C_c(\Z;C_c(X))$ of $C_0(X)\times_\tau \Z$ is contained in
$N_h^\ast N_h$, so that the Haar functional is well-defined
on $C_c(\Z;C_c(X))$. Indeed, take $g^{(m)}\in C_c(\Z;C_c(X))$, 
$g^{(m)}_n(x)=\de_{n,0} u^{(m)}(x)$, where
$u^{(m)}$ is a compactly supported approximate unit in
$C_0(X)$ with support in $|x|\leq q^{-m}$. Then $(g^{(m)})^\ast
\in N_h$, and $g^{(m)} f= f$ for $f\in C_c(\Z;C_c(X))$
and $m$ sufficiently large.
\enddemo

\demo{Remark}
We regard $H_h$ in this paper as the $q$-analogue of the
$L^2$-functions on $SU(1,1)$ with respect to the Haar measure.
It was shown by Baaj \cite{\Baaj}
and Quaegebeur and Verding \cite{\QuaeV} that $H_h$ is isomorphic
to $\ell^2(\Z^3)$. This corresponds nicely with the fact that 
$SU(1,1)$ is a three-dimensional Lie group.
\enddemo 

%%%%%%%%%%%%%%%%%%%%%%%%%%%%%%%%%%%%%%%%%%%%%%%%%%%%%%%%%%%%%%%%%%%%
\subhead \the\sectionnumber.3\ Invariance of the 
Haar functional\endsubhead
The weight $h$ is left and right invariant when considered
as a weight corresponding to the quantum group of plane motions,
see Baaj \cite{\Baaj, Thm.~4.2}. For this we have to have
the comultiplication of the quantum group of plane motions
on the $\text{C}^\ast$-algebra level, and this has been done
by Woronowicz \cite{\WoroEtwee}. This seems not to be possible
for the quantum $SU(1,1)$ group, see Woronowicz \cite{\WoroEtwee,
Thm. 4.1}. However we can introduce
the comultiplication for the quantum $SU(1,1)$ group
on the $\text{C}^\ast$-algebra level
in a weak form, and then the Haar weight is also left and right
invariant. For this we first have to encode
the comultiplication as in \eqtag{210} in terms of a
product for the matrix elements. See Baaj \cite{\Baaj}
for a similar procedure for the quantum group of plane motions.

\proclaim{Lemma \thname{270}}
For $x,m,k\in\Z$ define the normalised Wall function by
$$
f_m^x(k) = {{(-1)^k q^{(k-m)(1+x)}
(-q^{2-2k};q^2)_\infty^\hf }\over{
(-q^{2-2m},-q^{2-2k+2x};q^2)_\infty^\hf}} {{(q^{2+2x};q^2)_\infty}
\over{(q^2;q^2)_\infty}} \, {}_1\vp_1\left(
{{-q^{2+2x-2k}}\atop{q^{2+2x}}};q^2, q^{2+2k-2m}\right).
$$
The product of linear functionals $\om^x_{k,l}$ given by
$$
\om^x_{k,l}\star\om^y_{r,s} =
\de_{l-k-y,s-r+x}(-1)^{l-k-y}
\sum_{n=-\infty}^\infty  f^{s-l}_{n+l-k-y}(s) f^{r-k}_n(r)
\, \om^{r-s+l-k}_{n,n+l-k-y}
$$
is well-defined as a linear functional on
$C_0(X)\times_\tau \Z$.
\endproclaim

\demo{Remark \thname{217}} In order to motivate the definition of
Lemma \thtag{270} let us first consider the case of the
compact quantum $SU(2)$ group. The analogue of the
representations $\pt$ are in terms of bounded
operators on $\ell^2(\Zp)$ and the tensor product
decomposition $\pt\otimes \pp\cong (2\pi)^{-1}
\int_0^{2\pi} \pi_{e^{i\phi}}\, d\phi$ holds,
see \cite{\KoorAF}, \cite{\Sheu}. Moreover, the
Clebsch-Gordan coefficients are explicitly given in
terms of Wall polynomials, see \cite{\KoorAF}, and the
Clebsch-Gordan coefficients are determined by a spectral
analysis of the compact operator
$(\pt\otimes \pp)\De(\ga^\ast\ga)$. This is done by interpreting
this operator as a three-term recurrence operator, i.e. as
a Jacobi matrix, corresponding to the Wall polynomials. The
Clebsch-Gordan coefficients then determine the
product $\om\star\om'=(\om\otimes\om')\circ\De$
of the matrix elements  $\om,\om'$, see also
\cite{\Baaj, Prop. 4.3} for the quantum group of plane motions.
For the quantum $SU(1,1)$ group we can formally follow
the same method using the representations
$\pt$ as in \S \the\sectionnumber.1. In this case we have from
Proposition \thtag{220} and \eqtag{210}
$$
\multline
(\pt\otimes\pp)\De(\ga^\ast\ga) \, e_k\otimes e_l =
\bigl(q^{-2k}(1+q^{-2l})+q^{-2l}(1+q^{2-2k})\bigr) e_k\otimes e_l\\
+e^{i(\th-\psi)}q^{-k-l-1} \bigl((1+q^{-2k})(1+q^{-2l})\bigr)^\hf
e_{k+1}\otimes e_{l+1}\\ +
 e^{i(\psi-\th)}q^{-k-l+1} \bigl((1+q^{2-2k})(1+q^{2-2l})\bigr)^\hf
e_{k-1}\otimes e_{l-1}.
\endmultline
$$
Hence, the unbounded operator $(\pt\otimes\pp)\De(\ga^\ast\ga)$
leaves the subspace
$({\Cal D}(\Z)\times{\Cal D}(\Z))\cap H_x$ invariant, where
$H_x$ is the closure of 
$\text{span}\{e_{k-x}\otimes e_k\mid k\in\Z\}$, so that $H_x
\cong\H$ for any $x\in\Z$. Restricting the unbounded
symmetric operator $(\pt\otimes\pp)\De(\ga^\ast\ga)$ to
${\Cal D}(H_x)$ gives an unbounded symmetric three-term recurrence
operator in $\H$, i.e. a doubly infinite Jacobi matrix, of the
form
$$
\gather
L\, e_k = a_k \, e_{k+1} + b_k\, e_k + a_{k-1}\, e_{k-1}, \qquad
a_k>0, \ b_k\in\R,\\
a_k = q^{(x-1)-2k} \bigl( (1+q^{2(x-k)})(1+q^{-2k})\bigr)^\hf,
\qquad b_k = q^{-2k}(1+q^{2x})+q^{2x-4k}(1+q^2),
\endgather
$$
where $e_k$ is given by $e^{ik(\psi-\th)}e_{k-x}\otimes e_k$. 
This operator fits 
into the general framework as studied in Appendix A. 
Since we have solutions in terms of Wall functions, see
Gupta et al. \cite{\GuptIM, \S 5}, as well as the asymptotically
well-behaved solutions we can work out the details. We find
that the elements 
$$
F^x_m=\sum_{k\in\Z}(-1)^m e^{ix\psi}
e^{i(k-m)(\th-\psi)} f^x_m(k)\,  e_{k-x}\otimes
e_k\in\ell^2(\Z\times\Z)
$$
are eigenvectors of
$(\pt\otimes\pp)\De(\ga^\ast\ga)$ in $H_x$ for the eigenvalue
$q^{-2m}$, $m\in\Z$. The vectors are contained in the domain
of the adjoint $L^\ast$. 
Using contiguous relations we can formally
show that the action of $(\pt\otimes\pp)\De(a)$, 
$a\in A_q(SU(1,1))$,  
on $F^x_m$ is the same as the action of $\pi$ on the orthonormal
basis $\{f_{x,m}\}$ of $L^2(X,\mu)$ by identifying 
$F^x_m$ with $f_{x,m}$, see \eqtag{249}. 
However, $L$ has deficiency indices $(1,1)$ and there is
no self-adjoint extension of $L$ possible such that these
eigenvectors are contained in the domain of the self-adjoint
extension of $L$, so the $\{F^x_m\}$ are not orthogonal, see
\cite{\CiccKK}. 
Note that this observation corresponds to the
no-go theorem of Woronowicz \cite{\WoroEtwee, Thm. 4.1}. 
Ignoring this problem and regarding the
$\ast$-representations 
$(\pt\otimes\pp)\De$ in $\ell^2(\Z\times\Z)$ and
$\pi$ in $L^2(X,\mu)$ equivalent and using the
identity, see \cite{\CiccKK, (4.3)}, 
$$
e_l\otimes e_s =\sum_{m\in\Z} (-1)^m e^{i\psi(l-s)}
e^{i(m-s)(\th-\psi)} f_m^{s-l}(s) \, F^{s-l}_m,
\tag\eq{218}
$$
we formally can rewrite 
$$
(\om^x_{k,l}\otimes \om^y_{r,s})\circ \De = 
{1\over{4\pi^2}} \int_0^{2\pi}\int_0^{2\pi} 
e^{-ix\th}e^{-iy\psi} 
\langle (\pt\otimes\pp)\De(\cdot) e_l\otimes e_s, e_k\otimes e_r
\rangle \, d\th d\psi 
$$
as the linear combination of Lemma \thtag{270} 
using \eqtag{218}. 
So we cannot make this method rigorous, but the results of
\cite{\KoelSbig, \S 2}
suggest that we should consider $L$ on a bigger Hilbert space.
Finally, we note that for the quantum group of plane motions
the interwining operator consisting of the Clebsch-Gordan
coefficients can be used to find a multiplicative unitary, 
see Baaj \cite{\Baaj, \S 4}. For 
the quantum $SU(1,1)$ we only obtain a partial isometry and
we do not expect a multiplicative unitary from this construction.
\enddemo

\demo{Proof of Lemma \thtag{270}} 
Since $|\om_{k,l}^x(f)|\leq \|f\|$ for any
$f\in C_0(X)\times_\tau\Z$, it suffices to show that
$$
\sum_{n=-\infty}^\infty  f^{s-l}_{n+l-k-y}(s) f^{r-k}_n(r)
$$
is absolutely convergent. Since $\{ f_m^x(k)\}_{m\in\Z}
\in \ell^p(\Z)$, $1\leq p\leq \infty$,
see the next lemma and the remark following it,
this follows immediately.
\qed\enddemo

\proclaim{Lemma \thname{275}} For $x,m,k\in\Z$ the Wall functions of
Lemma \thtag{270} satisfy $f^{-x}_m(k)=f^x_m(k+x)$. 
Furthermore, for $x\geq 0$, 
$$
|f^x_m(k)| \leq
{{(-q^{2-2k},-q^{2-2m};q^2)_\infty^\hf (-q^{2+2x};q^2)_\infty}
\over{(-q^{2-2k+2x};q^2)_\infty^\hf (q^2;q^2)_\infty}}
\cases q^{(k-m)(1+x)}, &\text{$m\leq k$,}\\
q^{(m-k)(1+x)}q^{(m-k)(m-k-1)}, &\text{$m\geq k$.}
\endcases
$$
\endproclaim

\demo{Remark} The $\ell^p$-behaviour of
$\{ f^x_m(k)\}_{m=-\infty}^\infty$ follows from Lemma \thtag{275},
since for fixed $x$ and $k$ we see
$f^x_m(k) ={\Cal O}(q^{-m(1+|x|)})$ as $m\to-\infty$ and
$f^x_m(k) ={\Cal O}(q^{\hf m(m-1)}q^{m(1+x-2k)})$ 
as $m\to\infty$ using the theta product identity
$$
(aq^k,q^{1-k}/a;q)_\infty = (-a)^{-k} q^{-\hf k(k-1)}
(a,q/a;q)_\infty, \qquad a\in\C\backslash\{0\},\ k\in\Z.
\tag\eq{237}
$$
\enddemo

\demo{Proof} First observe
$$
\sum_{p=0}^\infty {{(q^{1-n};q)_\infty c_p}\over{(q,q^{1-n};q)_p}}
=\sum_{p=n}^\infty {{(q^{1-n+p};q)_\infty\, c_p}\over{(q;q)_p}}
= \sum_{p=0}^\infty {{(q^{1+n};q)_\infty c_{p+n}}\over{
(q,q^{1+n};q)_p}}
\tag\eq{238}
$$
for $n\in\Z$ provided that the sums are absolutely convergent. 
Applying this with $n=-x$ gives
$f^{-x}_m(k)=f^x_m(k+x)$. 

In order to estimate $f^x_m(k)$, $x\geq 0$, we use
the following limit transition of Heine's formula,
see \cite{\GaspR,(1.4.5)},
$$
(c;q)_\infty \, {}_1\vp_1 ( a;c;q,z) =
(z;q)_\infty \, {}_1\vp_1 ( az/c; z;q,c),
\tag\eq{234}
$$
to rewrite the ${}_1\vp_1$-series in the definition of the
Wall function as follows
$$
\multline
(q^{2+2x};q^2)_\infty\, {}_1\vp_1\left(
{{-q^{2+2x-2k}}\atop{q^{2+2x}}};q^2, q^{2+2k-2m}\right) = \\
(q^{2+2k-2m};q^2)_\infty \, {}_1\vp_1\left(
{{-q^{2-2m}}\atop{q^{2+2k-2m}}};q^2, q^{2+2x}\right) = \\
\sum_{l=0}^\infty (q^{2+2k-2m+2l};q^2)_\infty
{{(-q^{2-2m};q^2)_l}\over{(q^2;q^2)_l}}
(-1)^l q^{l(l-1)} q^{l(2+2x)}.
\endmultline
\tag\eq{276}
$$
Note that for $k-m\geq 0$ the sum starts at $l=0$, but for
$k-m\leq 0$ the sum actually starts at $l=m-k$, cf.
\eqtag{238}. In case $k\geq m$ we estimate the right hand
side of \eqtag{276} termwise to find the bound 
$$
(-q^{2-2m};q^2)_\infty  \sum_{l=0}^\infty {{q^{l(l-1)}}\over
{(q^2;q^2)_l}} q^{l(2+2x)} = (-q^{2-2m},-q^{2+2x};q^2)_\infty
$$
by \cite{\GaspR, (1.3.16)}. Combining this with the definition
of $f^x_m(k)$ gives the desired estimate in this case.

In case $m\geq k$ we rewrite the sum on the right hand side
of \eqtag{276}, by introducing $l=n+m-k$, as, cf. \eqtag{238},
$$
\multline
(-1)^{m-k} q^{(m-k)(m-k-1)} q^{(m-k)(2+2x)} (-q^{2-2m};q^2)_{m-k}
\\ \times \sum_{n=0}^\infty (q^{2+2m-2k+2n};q^2)_\infty
{{(-q^{2-2k};q^2)_n}\over{(q^2;q^2)_n}} (-1)^n q^{n(n-1)}
q^{2n(1+x+m-k)}
\endmultline
$$
and the sum is  estimated
by $(-q^{2-2k},-q^{2+2x+2m-2k};q^2)_\infty \leq
(-q^{2-2k},-q^{2+2x};q^2)_\infty$ in the same way.
This gives the result
for the case $m\geq k$.
\qed\enddemo

In $(C_0(X)\times_\tau\Z)^\ast$ the set $\{ \om^x_{k,l}\}$
is linearly independent. This follows by applying
it to the elements $g_{k,l,m}\in C_c(\Z;C_c(X))$,
see Proposition \thtag{259}. So we have
a well-defined product on linear functionals from
$B\subset (C_0(X)\times_\tau\Z)^\ast$, $B$ being the space of
finite linear combinations of the functionals $\om^x_{k,l}$.
For $\om\in B$ we extend the definition of the product to
$\om\star h$ by requiring that for any $f\in C_c(\Z;C_c(X))$
the expression
$\bigl(\om\star \sum_{k=-N}^N q^{-2k}\om^0_{k,k}\bigr)(f)$
converges as $N\to\infty$. By definition, the resulting
expression is $(\om\star h)(f)$. 
We choose $C_c(\Z;C_c(X))\subset C_0(X)\times_\tau\Z$
because $h$ is defined on $C_c(\Z;C_c(X))$, see
Remark \thtag{261}. A similar definition is
used for $h\star \om$.

\proclaim{Theorem \thname{280}}
Let $\xi\in L^2(X,\mu)$ be a finite linear combination of
the basis elements $f_{x,m}$, $x,m\in\Z$, see \eqtag{249}, and
let $\om_\xi(f)=\langle \pi(f)\xi,\xi\rangle_{L^2(X,\mu)}$
be the corresponding element from $B\subset
(C_0(X)\times_\tau\Z)^\ast$.
Then $h$ is right and left invariant in the sense that
$(\om_\xi\star h)(f) = \|\xi\|^2 h(f) = (h \star \om_\xi)(f)$
for all $f$ in the dense subspace $C_c(\Z;C_c(X))$
of $C_0(X)\times_\tau \Z$.
\endproclaim

The theorem can be extended using the
same argument to $\om_{\xi,\eta}\star h = \langle \xi,\eta\rangle 
h = h \star \om_{\xi,\eta}$ for 
$\om_{\xi,\eta}(f)=\langle \pi(f)\xi,\eta\rangle_{L^2(X,\mu)}$, 
where $\xi,\eta$ are finite linear combinations of 
the basis elements $f_{x,m}$, $x,m\in\Z$.

We start by proving a crucial special case.

\proclaim{Lemma \thname{290}}
$h\star\om^y_{r,s} = \de_{y,0}\de_{r,s}h$
and $\om^x_{k,l}\star h=\de_{x,0}\de_{k,l}h$
on $C_c(\Z;C_c(X))$.
\endproclaim

\demo{Proof} Take $f\in C_c(\Z;C_c(X))$ and consider
$$
\sum_{k=-N}^N q^{-2k} (\om_{k,k}^0\star \om^y_{r,s})(f)
=\de_{y,r-s}(-1)^y \sum_{n=-\infty}^\infty
\bigl( \sum_{k=-N}^N
q^{-2k} f^{s-k}_{n+s}(s) f^{r-k}_{n+r}(r)\bigr)
\om_{n+r,n+s}^y(f).
$$
Note that the sum over $n$ is finite.
Lemma \thtag{275}
implies $\{ q^{-k} f^{s-k}_{n+s}(s)
\}_{k=-\infty}^\infty \in \ell^2(\Z)$, since
$q^{-k}f^{s-k}_{n+s}(s)={\Cal O}(q^{-k(1+|n|)})$ as $k\to-\infty$
and  $q^{-k}f^{s-k}_{n+s}(s)={\Cal O}(q^{\hf k(k-1)} 
q^{k(1-2s-n)})$ as $k\to \infty$.
So we can take the limit $N\to\infty$. Recall that 
$$
\sum_{k=-\infty}^\infty q^{-2k} f^{s-k}_m(s) f^{r-k}_{m+r-s}(r) =
\de_{r,s}q^{-2m},
\tag\eq{295}
$$
which is \cite{\CiccKK, (4.3)} for the special case $\al$ and
$c$ replaced by $s-m$ and $q^{2-2s}$ in base $q^2$. 
Now we can use \eqtag{295} to find
$(h\star\om^y_{r,s})(f) = \de_{y,r-s}\de_{r,s}
\sum_{n=-\infty}^\infty q^{-2n}\om^0_{n,n}(f)$, which is the
desired result.

For the other statement we proceed analogously, now using
$$
\sum_{r=-\infty}^\infty f^{r-l}_{n+l-k}(r)\, f^{r-k}_n(r)
\, q^{-2r} = \de_{k,l}q^{-2n},
$$
which is the same sum as \eqtag{295} 
using $f^{-x}_m(k)=f^x_m(k+x)$, see Lemma \thtag{275}.
\qed\enddemo

\demo{Proof of Theorem \thtag{280}} 
Let $\xi=\sum_{x,s=-\infty}^\infty
\xi_{x,s}\, f_{x,s}\in L^2(X,\mu)$ with only finitely
many $\xi_{x,s}\not=0$, then
$\om_\xi = \sum_{r,s,y\in\Z}
\bigl( \sum_{x-x'=y}\xi_{x',s}\bar \xi_{x,r}
\bigr) \om^y_{r,s}$, so that for $f\in C_c(\Z;C_c(X))$ we
have
$$
(h\star \om_\xi)(f) =
\sum_{r,s,x,n,k=-\infty}^\infty  (-1)^{r-s}
\bar \xi_{x,r}\xi_{x-r+s,s} q^{-2k} f^{s-k}_{n+s}(s)
f^{r-k}_{n+r}(r) \, \om^{r-s}_{n+r,n+s}(f),
\tag\eq{299}
$$
provided that the sum is absolutely convergent. 
If this holds, we can use \eqtag{295} to find the result. 
The Cauchy-Schwarz inequality applied to \eqtag{295} 
gives
$$
\sum_{k=-\infty}^\infty
|q^{-2k} f^{s-k}_{n+s}(s) f^{r-k}_{n+r}(r)| \leq
q^{-2n-r-s}.
\tag\eq{296}
$$
Use of the estimate \eqtag{296} for the sum over
$k$ leads to the termwise estimate
$$
\multline
\sum_{r,s,x,n=-\infty}^\infty
|\xi_{x,r}| |\xi_{x-r+s,s}| q^{-2n-r-s}\, |\om^{r-s}_{n+r,n+s}(f)|
= \\
\sum_{r,s=-\infty}^\infty |\om^{r-s}_{r,s}(f)| q^{-s-r}
\sum_{x,n=-\infty}^\infty |\xi_{x+r,r-n}| |\xi_{x+s,s-n}|,
\endmultline
\tag\eq{297}
$$
for the right hand side of \eqtag{299}.
The sum over $x,n$ is estimated by $\| \xi\|^2$, and
$\om^{r-s}_{r,s}(f)=0$ for $r> N$, $|r-s|>M$ for some
$N,M\in\N$. Hence, the sum is absolutely convergent 
for $f\in C_c(\Z;C_c(X))$ and
the result follows.
For $\om_\xi\star h$ we proceed analogously.
\qed\enddemo

\demo{Remark \thname{2100}} We can rewrite the
Haar functional $h$ in a coordinate free way. Use 
the covariant representation $\pt$ of 
$C_0(X)\times_\tau\Z$ to get 
$\om_{k,l}^x(f)=(2\pi)^{-1}\int_0^{2\pi}
\langle \pt(f)e_l,e_k\rangle e^{-ix\th}\, d\th$ and introduce the
unbounded operator $Q\colon {\Cal D}(Q)\to\H$,
$e_k\mapsto q^{-k}e_k$, which is self-adjoint on its
maximal domain ${\Cal D}(Q)=\{ \sum_k c_ke_k\mid
\sum_k |c_k|^2 q^{-2k}<\infty\}$.
We can rewrite
$h$ defined in Theorem \thtag{260} by
$$
h(f) = {1\over{2\pi}} \int_0^{2\pi}
\text{\rm Tr}\vert_{\H}(\pt(f) Q^2) \, d\th
\tag\eq{2105}
$$
for any $f\in C_0(X)\times_\tau\Z$ such that
$\pt(f)Q^2$ is of trace class and $\th\mapsto
\text{\rm Tr}\vert_{\H}(\pt(f)Q^2)$ is integrable. 
Note that for any decomposable operator
$T\in{\Cal B}(L^2(\T;\H))$, i.e. $T=(2\pi)^{-1}\int_0^{2\pi}
T(e^{i\th})\, d\th$, with $T(e^{i\th})Q^2$ of trace class
in $\H$ and $\th\mapsto \text{Tr}\vert_{\H}(T(e^{i\th})Q^2)$
integrable we can define $h(T)$ by \eqtag{2105}. Note that
the order of the operators in the trace in \eqtag{2105} is important.
E.g. define the bounded operator $S$ on $\H$ by 
$Se_k=e_{-k}$ for $k\geq 0$ and $Se_k=0$ for $k<0$, then $Q^2S$
is of trace class and $SQ^2$ is unbounded. 
\enddemo

\demo{Remark \thname{2110}} Recall that every element from the
algebra $A_q(SU(1,1))$ can be written uniquely as
a sum of elements
of the form $\al^k\ga^l p(\ga^\ast\ga)$, $k,l\in\Zp$, 
$\al^k(\ga^\ast)^l p(\ga^\ast\ga)$, $k\in\Zp$, $l\in\N$,
$(\al^\ast)^k\ga^l p(\ga^\ast\ga)$, $l\in\Zp$, $k\in\N$,  and
$(\al^\ast)^k(\ga^\ast)^l p(\ga^\ast\ga)$, 
$k,l\in\N$, where $p$ is a polynomial, cf. 
Theorem~3.6. %harde referentie
We can give a meaning to
the Haar functional evaluated on such elements if we change
$p$ from polynomials to sufficiently decreasing functions.
Let us do this explicitly for an element of the first type.
Applying $\pt$ we see that
$\al^k\ga^l p(\ga^\ast\ga)$ corresponds
in the representation $\pt$ to the operator
$U^k e^{il\th} (-Q^2;q^{-2})_k^\hf Q^lp(Q^2)$ on $\H$,
where $U\colon\H\to\H$, $e_k\mapsto e_{k+1}$, is the
unilateral shift. 
Hence, for a function $p$ satisfying
$\sum_{r\in\Z} (-q^{-2r};q^{-2})_k^\hf q^{-r(l+2)}|p(q^{-2r})|
<\infty$ we see that the corresponding operator times $Q^2$
is of trace class on $\H$.
Its trace is non-zero only for $k=0$, and for $k=0$ we
have $\text{Tr}\vert_{\H}e^{il\th} Q^{l+2}p(Q^2)= e^{il\th}
\sum_{r\in\Z} q^{-r(l+2)} p(q^{-2r})$. Integrating
over the circle gives zero unless $l=0$, hence for
$p$ sufficiently rapidly decreasing we have
$$
h(\al^k\ga^l p(\ga^\ast\ga))=
\de_{k,0}\de_{l,0} \sum_{r=-\infty}^\infty
q^{-2r} p(q^{-2r}) = \de_{k,0}\de_{l,0} {1\over{1-q^2}}
\int_0^\infty p(x)\, d_{q^2}x,
$$
where the last equality defines the Jackson $q$-integral
on $(0,\infty)$. In a similar fashion the Haar functional
applied to any of the other types of elements of $A_q(SU(1,1))$
described in the beginning of this remark gives zero.
So we see that Theorems \thtag{260} and \thtag{280} with Remark
\thtag{2100} correspond precisely to
\cite{\MasuMNNSU}, \cite{\MasuW, Lemme~2.2}, \cite{\VaksK}.
So we have linked the Haar functional 
to the Jackson integral on $(0,\infty)$ when restricted
to the subalgebra corresponding to the self-adjoint
element $\ga^\ast\ga$, see \cite{\Kake}, \cite{\KakeMU},
\cite{\VaksK} for the further analysis.
\enddemo

%%%%%%%%%%%%%%%%%%%%%%%%%%%%%%%%%%%%%%%%%%%%%%%%%%%%%%%%%%%%%%%%%%%%
%N E W   S E C T I O N%
%%%%%%%%%%%%%%%%%%%%%%%%%%%%%%%%%%%%%%%%%%%%%%%%%%%%%%%%%%%%%%%%%%%%
\head \newsec The quantised universal enveloping algebra
and self-adjoint elements\endhead

In this section we gather the necessary algebraic results on the
quantised universal enveloping algebra $\U$, which is the
dual Hopf $\ast$-algebra to $A_q(SU(1,1))$ introduced in \S 2,
see \cite{\CharP} for generalities on quantised
universal enveloping algebras and Hopf $\ast$-algebras. The
proofs of all statements in this section are analogous to the
corresponding statements for the compact quantum $SU(2)$
group, see \cite{\KoelSIAM}, \cite{\KoelAAM},
\cite{\KoelFIC}, \cite{\KoorZSE}, and are skipped or only
indicated. 
The main idea is due to
Koornwinder \cite{\KoorZSE} resulting into a quantum group
theoretic interpretation of a two-parameter family of 
the Askey-Wilson polynomials
as spherical functions. Then Noumi and Mimachi, see 
\cite{\Noum}, \cite{\NoumMPJA}, \cite{\NoumMLNM}, have
given an interpretation of the full four parameter
family of Askey-Wilson polynomials, see also \cite{\KoelSIAM},
\cite{\KoelAAM}, \cite{\KoelFIC}. As indicated by the
results in \cite{\KoelVdJ} the algebraic methods apply
to $\U$ as well in case of the positive discrete series
representations. 

%%%%%%%%%%%%%%%%%%%%%%%%%%%%%%%%%%%%%%%%%%%%%%%%%%%%%%%%%%%%%%%%%%%%
\subhead \the\sectionnumber.1 
The quantised universal enveloping algebra\endsubhead 
This subsection is a reminder and is used to fix the notation. 
The material of this subsection is standard, and we refer
to e.g. \cite{\CharP} for further information. 
By $\UC$ we denote the algebra generated by
$A$, $B$, $C$ and $D$ subject to the relations, where
$0<q<1$,
$$
AD=1=DA, \quad AB=qBA,\quad AC=q^{-1}CA,\quad
BC-CB = {{A^2-D^2}\over{q-q^{-1}}}. 
\tag\eq{310}
$$
It follows from \eqtag{310} that the element
$$
\Om = {{q^{-1}A^2+qD^2-2}\over{(q^{-1}-q)^2}}+BC=
{{q^{-1}D^2+qA^2-2}\over{(q-q^{-1})^2}}+CB
\tag\eq{315}
$$
is a central element of $\UC$, the Casimir element. 
The algebra $\UC$ 
is in fact a Hopf algebra with comultiplication 
$\De\colon \UC\to \UC\otimes \UC$ given
by
$$
\De(A)=A\otimes A, \quad \De(B) =A\otimes B+B\otimes D,
\quad \De(C)=A\otimes C+C\otimes D,\quad \De(D)=D\otimes D.
$$
The Hopf algebra $\UC$ is in duality
with the Hopf algebra $A_q(SL(2,\C))$ of the previous section,
where the duality is incorporated by the representations of 
Theorem~3.1. %harde referentie 
There are two ways to introduce a $\ast$-operator in order
to make $\UC$ a Hopf $\ast$-algebra. The first one is defined
by its action on the generators as follows; 
$A^\ast=A$, $B^\ast=-C$,
$C^\ast=-B$, $D^\ast=D$. We call the corresponding Hopf
$\ast$-algebra $\U$. The other $\ast$-structure is given
by $A^\times=A$, $B^\times=C$, $C^\times=B$, $D^\times=D$.
The corresponding Hopf $\ast$-algebra is denoted by
$\Uc$. Note that $\Om^\ast=\Om=\Om^\times$. 

\proclaim{Theorem \thname{320}} {\rm (See \cite{\CharP, Ch.~10})}
For each spin $l\in\hZp$ there
exists a unique $(2l+1)$-dimensional re\-pre\-sen\-tation
of $\UC$ such that the spectrum of $A$ is contained in $q^{\hf\Z}$.
Equip $\C^{2l+1}$ with orthonormal basis $\{ e^l_n\}$,
$n=-l,-l+1,\ldots,l$ and denote the representation by $t^l$. The
action of the generators is given by
$$
\aligned
&t^l(A)\, e^l_n=q^{-n}e^l_n,\qquad t^l(D)\, e^l_n=q^n\, e^l_n, \\
&t^l(B)\, e^l_n=
{{\sqrt{(q^{-l+n-1}-q^{l-n+1})(q^{-l-n}-q^{l+n})}}\over
{q^{-1}-q}} \,  e^l_{n-1} \\
&t^l(C)\, e^l_n=
{{\sqrt{(q^{-l+n}-q^{l-n})(q^{-l-n-1}-q^{l+n+1})}}\over
{q^{-1}-q}}\,  e^l_{n+1},
\endaligned
\tag\eq{330}
$$
where $e^l_{l+1}=0=e^l_{-l-1}$.
\endproclaim

The representation $t^l$ of $\UC$ is not a 
$\ast$-representation of $\U$, but it is a $\ast$-representation
of $\Uc$.

Then $A_q(SL(2,\C))$ is spanned by the matrix elements
$X\mapsto t^l_{n,m}(X)=\langle t^l(X)e^l_m,e^l_n\rangle$, and
the link is given by
$$
t^\hf = \pmatrix t^\hf_{-\hf,-\hf} 
& t^\hf_{-\hf,\hf} \\
t^\hf_{\hf,-\hf} & t^\hf_{\hf,\hf}\endpmatrix =
\pmatrix \al & \be \\ \ga & \de\endpmatrix.
\tag\eq{340}
$$
Also
$$
t^1 = \pmatrix \al^2 &\sqrt{1+q^2}\be\al & \be^2 \\
\sqrt{1+q^2}\ga\al & 1+ (q+q^{-1})\be\ga & \sqrt{1+q^2}\de\be\\
\ga^2 & \sqrt{1+q^2}\de\ga & \de^2\endpmatrix.
\tag\eq{350}
$$

%%%%%%%%%%%%%%%%%%%%%%%%%%%%%%%%%%%%%%%%%%%%%%%%%%%%%%%%%%%%%%%%%%%%
\subhead \the\sectionnumber.2 
Self-adjoint elements in $A_q(SU(1,1))$\endsubhead 
The definition of the self-adjoint elements in
$A_q(SU(1,1))$ we give here is strongly motivated by the
paper by Koornwinder \cite{\KoorZSE} for the compact 
quantum $SU(2)$ group and the results of \cite{\KoelVdJ}
for the positive discrete series representations of
$\U$. We define
$$
Y_s=q^\hf B-q^{-\hf} C +{{s+s^{-1}}\over{q^{-1}-q}}(A-D)\in\U,
\tag\eq{360}
$$
then $Y_s$ is twisted primitive, i.e. $\De(Y_s)=
A\otimes Y_s+Y_s\otimes D$, 
and $Y_sA=(Y_sA)^\ast$ is self-adjoint for
$s\in\R\backslash\{0\}$, and without loss of generality
we assume $|s|\geq 1$. The convention is 
$$
Y_\infty = \lim_{s\to 0} s(q^{-1}-q)\, Y_s = 
\lim_{s\to \infty} s^{-1}(q^{-1}-q)\, Y_s = A-D.
\tag\eq{365}
$$
The definition of $Y_s$ is as in \cite{\KoelVdJ}.

Then $t^l(Y_sA)\in \text{Mat}_{2l+1}(\C)$ is completely
diagonalisible. Since the proof is completely
analogous to the proof of \cite{\KoorZSE, Thm.~4.3}, we
do not give the proof here.

\proclaim{Lemma \thname{370}}
The tridiagonal matrix 
$t^l(Y_sA)$ is completely diagonalisible with spectrum
$$
\la_y(s) = {{sq^{2y}+s^{-1}q^{-2y} - (s+s^{-1})}\over{q^{-1}-q}},
\qquad y\in\{ -l, -l+1,\ldots, l\}, 
$$
and corresponding eigenvector
$$
v^l_y(s) = \sum_{n=-l}^l {{(q^{4l};q^{-2})_{l-n}^\hf}\over
{(q^2;q^2)_{l-n}^\hf }}   q^{\hf (l-n)(l-n-1)} s^{l-n}
R_{l-n}(q^{2y-2l}+s^{-2}q^{-2y-2l};s^{-2},2l;q^2) \, e^l_n
$$
where
$$
R_n(q^{-x}+cq^{x-N};c,N;q)= {}_3\vp_2\left(
{{q^{-n},q^{-x},cq^{x-N}}\atop{q^{-N},0}};q,q\right)
$$
is a dual $q$-Krawtchouk polynomial.
\endproclaim
 
It is straightforward to check that $t^l(Y_sA)^\ast$,
where $\ast$ denotes the adjoint
of a $(2l+1)\times (2l+1)$-matrix, equals $-t^l(Y_{-s}A)$,
since $(Y_sA)^\times = -Y_{-s}A$.
Note that $\la_y(-s)=-\la_y(s)$, $\la_y(s)=\la_{-y}(s^{-1})$, 
and that $\la_y(s)\not=\la_{y'}(s)$ for $y\not=y'$,
$y,y'\in\{-l,\ldots,l\}$ 
if $s^2\notin q^{2\Z}$.  

Next we define the matrix
elements with respect to the basis of eigenvectors of
$t^l(Y_sA)$ by $a^l_{i,j}(s,t)(X) =
\langle t^l(X) v^l_j(t), v^l_i(-s)\rangle$ for
real $s,t$ satisfying $|s|,|t|\geq 1$, 
so that
$$
a^l_{i,j}(s,t)(XY_tA) = \la_j(t) a^l_{i,j}(s,t)(X),
$$
and
$$
\multline
a^l_{i,j}(s,t)(Y_sAX) =
\langle t^l(X) v^l_j(t), \bigl(t^l(Y_sA)\bigr)^\ast 
v^l_i(-s)\rangle
\\ = - \langle t^l(X) v^l_j(t), t^l(Y_{-s}A) v^l_i(-s)\rangle =
- \la_i(-s) a^l_{i,j}(s,t)(X) = \la_i(s) a^l_{i,j}(s,t)(X) .
\endmultline
$$
Or, using the notation $X.\xi,\xi.X \in A_q(SL(2,\C))$ 
with $X,Y\in \UC$ 
for the elements defined by $X.\xi(Y)=\xi(YX)$ and 
$\xi.X(Y)=\xi(XY)$, we have
$$
\gathered
(Y_tA).a^l_{i,j}(s,t) = \la_j(t) a^l_{i,j}(s,t), \qquad
a^l_{i,j}(s,t).(Y_sA) = \la_i(s) a^l_{i,j}(s,t), \\
Y_t.b^l_{i,j}(s,t) = \la_j(t) D.b^l_{i,j}(s,t), \qquad
b^l_{i,j}(s,t).Y_s = \la_i(s) b^l_{i,j}(s,t).D
\endgathered
\tag\eq{380}
$$
for $b^l_{i,j}(s,t) = A.a^l_{i,j}(s,t)$. Here we have
used $X.(Y.\xi)=(XY).\xi$, $(\xi.X).Y=\xi.(XY)$
and $(X.\xi).Y=X.(\xi.Y)$, i.e. the applications
$X.\xi$ and $\xi.X$ define mutual compatible left and
right actions of $\UC$ on $A_q(SL(2,\C))$. 

As before assume $s,t\in\R$, $|s|,|t|\geq 1$. 
Write $v^l_y(s) = \sum_{n=-l}^l v^l_y(s)_n e^l_n$, then 
$a^l_{i,j}(s,t)= \sum_{n,m=-l}^l v^l_j(t)_m 
\overline{ v^l_i(-s)_n} t^l_{n,m}$
and $b^l_{i,j}(s,t)= 
\sum_{n,m=-l}^l v^l_j(t)_m \overline{ v^l_i(-s)_n}
q^{-m}t^l_{n,m}$, since we have $A.t^l_{n,m}=q^{-m}t^l_{n,m}$.
The case $l=\hf$ of Lemma~\thtag{370} gives
$$
v^\hf_{-\hf}(s) = s^{-1} e^\hf_{-\hf} + e^\hf_\hf, \qquad
v^\hf_{\hf}(s) = s e^\hf_{-\hf} + e^\hf_\hf 
$$
and using \eqtag{340} we get
$$
\pmatrix b^\hf_{-\hf,-\hf}(s,t) & b^\hf_{-\hf,\hf}(s,t) \\
b^\hf_{\hf,-\hf}(s,t) & b^\hf_{\hf,\hf}(s,t)\endpmatrix =
\pmatrix -s^{-1}t^{-1}\al_{s,t} & -s^{-1}\be_{s,t} \\
t^{-1}\ga_{s,t}& \de_{s,t}\endpmatrix
\tag\eq{390}
$$
with
$$
\gathered
\al_{s,t} = q^{\hf}\al +q^{-\hf}t\be -q^{\hf}s\ga -
q^{-\hf}st\de, \quad 
\be_{s,t} = q^{\hf}t\al +q^{-\hf}\be - q^{\hf}st
\ga - q^{-\hf}s\de, \\
\ga_{s,t} = - q^{\hf}s\al - q^{-\hf}st\be +q^{\hf}\ga +
q^{-\hf}t\de, \quad
\de_{s,t} = - q^{\hf}st\al - q^{-\hf}s\be +
q^{\hf}t\ga + q^{-\hf}\de
\endgathered
\tag\eq{395}
$$

Similarly, using the vector spanning the kernel of $t^1(Y_sA)$; 
$$
v^1_0(s) = q^{-1}e^1_{-1} + {{s+s^{-1}}\over{\sqrt{1+q^2}}} e^1_0
+ e^1_1
$$
and \eqtag{350} we
see that $b^1_{0,0}(s,t)$ equals, 
up to an affine transformation,
$$
\aligned
\rst &= \hf \bigl( \al^2+\de^2+q\ga^2+q^{-1}\be^2
+(t+t^{-1})(q\de\ga+\be\al)\\ &\qquad -(s+s^{-1})(q\ga\al+\de\be)
-(t+t^{-1})(s+s^{-1})\be\ga\bigr) \\ &=  
\hf\bigl( 
\al^2+(\al^\ast)^2+q(\ga^2+(\ga^\ast)^2)  +
q(t+t^{-1})(\al^\ast\ga+\ga^\ast\al) \\ &\qquad - 
q(s+s^{-1})(\ga\al+\al^\ast\ga^\ast)- q(t+t^{-1})(s+s^{-1})
\ga\ga^\ast\bigr).
\endaligned
\tag\eq{3100}
$$
Remark that \eqtag{380} remains valid for $\rst$ instead of
$b^1_{0,0}(s,t)$, since $1\in A_q(SL(2,\C))$ satisfies
$Y_t.1=0=1.Y_s$. 
Observe that $\rst=\rho_{s^{\pm 1},t^{\pm 1}}=\rst^\ast$,
since $s$ and $t$ are real. 

\demo{Remark \thname{3101}}
There is also a certain symmetry between $s$ and $t$. To be 
explicit, let $A_q(SU(1,1))^{\text{opp}}$ be the opposite Hopf 
$\ast$-algebra, see e.g. \cite{\CharP}, then interchanging $\ga$
and $\ga^\ast$ gives a Hopf $\ast$-algebra isomorphism 
$\psi\colon A_q(SU(1,1))\to A_q(SU(1,1))^{\text{opp}}$ which maps
$\rst$ to $\rho_{-t,-s}$. 
\enddemo

We calculate the Haar weight on the subalgebra generated
by the self-adjoint element $\rst\in A_q(SU(1,1))$ 
in \S 5 explicitly. 

The following limit case plays an
important role in the sequel;
$$
\multline
\rit = \lim_{s\to 0}
{{-2s}\over q} \rst = \lim_{s\to\pm \infty} {-2\over{qs}} \rst
= \ga\al+q^{-1}\de\be+(t+t^{-1})q^{-1}\be\ga \\ 
= \al^\ast\ga^\ast+\ga\al+(t+t^{-1})\ga\ga^\ast.
\endmultline
\tag\eq{3105}
$$
This element also satisfies 
$\rit=\rho_{\infty,t^{\pm 1}}=\rit^\ast$,
and we calculate the Haar weight on the subalgebra generated
by the self-adjoint element $\rit\in A_q(SU(1,1))$ 
in \S 4 explicitly. We need the appropriate limit case
of \eqtag{395}:
$$
\gathered
\al_{\infty,t} = \lim_{s\to 0}\al_{s,t} =
q^{\hf}\al +q^{-\hf}t\be, \quad
\be_{\infty,t} = \lim_{s\to 0}\be_{s,t} =
q^{\hf}t\al +q^{-\hf}\be, \\
\ga_{\infty,t} = \lim_{s\to 0}\ga_{s,t} =
q^{\hf}\ga + q^{-\hf}t\de, \quad
\de_{\infty,t} = \lim_{s\to 0}\de_{s,t} =
q^{\hf}t\ga + q^{-\hf}\de.
\endgathered
\tag\eq{3107}
$$
Note that we can express the elements defined in
\eqtag{395} in terms of these elements by 
$$
\gathered
\al_{s,t} = \al_{\infty,t} - s \ga_{\infty,t}, \quad
\be_{s,t} =  \be_{\infty,t} - s \de_{\infty,t}  , \\
\ga_{s,t} = \ga_{\infty,t} - s \al_{\infty,t} , \quad
\de_{s,t} = \de_{\infty,t} - s \be_{\infty,t} .
\endgathered
\tag\eq{3108}
$$

%%%%%%%%%%%%%%%%%%%%%%%%%%%%%%%%%%%%%%%%%%%%%%%%%%%%%%%%%%%%%%%%%%%%
\subhead \the\sectionnumber.3 Cartan decomposition\endsubhead 
The matrix elements $t^l_{n,m}$, $l\in\hZp$, 
$n,m\in\{ -l,-l+1,\ldots,l\}$ form a linear basis for
$A_q=A_q(SL(2,\C))$. Put $A_q^l=
\text{span}_\C \{ t^l_{n,m} \mid n,m=-l,\ldots,l\}$
for $l\in \hZp$. Note that $b^l_{n,m}(s,t)$, 
$n,m\in\{ -l,-l+1,\ldots,l\}$, form a basis for
$A_q^l(SL(2,\C))$ as well if $s^2,t^2\not\in q^{2\Z}$.

\proclaim{Proposition \thname{3110}} Let $s^2,t^2\notin q^{2\Z}$.
\par
\noindent
{\rm (i)}  Let $\xi\in A_q^l$, $\l\in\hZp$,
 be a $(s,t)$-spherical element, i.e.
$Y_t.\xi=0=\xi.Y_s$,
and let $\eta\in A_q$ satisfy
$$
Y_t .\eta=\lambda\, D.\eta
\quad\text{and}\quad
\eta.Y_s = \mu\, \eta.D
\tag\eq{3120}
$$
for some $\lambda,\mu\in\C$. Then $\eta\xi$ satisfies
\eqtag{3120} for the same $\lambda$, $\mu$.
Moreover, if $\lambda,\mu\in\R$, then
$\eta^\ast\eta$ is a $(s,t)$-spherical element.\par 
\noindent
{\rm (ii)}
If $\eta\in A_q^l$ satisfies \eqtag{3120} for
some $\lambda, \mu\in\C$ and $\eta$ is non-zero, then
$\lambda = \lambda_j(t)$, $\mu=\lambda_i(s)$ for some
$i,j\in\{-l,-l+1,\ldots,l\}$ and $\eta$ is a multiple of
$b^l_{i,j}(s,t)$.
\endproclaim

It follows that $\rst$ generates the $\ast$-subalgebra of
$A_q(SL(2,\C))$ of $(s,t)$-spherical elements for 
$s^2,t^2\notin q^{2\Z}$. 

\proclaim{Proposition \thname{3130}} 
Let $\eta\in A_q$ satisfy \eqtag{3120} with
$\la=\la_j(t)$ and
$\mu=\la_i(s)$, then
\item{\rm (i)} $\al_{sq^{2i},tq^{2j}}\eta$
satisfies \eqtag{3120} with
$\la=\la_{j-1/2}(t)$ and $\mu=\la_{i-1/2}(s)$,
\item{\rm (ii)} $\be_{sq^{2i},tq^{2j}}\eta$
satisfies \eqtag{3120} with
$\la=\la_{j+1/2}(t)$ and $\mu=\la_{i-1/2}(s)$,
\item{\rm (iii)} $\ga_{sq^{2i},tq^{2j}}\eta$
satisfies \eqtag{3120} with
$\la=\la_{j-1/2}(t)$ and $\mu=\la_{i+1/2}(s)$,
\item{\rm (iv)} $\de_{sq^{2i},tq^{2j}}\eta$
satisfies \eqtag{3120} with
$\la=\la_{j+1/2}(t)$ and $\mu=\la_{i+1/2}(s)$.
\par\noindent
In case $s=\infty$ the result remains valid 
with $Y_\infty$ defined
in \eqtag{365} and $\la_j(\infty)=q^{2j}-1$.
\endproclaim

We skip the proofs of Propositions \thtag{3110} and \thtag{3130},
since they are completely analogous to the 
proofs of \cite{\KoelAAM, Prop.~6.4, Prop.~6.5},
see also \cite{\KoelSIAM, Prop.~2.3}. We note that Proposition
\thtag{3130} can also be proved by direct calculations using
\eqtag{390}, $Y_t=Y_{tq^{2j}}-\la_j(t)(A-D)$, $Y_t$ being
twisted primitive and $\la_i(s)+\la_{\pm\hf}(sq^{2i})=
\la_{i\pm\hf}(s)$. 

A direct consequence of Propositions \thtag{3110} and \thtag{3130}
and Lemma \thtag{370} 
is the product structure of the matrix elements
$b^l_{i,j}(s,t)$ with $\max( |i|,|j|)=l$. For this we
define elements $\Ga^{(i)}_{l,m}(s,t)\in A_q(SU(1,1))$
for $m\in\{-l,-l+1,\ldots,l\}$, $l\in\hf\N$ in terms
of products of elementary elements
with the convention $\prod_{i=0}^k \xi_i=\xi_k\xi_{k-1}\ldots
\xi_0$ and the empty product being $1$; 
$$
\aligned
\Ga^{(1)}_{l,m}(s,t) &= \prod_{i=0}^{l+m-1}
\de_{sq^{l-m+i},tq^{m-l+i}} 
\prod_{j=0}^{l-m-1} \ga_{sq^j,tq^{-j}}= C_1\, b^l_{l,m}(s,t),\\ 
\Ga^{(2)}_{l,m}(s,t) &= \prod_{i=0}^{l-m-1}
\al_{sq^{l+m-i},tq^{-l-m-i}} 
\prod_{j=0}^{l+m-1} \ga_{sq^j,tq^{-j}}= C_2\, b^l_{m,-l}(s,t), \\ 
\Ga^{(3)}_{l,m}(s,t) &= \prod_{i=0}^{l+m-1}
\de_{sq^{m-l+i},tq^{l-m+i}} 
\prod_{j=0}^{l-m-1} \be_{sq^{-j},tq^j} = C_3\, b^l_{m,l}(s,t), \\ 
\Ga^{(4)}_{l,m}(s,t) &= \prod_{i=0}^{l+m-1}
\be_{sq^{-l+m-i},tq^{-l+m+i}} 
\prod_{j=0}^{l-m-1} \al_{sq^{-j},tq^{-j}} = C_4\, b^l_{-l,m}(s,t)
\endaligned
\tag\eq{3135}
$$
for certain non-zero constants $C_i$. 
Initially, the second equality in each line
of \eqtag{3135} holds for $s^2,t^2\notin q^{2\Z}$, and this
condition can be removed by continuity. An analogous expression
as in \eqtag{3135} holds for the case $s=\infty$, where
$\Ga^{(i)}_{l,m}(\infty,t)=\lim_{s\to 0} \Ga^{(i)}_{l,m}(s,t)$
by \eqtag{3107}. 

The explicit expression of $b^l_{i,j}(s,t)$ for
$\max(|i|,|j|)=l$ in \eqtag{3135} 
and Propositions \thtag{3110} and \thtag{3130} imply 
the following Cartan-type decomposition of $A_q(SL(2,\C))$. 

\proclaim{Theorem \thname{3140}} Let $s^2,t^2\not\in q^{2\Z}$.
$A_q(SU(1,1))$ is a free right $\C[\rst]$-module, with 
$\C[\rst]$-basis given by
$$
\{1\}\cup\{ \Ga^{(1)}_{l,m}\}_{l\in\hf\N, m\in I_1^l}
\cup \{ \Ga^{(2)}_{l,m}\}_{l\in\hf\N, m\in I_2^l}
\cup \{ \Ga^{(3)}_{l,m}\}_{l\in\hf\N, m\in I_3^l}
\cup \{ \Ga^{(4)}_{l,m}\}_{l\in\hf\N, m\in I_4^l},
$$
where $I_1^l=I_4^l=\{1-l,2-l,\ldots,l\}$, 
$I_2^l=\{-l,1-l,\ldots,l\}$,
and $I_3^l=\{1-l,2-l,\ldots,l-1\}$. 
\endproclaim

So any element $\xi\in A_q(SU(1,1))$
can be written as a finite sum of the form
$$
\xi=p(\rst)+\sum_{i=1}^4 \sum_{l\in\hf\N}\sum_{m\in I_i^l}
\Ga^{(i)}_{l,m}(s,t)\,  p^{(i)}_{l,m}(\rst) 
\tag\eq{3145}
$$
for uniquely determined polynomials $p$, $p^{(i)}_{l,m}$.

\demo{Remark \thname{3142}} {\rm (i)} 
We also have a corresponding decomposition
for the case $s=\infty$, and for the case $(s,t)=(\infty,\infty)$
we are back to the case discussed in  
\cite{\MasuMNNSU}, \cite{\MasuW}, see 
also Remark \thtag{2110}. Note that 
$h(\xi)$ with $\xi$ as in \eqtag{3145} is not well-defined,
but it can be defined properly after replacing the polynomials
in \eqtag{3145} by sufficiently decreasing functions, cf.
Remark \thtag{2110}. For $(s,t)=(\infty,\infty)$
the Cartan decomposition is formally orthogonal with respect
to the Haar functional $h$ by
$\langle\xi_1,\xi_2\rangle=h(\xi_2^\ast\xi_1)$, see \cite{\Kake}.
For the general case $(s,t)$ this is not clear.

\noindent
{\rm (ii)} The Cartan decomposition of Theorem \thtag{3140} is the
decomposition of $A_q(SU(1,1))$ into common eigenspaces 
of the left action of $AY_t$ and right action of $Y_sA$
on $A_q(SU(1,1))$, i.e. of the left
and right infinitesimal action of a ``torus'' depending on a
parameter. Since the Casimir operator $\Om$ defined in
\eqtag{315} commutes with these actions, the Casimir
operator preserves the Cartan decomposition.

\noindent
{\rm (iii)} For $\xi$ as in \eqtag{3145} we have
$$
\pi(\xi)=p(\pi(\rst))+\sum_{i,l,m}
\pi(\Ga^{(i)}_{l,m}(s,t))\,  p^{(i)}_{l,m}(\pi(\rst))
\tag\eq{3146}
$$
as an unbounded operator on $L^2(X,\mu)$. Now 
$\pi(\rst)$ is a symmetric unbounded operator. Suppose that 
${\Cal D}(s,t)$ is the domain of a self-adjoint extension
of $\pi(\rst)$ which is preserved by $\pi(\Ga^{(i)}_{l,m}(s,t))$,
then the right hand side of \eqtag{3146} makes sense as an 
unbounded linear operator on $L^2(X,\mu)$ with domain
${\Cal D}(s,t)$ for all continuous functions $p$, $p^{(i)}_{l,m}$
by the functional calculus of unbounded operators. If the
functions $p$, $p^{(i)}_{l,m}$ are such that the right hand side
of \eqtag{3146} are in $\pi(N_h)$, the formal decomposition 
\eqtag{3146} of the corresponding unique element $\xi\in N_h$
is called the Cartan decomposition of $\xi$. 
\enddemo

It follows from \eqtag{3135}, \eqtag{380}, 
Proposition \thtag{3110} and since $\rst$ generates
the algebra of $(s,t)$-spherical elements
that $(\Ga^{(i)}_{l,m}(s,t))^\ast\Ga^{(i)}_{l,m}(s,t)$ is
a polynomial in $\rst$. From \eqtag{3135} and \eqtag{3100}
we see that the degree of this polynomial is $2l$. If we
use the one-dimensional $\ast$-representation of
$A_q(SU(1,1))$ sending $\al$ to $e^{\hf i\th}$ and
$\ga$ to zero, $\th\in\R$, we obtain
$$
\aligned 
\Ga^{(1)}_{l,m}(s,t)^\ast\Ga^{(1)}_{l,m}(s,t) &=  C_1
(qste^{i\th}, qste^{-i\th};q^2)_{l+m} 
(q{s\over t}e^{i\th}, q{s\over t}e^{-i\th};q^2)_{l-m}
\vert_{\cos\th=\rst}, \\
\Ga^{(2)}_{l,m}(s,t)^\ast\Ga^{(2)}_{l,m}(s,t) &=  C_2 s^{2l-2m}
(q{s\over t}e^{i\th}, q{s\over t}e^{-i\th};q^2)_{l+m} 
({q\over{st}}e^{i\th}, {q\over{st}}e^{-i\th};q^2)_{l-m}
\vert_{\cos\th=\rst}, \\
\Ga^{(3)}_{l,m}(s,t)^\ast\Ga^{(3)}_{l,m}(s,t) &=  C_3 s^{2l-2m}
(qste^{i\th}, qste^{-i\th};q^2)_{l+m} 
(q{t\over s}e^{i\th}, q{t\over s}e^{-i\th};q^2)_{l-m}
\vert_{\cos\th=\rst}, \\
\Ga^{(4)}_{l,m}(s,t)^\ast\Ga^{(4)}_{l,m}(s,t) &=  C_4 s^{4l}
(q{t\over s}e^{i\th}, q{t\over s}e^{-i\th};q^2)_{l+m} 
({q\over{st}}e^{i\th}, {q\over{st}}e^{-i\th};q^2)_{l-m}
\vert_{\cos\th=\rst}, 
\endaligned
\tag\eq{3137}
$$
for positive constants $C_i$ independent
of $\th$ and $s$, cf. \cite{\KoelAAM, \S 7}.
In \eqtag{3137} we also have the appropriate case
for $s=\infty$ using \eqtag{3105} and \eqtag{3107};
$$
\aligned
\Ga^{(1)}_{l,m}(\infty,t)^\ast\Ga^{(1)}_{l,m}(\infty,t) &=  C_1
(-tq^2\rit;q^2)_{l+m} 
(-q^2\rit/t;q^2)_{l-m}, \\
\Ga^{(2)}_{l,m}(\infty,t)^\ast\Ga^{(2)}_{l,m}(\infty,t) &=  C_2
(-q^2\rit/t;q^2)_{l+m} 
(-tq^{2+2m-2l}\rit;q^2)_{l-m}, \\
\Ga^{(3)}_{l,m}(\infty,t)^\ast\Ga^{(3)}_{l,m}(\infty,t) &=  C_3
(-q^2t\rit;q^2)_{l+m} 
(-q^{2-2l+2m}\rit/t;q^2)_{l-m}, \\
\Ga^{(4)}_{l,m}(\infty,t)^\ast\Ga^{(4)}_{l,m}(\infty,t) &=  C_4
(-q^{2-2l-2m}\rit/t;q^2)_{l+m} 
(-tq^{2-2l+2m}\rit;q^2)_{l-m},
\endaligned
\tag\eq{3138}
$$
for positive constants $C_i$.

%%%%%%%%%%%%%%%%%%%%%%%%%%%%%%%%%%%%%%%%%%%%%%%%%%%%%%%%%%%%%%%%%%%%
\subhead \the\sectionnumber.4 
Factorisation and commutation results\endsubhead 
In order to obtain recurrence relations in later sections 
we need factorisation and commutation relations
in the algebra $A_q(SU(1,1))$.
The following corollary is consequence of 
Proposition \thtag{3110} and
Proposition \thtag{3130}, see 
\cite{\KoelFIC, \S 2} for the analogous statement
for the case $A_q(SU(2))$. Since the proof is the same
we skip it. 

\proclaim{Corollary \thname{3150}}
The following factorisation and commutation relations hold;
$$
\gather
\be_{sq,tq^{-1}}\ga_{s,t} = -2st\rst + q^{-1}t^2 +q s^2, \quad
\ga_{sq^{-1},tq}\be_{s,t} = -2st\rst + qt^2 +q^{-1} s^2, \\
\al_{sq,tq}\de_{s,t} = -2qst\rst + 1 +q^2 s^2t^2, \quad
\de_{sq^{-1},tq^{-1}}\al_{s,t} = -2q^{-1}st\rst + 1 +q^{-2}s^2t^2,
\endgather
$$
and
$$
\gather
\al_{s,t}\rst = \rho_{sq^{-1},tq^{-1}}\al_{s,t}, \quad
\be_{s,t}\rst = \rho_{sq^{-1},tq}\be_{s,t}, \\
\ga_{s,t}\rst = \rho_{sq,tq^{-1}}\ga_{s,t}, \quad
\de_{s,t}\rst = \rho_{sq,tq}\de_{s,t}.
\endgather
$$
\endproclaim

Combining Corollary \thtag{3150} with \eqtag{3108} gives
$$
-2st\rst +q^{-1}t^2+qs^2 = \be_{sq,tq^{-1}}\ga_{s,t} =
(\be_{\infty,tq^{-1}}-sq\de_{\infty,tq^{-1}})
(\ga_{\infty,t}-s\al_{\infty,t}),
\tag\eq{3160}
$$
which is one of many ways of writing $\rst$ in products
of matrix elements $b^\hf_{i,j}(\infty,t)$. 

The limit case $s=\infty$, i.e. $s\to 0$, of
Corollary \thtag{3150} immediately gives the following.

\proclaim{Corollary \thname{3170}}
The following factorisation and commutation relations hold;
$$
\gather
q^{-1}\be_{\infty,tq^{-1}}\ga_{\infty,t} = t\rit + q^{-2}t^2, 
\quad
q^{-1}\ga_{\infty,tq}\be_{\infty,t} = t\rit + t^2, \\
q^{-2}\al_{\infty,tq}\de_{\infty,t} = t\rit +q^{-2}, \quad
\de_{\infty,tq^{-1}}\al_{\infty,t} = t\rit + 1,
\endgather
$$
and
$$
\gather
\al_{\infty,t}\rit = q\rho_{\infty,tq^{-1}}\al_{\infty,t}, \quad
\be_{\infty,t}\rit = q\rho_{\infty,tq}\be_{\infty,t}, \\
\ga_{\infty,t}\rit = q^{-1}\rho_{\infty,tq^{-1}}\ga_{\infty,t},
\quad
\de_{\infty,t}\rit = q^{-1}\rho_{\infty,tq}\de_{\infty,t}.
\endgather
$$
\endproclaim

%%%%%%%%%%%%%%%%%%%%%%%%%%%%%%%%%%%%%%%%%%%%%%%%%%%%%%%%%%%%%%%%%%%%
%N E W   S E C T I O N%
%%%%%%%%%%%%%%%%%%%%%%%%%%%%%%%%%%%%%%%%%%%%%%%%%%%%%%%%%%%%%%%%%%%%
\head \newsec The Haar functional on the
algebra generated by $\rit$\endhead

The Haar functional on the Cartan decomposition of 
Theorem \thtag{3140} for $(s,t)=(\infty,\infty)$ is related to
the Jackson integral on $(0,\infty)$, see 
\cite{\Kake}, \cite{\KakeMU}, \cite{\MasuW}, \cite{\VaksK}
and Remark \thtag{2110}. In this section we show that the
Haar functional on the Cartan decomposition of Theorem \thtag{3140}
for the case $s=\infty$, $t>q^{-1}$ finite, is related to the
Jackson integral on $[-d,\infty)$ for some $d>0$. 
The key ingredient is the spectral analysis of the 
unbounded symmetric operator $\pt(\rit)$ given in
\cite{\CiccKK} yielding an orthogonal basis
of eigenvectors of $\H$. To use the expression for $h$ 
of Remark \thtag{2100} we have to calculate the matrix elements
of $Q^2$ in this basis of eigenvectors in order to calculate
the trace. We also show that the
elements of \eqtag{3107} in the representation $\pt$
act as shift operators in the basis of eigenvectors. 
The results and approach are motivated by the results
for the quantum $SU(2)$ group case considered in 
\cite{\KoelV, \S 5} and we consider the Jackson integral on
$[-d,\infty)$ as the non-compact analogue of the Jackson integral on
$[-d,c]$. The proofs are more involved due to
the unboundedness of the operators. 

%%%%%%%%%%%%%%%%%%%%%%%%%%%%%%%%%%%%%%%%%%%%%%%%%%%%%%%%%%%%%%%%%%%%
\subhead \the\sectionnumber.1 
Spectral analysis of $\pt(\rit)$\endsubhead 
Using Proposition \thtag{220} and \eqtag{3105} we get
$$
\pt(\rit)\, e_k = (t+t^{-1}) q^{-2k}\, e_k +
 e^{i\th} q^{-1-k}\sqrt{1+ q^{-2k}}\, e_{k+1}
+ e^{-i\th} q^{-k} \sqrt{1+  q^{2-2k}}\, e_{k-1},
$$
and by going over to the orthonormal basis $f_k=e^{ik\th}e_k$
we obtain
$$
\pt(\rit)\, f_k = (t+t^{-1}) q^{-2k}\, f_k +
q^{-1-k}\sqrt{1+q^{-2k}}\, f_{k+1}
+ q^{-k} \sqrt{1+q^{2-2k}}\, f_{k-1}.
$$
This operator is unbounded and symmetric on the domain consisting
of finite linear combinations of the basis vectors. 
This unbounded symmetric operator has been studied in detail in 
\cite{\CiccKK} for $t\in\R\backslash\{ 0\}$. 
It turns out that the operator in question
is essentially self-adjoint for $|t|\geq q^{-1}$.
Then the spectrum consists
completely of point spectrum plus one accumulation point at zero,
which itself is not in the point spectrum. 
The domain of $\pt(\rit)$ is its maximal domain, i.e.
$\{ \sum_k c_ke_k\in\H \mid \sum_k c_k \pt(\rit) e_k\in\H\}$.
Proposition~4.1 %harde referentie 
is the analogue of \cite{\KoelV, Prop.~5.2}
and has been proved in \cite{\CiccKK}. 

\proclaim{Proposition \thname{410}}
Let $t\in\R$ satisfy $|t|>q$.
There exists an orthogonal basis of $\H$ of
the form $\{v^\th_p(t)\mid p\in\Zp\}
\cup\{ w^\th_p(t)\mid p\in\Z\}$ given by
$$
v^\th_p(t) = \sum_{k=-\infty}^\infty e^{ik\th}
V_k(-q^{2p}t^{-1};t) \, e_k, \quad 
w^\th_p(t) = \sum_{k=-\infty}^\infty e^{ik\th}
V_k(q^{2p}t;t)\, e_k, 
$$
where
$V_k(x;t) = (-q^{2-2k};q^2)_\infty^\hf q^{\hf k(k+1)}
(-t)^{-k}  \ {}_1\vp_1 (-(xt)^{-1};q^2t^{-2};q^2,
xq^{2k+2}t^{-1})$.
These vectors are all contained in the maximal
domain of $\pt(\rit)$. Moreover,
$\pt(\rit)\, v^\th_p(t)=-q^{2p}t^{-1}\, v^\th_p(t)$, $p\in\Zp$,
and $\pt(\rit)\, w^\th_p(t)=q^{2p}t\, w^\th_p(t)$, $p\in\Z$.
The lengths of the orthogonal basis vectors are
given by
$$
\align
\| v^\th_p(t)\|^2 &= t^2 q^{-2p} 
{{(q^2;q^2)_p}\over{(q^2t^{-2};q^2)_p}}
{{(q^2,-t^{-2},-q^2t^2;q^2)_\infty}
\over{(q^2t^{-2};q^2)_\infty}}, \quad p\in\Zp,\\
\| w^\th_p(t)\|^2 &= q^{-2p} {{(-q^{2+2p},q^2,q^2,-t^{-2},
-q^2t^2;q^2)_\infty}\over{(-q^{2p+2}t^2,q^2t^{-2},
q^2t^{-2};q^2)_\infty}}, \quad p\in\Z.
\endalign
$$
For $|t|\geq q^{-1}$, the operator $\pt(\rit)$ is
self-adjoint on its maximal domain.
\endproclaim

In particular we find that the spectrum of
$\pt(\rit)$ consists of
$\si(\pt(\rit))= 
\{ -q^{2p}t^{-1}\}_{p\in\Zp}\cup \{q^{2p}t\}_{p\in\Z} \cup
\{ 0\}$, which is independent of $\th$.  
This also follows from $T(e^{i\th})
\pt(\rit)T(e^{i\th})^\ast=\pi_1(\rit)$
with $T(e^{i\th})$ the unitary operator on $\H$ defined by
$e_k\mapsto e^{-ik\th}e_k$. Using \cite{\DixmvNA, Ch. II.2, \S 6}
we find for a bounded continuous function $g$ 
$$
g(\pi(\rit)) = {1\over{2\pi}} \int_0^{2\pi} 
g(\pt(\rit))\, d\th  = T^\ast\bigl(\id \otimes g(\pi_1(\rit))
\bigr) T,
\tag\eq{415}
$$
$T={1\over{2\pi}}\int_0^{2\pi} T(e^{i\th})\, d\th$, 
using $L^2(\T;\H)\cong L^2(\T)\otimes\H$ as tensor product of
Hilbert spaces. So $g(\pi(\rit))$ is a 
decomposable operator on ${1\over{2\pi}}\int_0^{2\pi} 
\H\, d\th = L^2(\T;\H)$, since
$T$ commutes with multiplication by a function from $L^2(\T)$, 
see \cite{\DixmvNA, Ch.II.2, \S 5}.

Note that Proposition \thtag{410} gives an orthogonal decomposition
$\H=V^\th(t)\oplus W^\th(t)$, with 
$V^\th(t)$, respectively $W^\th(t)$,
the closure of the linear span of the vectors
$v^\th_p(t)$, $p\in\Zp$, respectively $w^\th_p(t)$, $p\in\Z$. 
Using \cite{\DixmvNA, Ch.~II.1} we obtain the decomposition
$L^2(\T;\H)=V(t)\oplus W(t)$ with
$V(t)=(2\pi)^{-1}\int_0^{2\pi} V^\th(t)\, d\th$ and
$W(t)=(2\pi)^{-1}\int_0^{2\pi} W^\th(t)\, d\th$.  

%%%%%%%%%%%%%%%%%%%%%%%%%%%%%%%%%%%%%%%%%%%%%%%%%%%%%%%%%%%%%%%%%%%%
\subhead \the\sectionnumber.2 Calculating the trace\endsubhead 
In this subsection we calculate the trace
of $g(\pt(\rit))Q^2$, cf. \eqtag{2105}, for
sufficiently decreasing function $g$ using the 
basis described in Proposition \thtag{410}. 
We start with the 
partial analogue of \cite{\KoelV, Lemma~5.5}. 
The operator $Q^2$ is self-adjoint with respect to its maximal
domain ${\Cal D}(Q^2)= 
\{\sum_k c_ke_k\in\H\mid \sum_k |c_k|^2q^{-4k}<\infty\}$. 
The proof of Lemma~4.2 %harde referentie
is a lengthy calculation, and
can be skipped at first reading.

\proclaim{Lemma \thname{420}} Let $|t|> q^{-1}$, then
$v^\th_p(t),w^\th_p(t)\in {\Cal D}(Q^2)$.
Moreover,
$$
{{\langle Q^2v^\th_p(t), v^\th_p(t)\rangle}\over{
\langle v^\th_p(t), v^\th_p(t)\rangle}} = 
{{q^{2p}}\over{t^2-1}},\ p\in\Zp, \quad
{{\langle Q^2w^\th_p(t), w^\th_p(t)\rangle}\over{
\langle w^\th_p(t), w^\th_p(t)\rangle}} = 
{{q^{2p}}\over{1-t^{-2}}},\ p\in\Z.
$$
\endproclaim

\demo{Proof}
Since $(-q^{2-2k};q^2)_\infty^\hf q^{\hf k(k+1)}=
{\Cal O}(q^k)$ as $k\to\infty$, we see that
$V_k(x;t)= {\Cal O}((q/t)^k)$ as $k\to\infty$. In order
to have $v^\th_p(t),w^\th_p(t)\in {\Cal D}(Q^2)$
we need $|qt|^{-k}$ to be square summable for
$k\to\infty$. Hence, we need $|t|>q^{-1}$. 
Assuming that $|t|>q^{-1}$ we see that
$Q^2v^\th_p(t)$ and $Q^2w^\th_p(t)$ are well-defined
elements of $\H$.

The calculation of the diagonal elements of $Q^2$
in this basis is based on orthogonality properties
that follow from Proposition \thtag{410}. The idea of the
proof is taken from the proof of \cite{\KoelV, Lemma~5.5},
but since it is much more involved we give the
details. 
Let us consider the first case, and
for this we introduce a moment functional for
the $q$-Laguerre polynomials defined by  
$$
{\Cal L} (f) = \sum_{k=-\infty}^\infty
(-q^{2-2k};q^2)_\infty q^{k(k+1)}
t^{-2k} f(q^{2k}), 
$$
cf. \cite{\CiccKK, \S 4}. 
Let $P_p(x)= {}_1\vp_1(q^{-2p};qt^{-2};q^2,-xq^{2+2p}t^{-2})$
be the corresponding orthogonal polynomials, which
are $q$-Laguerre polynomials in which the usual
parameter $\al$ of the $q$-Laguerre polynomials
corresponds to $t$ via $t^{-2}=q^{2\al}$. Then we have
$$
{\Cal L}(P_pP_m)=\de_{mp}
t^2 q^{-2p} {{(q^2;q^2)_p}\over{(q^2t^{-2};q^2)_p}}
{{(q^2,-t^{-2},-q^2t^2;q^2)_\infty}
\over{(q^2t^{-2};q^2)_\infty}}, \quad m,p\in\Zp,
$$
which follows from Proposition \thtag{410}, cf. \cite{\CiccKK}.
It turns out that we can calculate the
general matrix element 
$\langle Q^2v^\th_p(t), v^\th_m(t)\rangle$ 
without any extra difficulty. Since $Q^2$ is self-adjoint
and the vectors $v^\th_p(t)$ are in its domain,
we can assume without loss of generality that $m\leq p$.
The matrix element can be expressed in terms of the
moment functional;
$$
\langle Q^2v^\th_p(t), v^\th_m(t)\rangle =
{\Cal L}\bigl( x^{-1}P_m(x) P_p(x)\bigr) =
\text{CT}(P_m) {\Cal L}\bigl( x^{-1}P_p(x)),
$$
where $\text{CT}(P_m)$ means the constant term of the
polynomial $P_m$. This is valid 
since $x^{-1}P_m(x)= \text{CT}(P_m)x^{-1} + \text{polynomial
of degree less than }m$, since 
we have ${\Cal L}(P'P_p)=0$ for any polynomial $P'$ of degree
less than $p$. It remains to calculate
$$
\align
&{\Cal L}\bigl( x^{-1}P_p(x)) =\sum_{k=-\infty}^\infty
(-q^{2-2k};q^2)_\infty q^{k(k-1)}
t^{-2k} {}_1\vp_1(q^{-2p};q^2t^{-2};
q^2,-q^{2k+2p+2}t^{-2})\\
&= \sum_{r=0}^p {{(q^{-2p};q^2)_r}\over{(q^2,q^2t^{-2};q^2)_r}}
q^{r(r-1)}t^{-2r}q^{2r(p+1)}
\sum_{k=-\infty}^\infty
(-q^{2-2k};q^2)_\infty q^{k(k-1)}
t^{-2k} q^{2rk},
\endalign
$$
where interchanging summations is allowed since the
sum converges absolutely. Using the theta product
identity \eqtag{237} and Ramanujan's 
${}_1\psi_1$-summation formula, see \cite{\GaspR, (5.2.1)},
we can evaluate the inner sum as 
$$
(t^{-2};q^2)_r t^{2r} q^{-r(r-1)}
{{(q^2,-t^{-2},-q^2t^2;q^2)_\infty}
\over{(t^{-2};q^2)_\infty}}.
$$
Since $\text{CT}(P_m)=1$ we obtain for $m\leq p$ 
$$
\align
\langle Q^2v^\th_p(t), v^\th_m(t)\rangle &=
{{(q^2,-t^{-2},-q^2t^2;q^2)_\infty}
\over{(t^{-2};q^2)_\infty}} \ 
{}_2\vp_1(q^{-2p},t^{-2};q^2t^{-2};
q^2,q^{2p+2}) \\ &=
{{(q^2,-t^{-2},-q^2t^2;q^2)_\infty}
\over{(t^{-2};q^2)_\infty}} 
{{(q^2;q^2)_p}\over{(q^2t^{-2};q^2)_p}}
\endalign
$$
by the $q$-Chu-Vandermonde summation, see 
\cite{\GaspR, (1.5.2)}. Using the norms given in Proposition
\thtag{410} we obtain
$$
{{\langle Q^2v^\th_p(t), v^\th_m(t)\rangle}\over{
\|v^\th_p(t)\| \|v^\th_m(t)\|}} = 
{{t^{-2}q^{p+m}}\over{1-t^{-2}}} 
\left( {{(q^2;q^2)_p (q^2t^{-2};q^2)_m}\over
{(q^2;q^2)_m (q^2t^{-2};q^2)_p}}\right)^\hf.
\tag\eq{425}
$$
Now take $m=p$ to find the first statement. 

For the other statement we proceed in the same way. However,
this time we cannot get rid of a summation so easily since
there are no orthogonal polynomials around.  
We consider the functions
$$
\align
M_p(x)&={}_1\vp_1(-q^{-2p}t^{-2};q^2t^{-2};q^2,xq^{2p+2}) \\
&={{(xq^{2p+2};q^2)_\infty}\over{(q^2t^{-2};q^2)_\infty}}
{}_1\vp_1( -x;xq^{2p+2};q^2,q^2t^{-2})
\endalign
$$
by the transformation \eqtag{234}. Using \eqtag{237} we find 
$$
\multline
\langle Q^2w^\th_p(t),w^\th_m(t)\rangle = 
{\Cal L}\bigl( x^{-1}M_p(x)M_m(x)\bigr)\\
= {{(-1,-q^2;q^2)_\infty}\over{(q^2t^{-2};q^2)_\infty^2}}
\sum_{k=-\infty}^\infty {{t^{-2k}}\over{(-q^{2k};q^2)_\infty}}
(q^{2k+2p+2};q^2)_\infty \, {}_1\vp_1\left( {{-q^{2k}}\atop
{q^{2k+2p+2}}};q^2,{{q^2}\over{t^2}}\right) \\
\times
(q^{2k+2m+2};q^2)_\infty \, {}_1\vp_1\left( {{-q^{2k}}\atop
{q^{2k+2m+2}}};q^2,{{q^2}\over{t^2}}\right)
\endmultline
\tag\eq{430}
$$

Now consider the following generating functions, 
see \cite{\CiccKK, Lemma~5.1},
$$
\sum_{k=-\infty}^\infty z^kb^k {{(q^{k+1};q)_\infty}\over
{(aq^k/b;q)_\infty}} \, {}_1\vp_1\left( {{aq^k/b}\atop{q^{k+1}}};
q,bx\right) = {{(q,az,x/z;q)_\infty}\over{(a/b,bz;q)_\infty}},
\quad 0<|z|<|b|^{-1}
$$
and
$$
\sum_{n=-\infty}^\infty w^n (q^{n+1};q)_\infty
\, {}_1\vp_1\left( {{dq^{n+1}/y}\atop{q^{n+1}}};
q,y\right) = {{(d,q,y/w;q)_\infty}\over{(w,d/w;q)_\infty}}, 
\quad |d|<|w|<1,
$$
to see that for $0<|z|<|t|^2$ we have 
$$
{{(q^2,q^2/z,-t^{-2}q^{-2p}z;q^2)_\infty}\over{
(-q^{-2p},zt^{-2};q^2)_\infty}}z^{-p}t^{2p} = 
\sum_{k=-\infty}^\infty z^k t^{-2k}
{{(q^{2k+2p+2};q^2)_\infty}\over{(-q^{2k};q^2)_\infty}}
{}_1\vp_1\left( {{-q^{2k}}\atop{q^{2k+2p+2}}};q^2,
{{q^2}\over{t^2}}\right), 
$$
and for $1<|z|<|t|^2q^{2m}$ we have 
$$
{{(-t^{-2}q^{-2m},q^2,zq^2t^{-2};q^2)_\infty}\over{
(1/z,-t^{-2}q^{-2m}z;q^2)_\infty}} z^m
=\sum_{n=-\infty}^\infty z^{-n}
(q^{2n+2m+2};q^2)_\infty
{}_1\vp_1\left( {{-q^{2n}}\atop{q^{2n+2m+2}}};q^2,
{{q^2}\over{t^2}}\right).
$$
Note $|t|>q^{-1}$ 
by assumption and assume first the
additional condition $|t|^2q^{2m}>1$,
so that multiplying these generating functions is valid
in the annulus $1<|z|<\min(|t|^2,|t|^2q^{2m})$. 
The constant term
is the series in \eqtag{430}. Hence, this series
equals the constant term of
$$
{{(q^2,q^2,-t^{-2}q^{-2m};q^2)_\infty t^{2p}}\over
{(-q^{-2p};q^2)_\infty}}
{{z^{m-p}(-t^{-2}q^{-2p}z;q^2)_{p-m}}
\over{(1-1/z)(1-zt^{-2})}},
$$
assuming without loss of generality that $p\geq m$. 
Since $1<|z|<t^2$ we have 
$$
{1\over{(1-1/z)(1-zt^{-2})}} = \sum_{k=0}^\infty z^{-k}
\sum_{p=0}^\infty z^pt^{-2p} = 
{1\over{1-t^{-2}}}\bigl( \sum_{l=-\infty}^{-1} z^l 
+\sum_{l=0}^\infty t^{-2l}z^l\bigr) 
$$
and by the $q$-binomial formula \cite{\GaspR, (1.3.14)} 
we have
$$
z^{m-p}(-t^{-2}q^{-2p}z;q^2)_{p-m}= 
\sum_{k=0}^{p-m} {{(q^{2m-2p};q^2)_k}\over{(q^2;q^2)_k}}
z^{m-p+k} (-t^{-2}q^{-2m})^k.
$$
So 
$$
\multline
\text{CT}\left( {{z^{m-p}(-t^{-2}q^{-2p}z;q^2)_{p-m}}
\over{(1-1/z)(1-zt^{-2})}}\right) = 
\sum_{k=0}^{p-m} {{(q^{2m-2p};q^2)_k}\over{(q^2;q^2)_k}}
(-t^{-2}q^{-2m})^k {{t^{-2(p-m-k)}}\over{1-t^{-2}}} \\ =
{{t^{2(m-p)}}\over{1-t^{-2}}} (-q^{-2p};q^2)_{p-m},
\endmultline
$$
again by the $q$-binomial formula \cite{\GaspR, 1.3.14)}. 
We conclude that for $p\geq m$
$$
\multline
\sum_{k=-\infty}^\infty {{t^{-2k} (q^{2k+2p+2};q^2)_\infty}
\over{(-q^{2k};q^2)_\infty}}
{}_1\vp_1\left( {{-q^{2k}}\atop{q^{2k+2p+2}}};q^2,
{{q^2}\over{t^2}}\right) \\ \qquad\times
(q^{2k+2m+2};q^2)_\infty
{}_1\vp_1\left( {{-q^{2k}}\atop{q^{2k+2m+2}}};q^2,
{{q^2}\over{t^2}}\right) =
{{(q^2,q^2,-t^{-2}q^{-2m};q^2)_\infty}\over
{(-q^{-2m};q^2)_\infty}}
{{t^{2m}}\over{1-t^{-2}}}.
\endmultline
\tag\eq{440}
$$
This identity is valid under the extra assumption 
$|t|^2q^{2m}>1$. The right hand side of \eqtag{440}
is analytic in $t$ for $|t|>1$. Each summand on the left
hand side of \eqtag{440} is analytic in $t$ for $|t|>1$.
As $k\to\infty$ the summand behaves like
$t^{-2k}$, and as $k\to-\infty$ the summand
behaves like $t^{2k}q^{2(p+m)k}q^{k(k-1)}$
using \eqtag{238}, so that we obtain
uniform convergence on compact sets for the 
left hand side of \eqtag{440}. Hence, \eqtag{440} holds
for all $t$ with $|t|>1$.

Since \eqtag{440} gives the evaluation of the sum in
\eqtag{430} we obtain for $p\geq m$ 
$$
\align
\langle Q^2w^\th_p(t),w^\th_m(t)\rangle &= t^{2m}
{{(-1,-q^2;q^2)_\infty}\over{(q^2t^{-2};q^2)_\infty^2}}
{{(q^2,q^2,-t^{-2}q^{-2m};q^2)_\infty}\over{
(-q^{-2m};q^2)_\infty(1-t^{-2})}}\\
&= {{(q^2,q^2,-t^{-2},-q^2t^2,-q^{2+2m};q^2)_\infty}
\over{(q^2t^{-2},q^2t^{-2},-t^2q^{2+2m};q^2)_\infty(1-t^{-2})}}
\endalign
$$
and using the norm of $w^\th_p(t)$
given in Proposition \thtag{410} we obtain for $p\geq m$ 
$$
{{\langle Q^2w^\th_p(t),w^\th_m(t)\rangle }\over{
\|w^\th_p(t)\| \|w^\th_m(t)\|}} = 
{{q^{p+m}}\over{1-t^{-2}}} 
{{(-t^2q^{2+2p},-q^{2+2m};q^2)_\infty^\hf}\over
{(-t^2q^{2+2m},-q^{2+2p};q^2)_\infty^\hf}}.
\tag\eq{445}
$$
Now take $m=p$. 
\qed\enddemo

\demo{Remark} By an analogous computation we obtain
$$
\langle Q^2 v_p^\th(t), w_m^\th(t)\rangle = 
{{(-1,-q^2,q^2,q^2,-t^{-2}q^{-2m};q^2)_\infty}\over
{(q^2t^{-2},t^{-2},-q^{-2m};q^2)_\infty}} t^{2m},
\qquad p\in\Zp,\ m\in\Z,
$$
so that all the matrix elements of $Q^2$ with respect to
the orthogonal basis of Proposition \thtag{410} are known.
\enddemo

\proclaim{Corollary \thname{450}} Let $|t|>q^{-1}$
and let $g$ be a bounded continuous function on the spectrum
of $\pi_1(\rit)$ such that $g(\pt(\rit))Q^2$ is of trace class, then
$$
(1-q^2) \text{\rm Tr}\vert_{\H} \bigl(g(\pt(\rit))Q^2\bigr)
= {1\over{t-t^{-1}}} \int_{-t^{-1}}^{\infty(t)} g(x)\, d_{q^2}x
= (1-q^2) h\bigl( g(\pi(\rit))\bigr).
$$
\endproclaim

Here we use the notation for the Jackson $q$-integral, see
\cite{\GaspR, \S 1.11}, $cd>0$, 
$$
\int_{-d}^{\infty(c)} g(x)\,d_qx =
(1-q)d\sum_{p=0}^\infty g(-dq^p)q^p + 
(1-q)c \sum_{p=-\infty}^\infty g(cq^p)q^p.
\tag\eq{451}
$$
Note that any finitely supported $g$ gives 
a trace class operator $g(\pt(\rit))Q^2$, and 
Corollary \thtag{450} implies that it suffices to take $g$
satisfying 
$\int_{-t^{-1}}^{\infty(t)} |g(x)|\, d_{q^2}(x)<\infty$. 

\demo{Proof} Calculate the trace with respect to the
orthogonal basis of Proposition \thtag{410} using Lemma \thtag{420}
for the first equality. Use Remark \thtag{2100} to get the
second equality from the first. 
\qed\enddemo

%%%%%%%%%%%%%%%%%%%%%%%%%%%%%%%%%%%%%%%%%%%%%%%%%%%%%%%%%%%%%%%%%%%%
\subhead \the\sectionnumber.3 Shift operators\endsubhead 
In this subsection we 
consider the unbounded operators $\pt(\al_{\infty,t})$,
$\pt(\be_{\infty,t})$, 
$\pt(\ga_{\infty,t})$ and $\pt(\de_{\infty,t})$ on $\H$. 
By Proposition \thtag{220} these operators are initially
defined on ${\Cal D}(\Z)$, but we see, using \eqtag{3107},
that these operators are defined on
${\Cal D}(Q) = \{ \sum_k c_ke_k\in\H\mid \sum_k
|c_k|^2q^{-2k}<\infty\}$ and that the actions 
are given by the same formulas. 
The commutation relations of Corollary \thtag{3170} 
show that we may expect that these operators act as shift
operators in the basis of eigenvectors of $\pt(\rit)$
of Proposition \thtag{410}. 
The next proposition shows that this is the case. 

\proclaim{Proposition \thname{460}} Let $|t|>1$, then
$v^\th_p(t), w^\th_p(t)\in {\Cal D}(Q)$. Moreover, 
$$
\align
\pt(\al_{\infty,t})\, v^\th_p(t) &=
q^\hf e^{-i\th}t^{-1} {{1-q^{2p}}\over{1-q^2t^{-2}}}
v^\th_{p-1}(tq^{-1}), \quad p\in\Zp, \\ 
\pt(\al_{\infty,t})\, w^\th_p(t) &=q^\hf e^{-i\th}t^{-1}
{{1+t^2q^{2p}}\over{1-q^2t^{-2}}}
w^\th_p(tq^{-1}), \quad p\in\Z,\\
\pt(\be_{\infty,t})\, v^\th_p(t) &= -t^2q^\hf e^{-i\th}(1-t^{-2})
v^\th_p(tq), \quad p\in\Zp,\\
\pt(\be_{\infty,t})\, w^\th_p(t) &= -t^2q^\hf e^{-i\th}(1-t^{-2})
w^\th_{p-1}(tq), \quad p\in\Z,\\
\pt(\ga_{\infty,t})\, v^\th_p(t) &=
-q^\hf e^{i\th} {{1-q^{2+2p}t^{-2}}\over{1-q^2t^{-2}}}
v^\th_p(tq^{-1}), \quad p\in\Zp,\\
\pt(\ga_{\infty,t})\, w^\th_p(t) &=  -q^\hf e^{i\th}
{{1+q^{2+2p}}\over{1-q^2t^{-2}}} w^\th_{p+1}(tq^{-1}), 
\quad p\in\Z,\\
\pt(\de_{\infty,t})\, v^\th_p(t) &= e^{i\th}
tq^\hf (1-t^{-2}) v^\th_{p+1}(tq), \quad p\in\Zp,\\
\pt(\de_{\infty,t})\, w^\th_p(t) &= e^{i\th}
tq^\hf (1-t^{-2}) w^\th_p(tq), \quad p\in\Z.
\endalign
$$
\endproclaim

\demo{Proof} To see that for $|t|>1$ we have 
$v^\th_p(t), w^\th_p(t)\in {\Cal D}(Q)$
we proceed as in the first paragraph of the proof of
Proposition \thtag{420}. Note that the commutation relations
of Corollary \thtag{3170} suggest that the operators
$\pt(\al_{\infty,t})$ and $\pt(\ga_{\infty,t})$,
respectively $\pt(\be_{\infty,t})$ and $\pt(\de_{\infty,t})$,
map eigenvectors of $\pt(\rit)$ into eigenvectors of
$\pt(\rho_{\infty,t/q})$, respectively 
$\pt(\rho_{\infty,tq})$, with a possible $q$-shift in the
eigenvalue. However, this is not sufficient since we 
do not have a priori estimates implying that 
$\pt(\al_{\infty,t})\, v^\th_p(t)\in {\Cal D}
(\pt(\rho_{\infty,t/q}))$. So we have to prove it in a
direct manner.

Let us prove the first two statements. 
{}From Proposition \thtag{220} and \eqtag{3107} we get 
$$
\pt(\al_{\infty,t})\, e_k =
q^\hf \sqrt{1+q^{-2k}}\, e_{k+1} +
q^{\hf -k}t e^{-i\th}\, e_k
$$  
so that we get, using the notation of Proposition \thtag{410},
$$
\pt(\al_{\infty,t}) \sum_{k=-\infty}^\infty 
e^{ik\th} V_k(x;t)e_k =
\sum_{k=-\infty}^\infty q^\hf e^{i(k-1)\th}\lbrace
\sqrt{1+q^{2-2k}} V_{k-1}(x;t) +
q^{-k}t  V_k(x;t)\rbrace e_k,
$$
for $x=-q^{2p}/t$, $p\in\Zp$, or $x=tq^{2p}$, $p\in\Z$.
We use $(z;q)_\infty{}_2\vp_1(a,b;0;q,z)=
(bz;q)_\infty {}_1\vp_1(b;bz;q,az)$, see \cite{\GaspR, (1.4.5)},
to evaluate the term in curly brackets. 
So in terms of a ${}_2\vp_1$-series we have 
$$
V_k(x;t)=(-q^{2-2k};q^2)^\hf_\infty q^{\hf k(k+1)}
(-t)^{-k} {{(-q^2xt^{-1};q^2)_\infty}\over
{(q^2t^{-2};q^2)_\infty}} 
{}_2\vp_1\left( {{{-1\over{xt}},
-q^{2k}}\atop{0}};q^2,-q^2{x\over t}\right).
$$
The contiguous relation
${}_2\vp_1(aq,b;0;q,z)-{}_2\vp_1(a,b;0;q,z)=
az(1-b){}_2\vp_1(aq,bq;0;q,z)$, see \cite{\GaspR, Exer. 1.9(ii)},
gives
$$
\sqrt{1+q^{2-2k}} V_{k-1}(x;t) + q^{-k}t V_k(x;t) =
t^{-1} {{1+xt}\over{1-q^2t^{-2}}} V_k(xq^{-1};tq^{-1})
$$
for $|x|< |t|q^{-2}$, and which is valid for
$x\not=0$ by analytic continuation.  
This gives 
$$
\pt(\al_{\infty,t})\, \sum_{k=-\infty}^\infty e^{ik\th}
V_k(x;t)\, e_k
= q^\hf e^{-i\th}t^{-1} {{1+xt}\over{1-q^2t^{-2}}}
\sum_{k=-\infty}^\infty e^{ik\th}V_k(xq^{-1};tq^{-1})\, e_k.
$$
for $x=-q^{2p}t^{-1}$, $p\in\Zp$ and $x=tq^{2p}$, $p\in\Z$, 
which is the desired result for $\pt(\al_{\infty,t})$. 

The statements for $\pt(\de_{\infty,t})$ follow from 
the ones for $\pt(\al_{\infty,t})$ already proved
and Proposition \thtag{410}.
{}From the factorisation in Corollary \thtag{3170} we get
$\pt(\de_{\infty,tq^{-1}})\pt(\al_{\infty,t})=t\pt(\rit)+1$.
Hence, the result for $\pt(\de_{\infty,t})$ follows for
$|t|>q^{-1}$. The result for $|t|>1$ follows by component wise
analytic continuation in $t$. 

Next we consider $\pt(\ga_{\infty,t})$. This derivation is
completely similar to the one for $\pt(\al_{\infty,t})$, but
now we use 
${}_2\vp_1(a,b;c;q,z)-(1-a){}_2\vp_1(aq,b;c;q,z) =
a {}_2\vp_1(a,b;c;q,qz)$, which is directly verified.
The statements for $\pt(\be_{\infty,t})$ follow from 
the ones for $\pt(\ga_{\infty,t})$ together with 
the factorisation in Corollary \thtag{3170}. 
\qed\enddemo 

%%%%%%%%%%%%%%%%%%%%%%%%%%%%%%%%%%%%%%%%%%%%%%%%%%%%%%%%%%%%%%%%%%%%
\subhead \the\sectionnumber.4 The Haar functional\endsubhead 
In Remark \thtag{2100} we have defined $h(T)\in[0,\infty]$ for
any decomposable bounded operator acting on
$L^2(\T;\H)$ such that $T(e^{i\th})Q^2$ is of trace class.
We want to give an explicit form for the
sesquilinear form $\langle \xi_1,\xi_2\rangle =
h(\pi(\xi_2^\ast\xi_1))$ per bi-$K$-type in the 
Cartan decomposition of Theorem \thtag{3140} for the case 
$s\to\infty$. We have to adjust the definition of the
sesquilinear form in order to apply it to sufficiently decreasing
functions in Theorem \thtag{3140}. For a formal element of the
form $\xi^{(i)}=\Ga^{(i)}_{l,m}(\infty,t)g(\rit)$ for a bounded
continuous function $g$ on the spectrum of 
$\pi_1(\rit)$ we define the corresponding quadratic form by 
$$
\langle \xi^{(i)},\xi^{(i)}\rangle =
h\Bigl( \bar g(\pi(\rit)) \pi\bigl(\Ga^{(i)}_{l,m}(\infty,t)^\ast 
\Ga^{(i)}_{l,m}(\infty,t)\bigr) g(\pi(\rit))\Bigr).
$$
By \eqtag{415} and \eqtag{3138} we regard the operator in 
parentheses as a decomposable operator for suitable functions $g$ 
and we assume it satisfies
the conditions of Remark \thtag{2100}.
For $a,b\geq 0$ and $cd>0$ we put for functions $f$ and $g$
$$
\langle f, g\rangle_{c,d}^{a,b} = 
\int_{-d}^{\infty(c)} f(x) \overline{g(x)} \, 
{{(-q^2x/c,-q^2x/d;q^2)_\infty}\over{ 
(-q^{2+2a}x/c,-q^{2+2b}x/d;q^2)_\infty}}\, d_{q^2}x,
\tag\eq{465}
$$
provided that the $q$-integral is absolutely convergent.

\proclaim{Theorem  \thname{470}} 
Let $g$ be a bounded continuous function on the 
spectrum of $\pi_1(\rit)$, and 
$l\in\hZp$, $m\in\{ -l,-l+1,\ldots,l\}$, $|t|>q^{-1}$. 
Assume that
$$
\bar g(\pt(\rit)) \pt\bigl(\Ga^{(i)}_{l,m}(\infty,t)^\ast 
\Ga^{(i)}_{l,m}(\infty,t)\bigr) g(\pt(\rit))Q^2
$$ 
is of trace class on $\H$. Then
$$
\gather
\langle \xi^{(1)},\xi^{(1)}\rangle = 
C_1 \langle g, g\rangle_{t,t^{-1}}^{l-m,l+m}, \qquad 
\langle \xi^{(2)},\xi^{(2)}\rangle = 
C_2 \langle g, g\rangle_{t,t^{-1}q^{2l-2m}}^{l+m,l-m},\\
\langle \xi^{(3)},\xi^{(3)}\rangle = 
C_3 \langle g, g\rangle_{tq^{2l-2m},t^{-1}}^{l-m,l+m}, \qquad
\langle \xi^{(4)},\xi^{(4)}\rangle = 
C_4 \langle g, g\rangle_{tq^{2l+2m},t^{-1}q^{2l-2m}}^{l+m,l-m}, 
\endgather
$$
for positive constants $C_i$ independent of $g$.
\endproclaim

\demo{Proof}
We have $\pt\bigl(\Ga^{(i)}_{l,m}(\infty,t)^\ast 
\Ga^{(i)}_{l,m}(\infty,t)\bigr) = C_ip_{2l}^{(i)}(\pt(\rit))$
as an unbounded operator defined on ${\Cal D}(\Z)$
by \eqtag{3138} and Proposition \thtag{220} for
explicit polynomials $p^{(i)}_{2l}$ of degree $2l$.
We can now use Proposition \thtag{410} to extend this
operator to an unbounded self-adjoint operator for $|t|>q^{-1}$.
The last statement is a consequence of Corollary \thtag{450}.
\qed\enddemo

\demo{Remark \thname{472}} (i) The sesquilinear form
in \eqtag{465} is positive semi-definite if $q^2d/c<1$.
For $i=1,2$ this condition is satisfied. For $i=3$ we need
$q^{2-2l+2m}<t^2$ and for $i=4$ we need $q^{2-4m}<t^2$ for
positive semi-definiteness of the quadratic form in Theorem
\thtag{470}. To explain this phenomenon
we note that in general
$$
\text{\rm Tr}\vert_{\H} \bigl(\bar g(\pt(\rit))
\pt\bigl(\Ga^{(i)}_{l,m}(\infty,t)^\ast 
\Ga^{(i)}_{l,m}(\infty,t)\bigr) g(\pt(\rit)) Q^2\bigr)
\tag\eq{473}
$$
is not equal to the square of the Hilbert-Schmidt norm
of the operator $S$ with 
$S=\pt(\Ga^{(i)}_{l,m}(\infty,t)) g(\pt(\rit)) Q$.
First of all, the trace is not cyclic for
unbounded operators, cf. Remark \thtag{2100}.
Secondly, it is not true that, under the conditions
of Theorem \thtag{470}, $S$ extends to a bounded operator
on $\H$. Initially, $S$ is only defined on
$\{ v\in {\Cal D}(\Z)\mid g(\pt(\rit))v\in {\Cal D}(\Z)\}$,
which does not even have to be dense in $\H$. However,
if we assume that $S$ and also $T=SQ^{-1}$ are defined on finite
linear combinations of the orthogonal basis of $\H$
given in Proposition \thtag{410}, and then have an extension
to a bounded operator on $\H$, then \eqtag{473} indeed equals
the square of the Hilbert-Schmidt norm of $S$, and
positivity follows. Let us consider the case $i=3$, the other
cases are similar. By Proposition \thtag{410}, Proposition
\thtag{460} and \eqtag{3135} for $s=\infty$, we see that we
can calculate $Tv^\th_p(t)$ and $Tw^\th_p(t)$
explicitly for $|t|>q^{1-2l}$. From this we can determine
conditions on $g$ that imply that $T$ extends to a bounded
operator on $\H$, e.g. it suffices to consider compactly
supported $g$ with $0\not\in\text{supp}(g)$. For such
$g$ we can also extend $S$ to a bounded operator on $\H$,
since we can estimate the growth of
$\| Qv^\th_p(t)\|$ and $\| Qw^\th_p(t)\|$ from
Lemma \thtag{420}. Note that $|t|>q^{1-2l}$ implies
$|t|^2>q^{2-2l+2m}$, since $|m|\leq l$.

(ii) In case $q^2d/c<1$, $\langle\cdot,\cdot\rangle_{c,d}^{a,b}$
gives an inner product. The corresponding Hilbert space is a
weighted $L^2$-space, on which the big $q$-Jacobi function
transform lives, see \cite{\KoelSbig} and \S 6 for a quantum group
theoretic interpretation.
\enddemo

\demo{Remark \thname{475}} It would be desirable 
to interpret  
the elements $\pi(\Ga^{(i)}_{l,m}(\infty,t))g(\pi(\rit))$
in Theorem \thtag{470} in terms of affiliated elements 
for the $\text{C}^\ast$-algebra 
$\pi(C_0(X)\times_\tau\Z)$ or more generally, in terms of regular
operators for Hilbert $\text{C}^\ast$-modules,   
see Woronowicz \cite{\WoroEtwee} for the notion
of affiliated elements and \cite{\Lanc},  
\cite{\Kust} for regular operators. This would give rise
to an interpretation of $\Ga^{(i)}_{l,m}(\infty,t)g(\rit)$
as a uniquely defined element affiliated to the 
$\text{C}^\ast$-algebra 
$C_0(X)\times_\tau\Z$, see Kustermans \cite{\Kust}. 
In general this seems not to be
possible due to the fact that the density requirements 
in either the
definition of affiliated element in \cite{\WoroEtwee} or in the
definition of regular operator in \cite{\Lanc} is not met. 
Let $f\in C_c(\Z;C(X))$ such that it is supported in precisely one
point. It is
straightforward to check that multiplication by $\pi(f)$ is 
a regular operator of the $\text{C}^\ast$-algebra
$\pi(C_0(X)\times_\tau\Z)$ viewed as a 
Hilbert $\text{C}^\ast$-module
over itself, see Lance \cite{\Lanc, Ch.~9} and Woronowicz
\cite{\WoroEtwee, \S 3.C}. 
However, it is not clear if this remains true
for $f\in C_c(\Z;C(X))$
supported in more that just one
point, such as for $f$ corresponding to $\pi(\rit)$. 

Moreover, it is also unclear that for functions 
$g\in C_0(\R)$ the operator of multiplication
by $g(\pi(\rit))$ is a multiplier of the $\text{C}^\ast$-algebra
$\pi(C_0(X)\times_\tau\Z)$. The problem is that it is not clear if
this multiplication operator preserves 
the $\text{C}^\ast$-algebra $\pi(C_0(X)\times_\tau\Z)$. 
However, in this particular case $s=\infty$, estimates
can be obtained that show that $g(\pi(\rit))\in
\pi(C_0(X)\times_\tau\Z)$ for $g$ finitely supported on the
spectrum. 
\enddemo

%%%%%%%%%%%%%%%%%%%%%%%%%%%%%%%%%%%%%%%%%%%%%%%%%%%%%%%%%%%%%%%%%%%%  
%N E W   S E C T I O N%
%%%%%%%%%%%%%%%%%%%%%%%%%%%%%%%%%%%%%%%%%%%%%%%%%%%%%%%%%%%%%%%%%%%%

\head \newsec The Haar functional on the
algebra generated by $\rst$\endhead

In this section we calculate explicitly the Haar functional
related to the Cartan decomposition of
Theorem \thtag{3140}. The result is given in terms of a
non-compact analogue of the Askey-Wilson measure,
and it is obtained using the spectral analysis of
$\pt(\rst)$ and \eqtag{2105}.  This operator is
considered in two invariant complementary subspaces
$V^\th(t)$ and $W^\th(t)$ of $\H$. The spectral 
decomposition of $\pt(\rst)$ on $V^\th(t)$ is
obtained using orthogonal polynomials and is
analogous to \cite{\KoelV, \S 6}. On $W^\th(t)$
the spectral analysis is related to the 
little $q$-Jacobi function transform. Matching
these two results involves non-trivial summation
formulas for basic hypergeometric series. The 
main new summation formula has been proved by Mizan
Rahman, and its proof is presented in Appendix B. 
Recall the basic assumption that 
$s$ and $t$ are real parameters, and we also
assume that $|s|, |t|>q^{-1}$.

%%%%%%%%%%%%%%%%%%%%%%%%%%%%%%%%%%%%%%%%%%%%%%%%%%%%%%%%%%%%%%
\subhead \the\sectionnumber.1 Spectral
analysis of $\pt(\rst)\vert_{V^\th(t)}$\endsubhead 
In this subsection we calculate the spectral measure
for $\pt(\rst)\vert_{V^\th(t)}$, which is a bounded
operator that can be viewed as a Jacobi matrix. This
enables us to link it to the Al-Salam and Chihara
polynomials. 
The analysis in this part follows \cite{\KoelV, \S 6}.

The operator $\pt(\rst)$ is an unbounded five-term recurrence 
operator in the standard basis $\{e_k\mid k\in\Z\}$ of $\H$ 
densely defined on ${\Cal D}(\Z)$ by Proposition \thtag{220}
and \eqtag{3100}. We can extend the domain of $\pt(\rst)$
to ${\Cal D}(Q^2)$, since $\rst$ consists of quadratic elements
in the generators $\al$ and $\ga$. Since 
$v_p^\th(t)\in {\Cal D}(Q^2)$ for $|t|>q^{-1}$ we see that 
$\pt(\rst)v_p^\th(t)$ is well-defined. 
It follows from \eqtag{3160} and 
Proposition \thtag{460}
that $\pt(\rst)$ is a three-term recurrence operator in 
the basis $v^\th_p(t)$, $p\in\Zp$; 
$$
\multline
2\pt(\rst)v^\th_p(t) = 
- qe^{2i\th}(1-q^{2+2p}t^{-2})v^\th_{p+1}(t) + \\ 
q^{1+2p}t^{-1}(s+s^{-1}) v^\th_p(t) 
-q^{-1}e^{-2i\th}(1-q^{2p})v^\th_{p-1}(t).
\endmultline
\tag\eq{505}
$$
Note that $\pt(\rst)\vert_{V^\th(t)}$ is a bounded operator. 
By going over to the orthonormal basis 
$f_p=(-e^{2i\th})^p v^\th_p(t)/\|v^\th_p(t)\|$, $p\in\Zp$, see
Proposition \thtag{410}, we obtain 
$$
\gathered
2\pt(\rst)f_p = a_{p+1}f_{p+1}+b_pf_p + a_p f_{p-1}, 
\qquad p\in\Zp, \\
a_p= \sqrt{(1-q^{2p}t^{-2})(1-q^{2p})}, \qquad b_p
=q^{1+2p}t^{-1}(s+s^{-1}),
\endgathered
\tag\eq{510} 
$$
which is, by Favard's Theorem, a three-term recurrence for
orthonormal polynomials since $|t|>1$.  Note that $a_0=0$,
so that \eqtag{510} is a well-defined operator. 
The spectral measure can be
determined completely in terms of the orthogonality measure of
the corresponding orthonormal polynomials, see
e.g. \cite{\Akhi}, \cite{\KoelVdJ}, \cite{\KoelV}, \cite{\Simo}.
The polynomials can be identified
with the Al-Salam and Chihara polynomials.

We recall that the Al-Salam and Chihara polynomials,
originally introduced by Al-Salam and Chihara 
in \cite{\AlSaC}, are orthogonal polynomials with respect  to an
absolutely continuous measure on $[-1,1]$ plus a 
finite set, possibly empty, of discrete mass points as
established by Askey and Ismail \cite{\AskeI}. The Al-Salam and
Chihara polynomials are a subclass of the Askey-Wilson
polynomials by setting two parameters of the four parameters of
the Askey-Wilson polynomials equal to zero, see
Askey and Wilson \cite{\AskeW}, or \cite{\GaspR, Ch.~7}. 

The Al-Salam and Chihara polynomials are defined by
$$
s_m(\cos\psi;a,b|q) = a^{-m} (ab;q)_m\,
{}_3\vp_2 \left( {{q^{-m},ae^{i\psi},ae^{-i\psi}}\atop{ab,\ 0}};
q,q\right).
\tag\eq{515}
$$
Let $S_m(x;a,b|q)=s_m(x;a,b|q)/\sqrt{(q,ab;q)_m}$ denote the
orthonormal Al-Salam and Chihara polynomials, which satisfy the
three-term recurrence relation
$$
\gathered
2x\, S_n(x) = a_{n+1}\, S_{n+1}(x) + q^n(a+b)\, S_n(x) +
a_n\, S_{n-1}(x),\\ a_n = \sqrt{ (1-abq^{n-1})(1-q^n)}.
\endgathered
\tag\eq{520}
$$
We assume $ab<1$, so that $a_n>0$. 
Since the coefficients $a_n$ and $b_n$ are bounded, the
corresponding moment problem is determined and the 
orthonormal Al-Salam and Chihara polynomials
form an orthonormal basis of $L^2(\R, dm(\cdot;a,b|q^2))$, with
$dm(\cdot;a,b|q^2))$ the normalised orthogonality measure.  The
explicit form of the orthogonality measure is 
originally obtained by Askey and Ismail \cite{\AskeI},
and it is a special case of the Askey-Wilson measure. 
Since the Askey-Wilson measure is needed in the
next subsection we recall it here, see 
\cite{\AskeW, \S 2}, \cite{\GaspR, Ch.~6}: 
$$
\int_{\R} p(x)dm(x;a,b,c,d| q) =
{1\over {h_0 2\pi}}\int_{0}^{\pi}
p(\cos\theta)w(e^{i\th})\,d\theta +
{1\over{h_0}} \sum_{k}p(x_k)w_k.
\tag\eq{530}
$$
Here we use the notation $w(z)=w(z;a,b,c,d|q)$,
$h_0=h_0(a,b,c,d|q)$ and
$$
\aligned
h_0(a,b,c,d|q) & =
{{(abcd;q)_\infty}\over{(q,ab,ac,ad,bc,bd,cd;q)_\infty}}, \\
w(z;a,b,c,d|q) & = {{(z^2,z^{-2};q)_\infty} \over
{(az,a/z, bz,b/z,
cz,c/z, dz,d/z;q)_\infty}},
\endaligned
\tag\eq{540}
$$
and we suppose that $a$, $b$, $c$ and $d$ are real 
or $\bar a=b$, and $c,d\in\R$ and such that all
pairwise products are less than $1$.
The sum in \eqtag{530} is over the points $x_k$
of the form $\mu(eq^k)=(eq^k + e^{-1}q^{-k})/2$ with $e$ any of
the parameters $a$, $b$, $c$ or $d$
whose absolute value is larger than one and such that
$|eq^k|>1$, $k\in\Zp$. The corresponding mass
$w_k$ is the residue of $z\mapsto w(z)/z$ at 
$z = eq^k$. The value of $w_k$ in case $e=a$ is
given in \cite{\AskeW, (2.10)}, \cite{\GaspR, (6.6.12)}. 
Explicitly,
$$
\multline
w_k(a;b,c,d|q) =
{{(a^{-2};q)_\infty}\over{(q,ab,b/a,ac,c/a,ad,d/a;q)_\infty}} \\
\times {{(1-a^2q^{2k})}\over{(1-a^2)}} {{(a^2,ab,ac,ad;q)_k}\over
{(q,aq/b,aq/c,aq/d;q)_k}} \Bigl( {q\over{abcd}}\Bigr)^k.
\endmultline
\tag\eq{545}
$$
The orthogonality measure for the Al-Salam and Chihara
polynomials is obtained by taking $c=d=0$ 
in \eqtag{530}, so $dm(\cdot;a,b|q)=
dm(\cdot;a,b,0,0|q)$. 

Now compare \eqtag{510} with \eqtag{520} in base $q^2$ with $a$
and $b$ replaced by $qst^{-1}$ and $qs^{-1}t^{-1}$. This shows
that we can realise $\pt(\rst)\vert_{V^\th(t)}$ 
as a multiplication operator on the weighted $L^2$-space 
corresponding to the orthogonality measure 
$dm(\cdot;qst^{-1},qs^{-1}t^{-1}|q^2)$ 
using the unitary isomorphism $V^\th(t)\to L^2(\R, 
dm(\cdot;qst^{-1},qs^{-1}t^{-1}|q^2))$ mapping
$f_p$ to the corresponding $p$-th orthonormal
polynomial $S_p(\cdot;qst^{-1},qs^{-1}t^{-1}|q^2)$, 
see e.g. Akhiezer 
\cite{\Akhi, Ch.~1}, Simon \cite{\Simo}. This
proves the following proposition.   

\proclaim{Proposition \thname{550}} Let $s,t\in\R$ with 
$|s|>1$, $|t|>q^{-1}$. 
The spectrum of the bounded self-adjoint operator
$\pt(\rst)\vert_{V^\th(t)}$ consists of the continuous
spectrum $[-1,1]$ and the finite discrete spectrum, 
possibly empty,
$\{ q^{1+2k}st^{-1}\mid |q^{1+2k}st^{-1}|>1, \, k\in\Zp\}$.
Explicitly, with $f_p=(-e^{2i\th})^p v^\th_p(t)/\|v^\th_p(t)\|$,
$$
\langle \pt(\rst) f_n,f_m\rangle 
= \int_\R x \bigl(S_nS_m\bigr)(x; qst^{-1},qs^{-1}t^{-1}|q^2)
dm(x;qst^{-1},qs^{-1}t^{-1}|q^2).
$$
\endproclaim

\proclaim{Proposition \thname{560}} Let $s,t\in\R$ with 
$|s|,|t|>q^{-1}$ and let $P$ be the orthogonal projection
onto $V^\th(t)$ along the decomposition
$\H=V^\th(t)\oplus W^\th(t)$. Then $PQ^2\vert_{V^\th(t)}\colon
V^\th(t)\to V^\th(t)$ is bounded. Let $f$ be a
continuous function on the spectrum of 
$\pt(\rst)\vert_{V^\th(t)}$, and assume that
$f\bigl(\pt(\rst)\vert_{V^\th(t)}\bigr)PQ^2$ is of trace class
on $V^\th(t)$. Then its trace is integrable over $[0,2\pi]$ as
function of $\th$ and 
$$
\multline
{1\over{2\pi}} \int_0^{2\pi}
\text{\rm Tr}\vert_{V^\th(t)}(f(\pt(\rst)\vert_{V^\th(t)})
PQ^2)\, d\th =
\\ {1\over{2\pi}} \int_0^\pi f(\cos\th)
{{(1-q^2/t^2)(1-e^{\pm 2i\th})}\over{
(t^2-1)(1-{{qs}\over t}e^{\pm i\th})
(1-{q\over{st}}e^{\pm i\th})}}
\, {}_8W_7({{q^2}\over{t^2}};q^2,{{qs}\over t}e^{\pm i\th},
{q\over{st}}e^{\pm i\th};q^2,q^2)\, d\th \\
+ \sum_{k\in\Zp,\, |q^{1+2k}s/t|>1} 
{{w_k(qs/t;q/st,qt/s,qst|q^2)}\over{h_0(qs/t,q/st,qt/s,qst|q^2)}}
{-1\over{1-q^2}} f\bigl( \mu(q^{1+2k}s/t)\bigr), 
\endmultline
$$
where $\mu(z)=\hf(z+z^{-1})$. 
The $\pm$-sign means that we have to take all
possibilities. 
\endproclaim

The positivity of the weight for the discrete mass points 
in Proposition \thtag{560} follows from 
$$
{-1\over{1-q^2}}
{{w_k(qs/t;q/st,qt/s,qst|q^2)}\over{h_0(qs/t,q/st,qt/s,qst|q^2)}}
= {{(q^{2+4k}s^2-t^2)q^{-2-2k}}\over{(s^2-1)(t^2-1)}},
$$
using \eqtag{540}, \eqtag{545}. For $|q^{1+2k}s/t|>1$ this is
positive. 

Note that in the
${}_8W_7$-series, see \S 1, 
a lot of cancellation occurs, 
$$
\multline
{}_8W_7({{q^2}\over{t^2}};q^2,{{qs}\over t}e^{i\th},{{qs}\over
t}e^{-i\th}, {q\over{st}}e^{i\th}, 
{q\over{st}}e^{-i\th};q^2,q^2) = \\ 
\sum_{k=0}^\infty {{1-t^{-2}q^{2+4k}}\over{1-t^{-2}q^2}}
{{(1-{{qs}\over t}e^{i\th})(1-{{qs}\over t}e^{-i\th})
(1-{{q}\over {st}}e^{i\th})(1-{{q}\over {st}}e^{-i\th})}\over
{(1-{{q^{2k+1}s}\over t}e^{i\th})(1-{{q^{2k+1}s}\over
t}e^{-i\th}) (1-{{q^{2k+1}}\over {st}}e^{i\th})(1-{{q^{2k+1}}
\over {st}}e^{-i\th})}} q^{2k}.
\endmultline
$$
This can be used to integrate the function $f=1$ explicitly
over the interval $[0,\pi]$ in Proposition \thtag{560}.  

\demo{Proof of Proposition \thtag{560}} 
Since \eqtag{425} implies 
$$
\big\vert{{\langle Q^2 v^\th_n(t),v^\th_m(t)\rangle}\over{
{\|v^\th_n\| \|v^\th_m\|}}}\big\vert \leq Cq^{m+n}, 
\qquad n,m\in\Zp, 
$$
using that $Q^2$ is self-adjoint and that
$v^\th_n(t)\in{\Cal D}(Q^2)$ for
$|t|>q^{-1}$, it follows that $PQ^2\vert_{V^\th(t)}$ is bounded. 
The rest of the 
proof of Proposition \thtag{560} is completely analogous
to the proof of \cite{\KoelV, Prop.~6.3},
see also the proof of Proposition~5.4. %harde referentie
Use Proposition \thtag{550} and \eqtag{425} to calculate the
trace formally as a double sum involving an integral of
a product of $f$ and two Al-Salam and Chihara polynomials. Since
the $\th$-dependence in 
$\text{\rm Tr}\vert_{V^\th(t)}(f(\pt(\rst))Q^2)$ 
is easy, integration over $[0,2\pi]$ reduces to a single
sum. Interchanging summation and integration gives an integral
involving the Poisson kernel for the Al-Salam and Chihara
polynomials. This can be justified using 
the same estimate as in the beginning of the proof 
and the asymptotics for the
Al-Salam and Chihara polynomials, see \cite{\AskeI, \S 3.3}. The
Poisson kernel for the
Al-Salam and Chihara polynomials is given in terms of
a very-well-poised ${}_8\vp_7$-series by
Askey, Rahman and Suslov \cite{\AskeRS, (14.8)}, see
also \cite{\IsmaS, \S 4} and \cite{\VandJ} for other
derivations. 
The very-well-poised ${}_8\vp_7$-series is
summable for points in the discrete spectrum.  See 
\cite{\KoelV, \S 6} for details. 
\qed\enddemo

%%%%%%%%%%%%%%%%%%%%%%%%%%%%%%%%%%%%%%%%%%%%%%%%%%%%%%%%%%%%%%%%%
\subhead \the\sectionnumber.2 Spectral
analysis of $\pt(\rst)\vert_{W^\th(t)}$\endsubhead 
In this subsection we calculate the spectral measure
for $\pt(\rst)\vert_{W^\th(t)}$, which is an unbounded
operator that can be viewed as doubly infinite 
Jacobi matrix. The operator has been studied by 
Kakehi \cite{\Kake}, Kakehi, Masuda and Ueno \cite{\KakeMU}
in connection with the spherical Fourier transform on the
quantum $SU(1,1)$, i.e. 
corresponding to the Cartan decomposition of 
Theorem \thtag{3140} for the case $(s,t)=(\infty,\infty)$. 
This is the little $q$-Jacobi function transform, and 
it is discussed in Appendix A. 

As in the previous subsection, cf. \eqtag{505}, we can apply 
$\pt(\rst)$ to $w_p^\th(t)$ for $|t|>q^{-1}$.  
{}From \eqtag{3160} and Proposition \thtag{460}
we see that $\pt(\rst)$ is a three-term recurrence operator in 
the basis $w^\th_p(t)$, $p\in\Z$, of $W^\th(t)$; 
$$
\multline
2\pt(\rst)w^\th_p(t) = 
- qe^{2i\th}(1+q^{2+2p})w^\th_{p+1}(t)  \\ 
- q^{1+2p}t(s+s^{-1}) w^\th_p(t)
- q^{-1}e^{-2i\th}(1+t^2q^{2p})w^\th_{p-1}(t).
\endmultline
\tag\eq{565}
$$
By going over to the orthonormal basis 
$f_{-p}=(-e^{2i\th})^p w^\th_p(t)/\|w^\th_p(t)\|$, 
$p\in\Z$, see
Proposition \thtag{410}, we obtain 
$$
\gathered
2\pt(\rst)f_p = a_{p+1}f_{p+1}+b_pf_p + a_p f_{p-1}, 
\qquad p\in\Z, \\
a_p= \sqrt{(1+q^{2-2p}t^2)(1+q^{2-2p})}, \qquad b_p
= - q^{1-2p}t(s+s^{-1}).
\endgathered
\tag\eq{570} 
$$
This is an unbounded symmetric operator that has been studied in
\cite{\Kake}, \cite{\KakeMU}, see also Appendix A, 
Theorem~A.5. %harde referentie
So the spectral measure of the
operator $\pt(\rst)\vert_{W^\th(t)}$ is determined
in terms of little $q$-Jacobi functions. 

Put
$$
\phi_n(x;s,t|q^2) = {}_2\vp_1\left( {{qs^{-1}t^{-1}z,
qs^{-1}t^{-1}z^{-1}}\atop{q^2s^{-2}}};q^2, -q^{2n}\right), \qquad
x=\mu(z)=\hf(z+z^{-1}), 
$$
for the little $q$-Jacobi function adapted to our 
situation, see Appendix A for its definition in case $n\leq 0$.
Put, cf. \thetag{A.13}, %harde referentie 
$$
w_n = (qs^{-1}t^{-1})^n \sqrt{{{(-q^{2-2n}t^2;q^2)_\infty}
\over{(-q^{2-2n};q^2)_\infty}}},
\tag\eq{580}
$$
so that the little $q$-Jacobi function transform is given
by
$$
({\Cal G}u)(x) = \sum_{n=-\infty}^\infty
w_n \phi_n(x;s,t|q^2) u_n, \qquad u=\sum_{n=-\infty}^\infty 
u_n f_n \in W^\th(t),
$$
initially defined for finite sums and extended to $W^\th(t)$ by
continuity, see Appendix A. 
In order to describe the spectral measure we introduce 
the following measure 
$$ 
\multline
\int_\R f(x)  \, d\nu(x;a,b;d|q)= 
h_0(a,b,q/d,d|q)\int_\R f(x) \, dm(x;a,b,q/d,d|q) \\
+ \sum_{ k\in\N} 
f\bigl(\mu(dq^{-k})\bigr) 
\text{Res}_{z=dq^{-k}} {{w(z;a,b,q/d,d|q)}\over z}
\endmultline
\tag\eq{590}
$$
using the notation of \eqtag{540}. 
Observe that
$$
\text{Res}_{z=dq^{-k}} {{w(z;a,b,q/d,d|q)}\over z} = 
{{-d^{2(k-1)}q^{-k(k-1)}(1-d^2q^{-2k})}\over{
(q,q,adq^{-k},bdq^{-k},aq^k/d,bq^k/d;q)_\infty}}
\tag\eq{595}
$$
by a straightforward calculation, 
wich equals $- w_{k-1}({q\over d};a,b,d|q)$ using \eqtag{545}.  
It follows that this
measure is supported on $[-1,1]$ plus a finite, possibly
empty, set of discrete mass points of the form
$\{ \mu(eq^k)\mid k\in\Zp, |eq^k|>1\}$, where $e$ is
$a$ or $b$, plus an infinite set of discrete mass points
of the form $\{ \mu(dq^{-k})\mid k\in\Z, |dq^{-k}|>1\}$. 

\demo{Remark}
The measure $d\nu(x;a,b;d|q)$ 
defined in \eqtag{590} is positive 
for $ab<1$, $ad<0$, $bd<0$ or for
$a$, $b$ in complex conjugate pair with $ab<1$,  
and has unbounded support. 
The measure in \eqtag{590}
can be obtained from the standard Askey-Wilson measure
by a limiting procedure; consider the Askey-Wilson measure with
parameters $a$, $b$, $cq^l$ and $dq^{-l}$, and let
$l\to\infty$. The parameter $c$ disappears in the limit. 
Then we formally obtain the measure 
of \eqtag{590}, and in this way we formally obtain the
little $q$-Jacobi function transform as a limit case of the 
orthogonality relations for the Askey-Wilson polynomials.  In the
corresponding quantum group theoretic setting this
corresponds to the limit transition of the compact 
quantum $SU(2)$ group to the quantum $E(2)$ group of
orientation and distance preserving motions of the Euclidean 
plane. 
In that case, \eqtag{590} gives an expression for
the Haar functional on a certain subalgebra, and
$\Z$ labels the representations of the quantum
group of plane motions, see \cite{\KoelPhD, Ch.~3}. 
\enddemo 

The results of \cite{\Kake} on the little $q$-Jacobi
function transform imply the following proposition,  
see also Theorem~A.5, case (3). %harde referentie 

\proclaim{Proposition \thname{5100}} Let $s,t\in\R$ with 
$|s|,|t|>q^{-1}$. The
operator $\pt(\rst)\vert_{W^\th(t)}$ is essentially
self-adjoint and its spectral decomposition 
is given by 
$$
\langle \pt(\rst)f_n,f_m\rangle = C 
\int_\R x w_nw_m \bigl(\phi_n\phi_m\bigr)(x;s,t|q^2) \, 
d\nu(x;q/st,qt/s;-qst|q^2)
$$
with $C=(q^2s^{-2},-1,-q^2;q^2)_\infty^2$.
The support of 
$d\nu(\cdot;q/st,qt/s;-qst|q^2)$ is the spectrum 
of  $\pt(\rst)\vert_{W^\th(t)}$.
\endproclaim

\proclaim{Proposition \thname{5110}} Let $|s|,|t|>q^{-1}$,  
and assume $s^2t^{\pm 2}\not\in q^{2\Z}$. Let $f$ be a
continuous, compactly supported function 
on the spectrum of $\pt(\rst)\vert_{W^\th(t)}$, 
integrable with respect to 
the measure $d\nu(\cdot;qs^{-1}t^{-1},qts^{-1};-qst|q^2)$, 
and such that $f\bigl(\pt(\rst)\vert_{W^\th(t)}\bigr)
(1-P)Q^2$, with $P$ as in Proposition \thtag{560}, 
is of trace class on $W^\th(t)$. Then 
$$
\th\mapsto \text{\rm Tr}\vert_{W^\th(t)}
(f(\pt(\rst)\vert_{W^\th(t)})(1-P)Q^2)
$$ 
is integrable over $[0,2\pi]$ and 
$$
\multline
{1\over{2\pi}}\int_0^{2\pi}
\text{\rm Tr}\vert_{W^\th(t)}
(f(\pt(\rst)\vert_{W^\th(t)})(1-P)Q^2) d\th =
\\ {{(q^2s^{-2},-1,-q^2;q^2)_\infty^2}\over{1-t^{-2}}} 
\int_\R f(x) R_{q^{-2}}(x;s,t|q^2)
\, d\nu(x;qs^{-1}t^{-1},qts^{-1};-qst|q^2)
\endmultline
$$
where
$R_u(x;s,t|q^2) = \sum_{n=-\infty}^\infty 
u^n  w_n^2 | \phi_n(x;s,t|q^2)|^2$
is absolutely convergent for $u=q^{-2}$, uniformly 
for $x$ in compacta of the support of 
$d\nu(\cdot;qs^{-1}t^{-1},qts^{-1};-qst|q^2)$. 
\endproclaim

\demo{Proof} Using 
$f_{-p}=(-e^{2i\th})^p w^\th_p(t)/\|w^\th_p(t)\|$
and the spectral decomposition of 
Proposition \thtag{5100} we calculate the trace as
$$
\multline
\text{Tr}\vert_{W^\th(t)} (f(\pt(\rst))Q^2) 
= (q^2s^{-2},-1,-q^2;q^2)_\infty^2
\sum_{n,m=-\infty}^\infty (-e^{2i\th})^{(m-n)} \\ 
\times 
{{\langle Q^2\, w^\th_{-n}(t), w^\th_{-m}(t)\rangle}
\over{\|w^\th_{-n}(t)\| \|w^\th_{-m}(t)\|}} w_nw_m 
\int_\R f(x) \phi_n(x)\phi_m(x)
\, d\nu(x;qs^{-1}t^{-1},qts^{-1};-qst|q^2).
\endmultline
\tag\eq{5120}
$$
Note that by symmetry in $n$ and $m$, since $Q^2$ is self-adjoint
and $w^\th_n(t)\in {\Cal D}(Q^2)$ for $|t|>q^{-1}$, we
may restrict to $n\leq m$. We estimate the double sum 
$$
\sum_{n=-\infty}^\infty \sum_{m=n}^\infty \big\vert
{{\langle Q^2\, w^\th_{-n}(t), w^\th_{-m}(t)\rangle}
\over{\|w^\th_{-n}(t)\| \|w^\th_{-m}(t)\|}} w_nw_m 
\phi_n(x)\phi_m(x) \big\vert
$$
and this suffices for most points of the spectrum. 
The weight function is only needed for the
cases $x=\pm 1$.

Using \eqtag{445} we obtain
$$
\Bigl\vert {{\langle Q^2\, w^\th_{-n}(t), w^\th_{-m}(t)\rangle}
\over{\|w^\th_{-n}(t)\| \|w^\th_{-m}(t)\|}} w_nw_m \Bigr\vert=
{{|st|^{-n-m}}\over{1-t^{-2}}} 
{{(-t^2q^{2-2n};q^2)_\infty}\over{(-q^{2-2n};q^2)_\infty}} \leq
\cases C|ts|^{-n-m}, &\text{$n\leq 0$,}\\
C|t/s|^n|st|^{-m}, &\text{$n\geq 0$,}
\endcases
$$
using \eqtag{580} and \eqtag{237}. 
{}From the definition of $\phi_n(\cdot;s,t|q^2)$ it is immediate
that, for $x$ in compact subsets of $\text{supp}(d\nu)$,
$|\phi_n(\cdot;s,t|q^2)|$ is uniformly bounded for $n\geq 0$.
Hence, the $\sum_{n=0}^\infty \sum_{m=n}^\infty$ part of the 
double sum can be majorised by 
$C\sum_{n=0}^\infty \sum_{m=n}^\infty |st|^{-m}|t/s|^n<\infty$,
uniformly for $x$ in compact subsets of $\text{supp}(d\nu)$. 

It remains to consider $\sum_{n=-\infty}^0 \sum_{m=n}^\infty$.
For this we have to estimate $\phi_n(x;s,t|q^2)$
for $n\leq 0$ for $x$ in the support of
the measure. 
Using the $c$-function expansion, see \cite{\GaspR, (4.3.2)}
or \thetag{A.10}, %harde referentie
or e.g. \cite{\Kake}, \cite{\KakeMU}, \cite{\KoelSbig}, we find 
$$
\aligned
\phi_n(\hf(z+z^{-1});s,t|q^2) &=  c(z) \Phi_n(z;s,t;q^2) + 
c(z^{-1}) \Phi_n(z^{-1};s,t;q^2), \\
c(z) &= {{(qt/sz,q/stz, -qz/st,
-qst/z;q^2)_\infty}\over{(z^{-2},q^2/s^2,-1,
-q^2;q^2)_\infty}},\\
\Phi_n(z;s,t|q^2) &= \left( {{qz}\over{st}}\right)^{-n}
\, {}_2\vp_1\left( {{qz/st,qzs/t}\atop{q^2z^2}};
q^2, -t^2 q^{2-2n}\right)
\endaligned
\tag\eq{5140}
$$
for $z^2\notin q^{2\Z}$. See also \S A.2 for the 
definition of $\Phi_n(z;s,t|q^2)$ in case $|t^2 q^{2-2n}|\geq 1$. 
For $n$ sufficiently negative 
we can estimate $\Phi_n(e^{i\psi};s,t|q^2)$
by $|st/q|^n$ times a constant for $0\leq \psi\leq\pi$ by 
continuity. 
Hence, for $n,m\leq 0$ we get, using that 
$(e^{2i\psi},e^{-2i\psi};q^2)_\infty$ is part of the weight
function of the measure, 
$$
|(e^{2i\psi},e^{-2i\psi};q^2)_\infty
\phi_n(\cos\psi;s,t|q^2)\phi_m(\cos\psi;s,t|q^2)
|\leq C |st/q|^{n+m}.
$$
So, on the interval $[-1,1]$ the sum 
$\sum_{n=-\infty}^0 \sum_{m=n}^\infty$, after multiplication
by $(e^{\pm 2i\psi};q^2)_\infty$, is estimated by
$$
\sum_{n=-\infty}^0 q^{-n} 
\bigl( \sum_{m=n}^0 q^{-m}  + 
\sum_{m=0}^\infty |st|^{-m}\bigr)<\infty.
$$
This deals with the convergence of \eqtag{5120}  
on the absolutely continuous part. 

For the discrete part we observe that for
$z_k=tq^{1+2k}/s$, $k\in\Zp$, $|z_k|>1$, or for 
$z_k=-stq^{-1-2k}$, $k\in\Z$, $|z_k|>1$, we see that
$c(z_k)=0$ and by \eqtag{5140} 
$|\phi_n(\mu(z_k))|\leq C
|q/stz_k|^{-n}$ for $n\leq 0$. This estimate then shows
that the double sum is absolutely convergent, and uniform
for $x$ in compact subsets of the support of the 
measure. Note that 
we have used $s^2/t^2, s^2t^2 \notin q^{2\Z}$ to avoid
zeroes in the denominator of the $c$-function of \eqtag{5140}
at $z=z_k$. 

The above proves that the double sum in \eqtag{5120} 
is absolutely convergent. Since the estimates are
uniformly in $\theta$, we see that 
$\th\mapsto \text{\rm Tr}\vert_{W^\th(t)}
(f(\pt(\rst)\vert_{W^\th(t)})(1-P)Q^2)$ is
integrable over $[0,2\pi]$ and moreover that we may
integrate term by term in \eqtag{5120} and interchange
summation. This gives using Lemma \thtag{420}, 
$$
\multline
{1\over{2\pi}}
\int_0^{2\pi}\text{Tr}\vert_{W^\th(t)}
(f(\pt(\rst)\vert_{W^\th(t)})(1-P)Q^2) \, d\th =
{{(q^2s^{-2},-1,-q^2;q^2)_\infty^2}\over{(1-t^{-2})}} \\ \times 
\int_\R f(x) \sum_{n=-\infty}^\infty  q^{-2n}
w_n^2 |\phi_n(x;s,t|q^2)|^2 
\, d\nu(x;qs^{-1}t^{-1},qts^{-1};-qst|q^2)  
\endmultline
$$
which is the expression stated. 

Using the explicit form of $R_u(x;s,t|q^2)$ and the estimates
already in use on the little $q$-Jacobi functions we
immediately obtain that the sum in $R_{q^{-2}}$ is
absolutely convergent both for $x$
in the absolutely continuous part and for $x$ in the discrete
part of the measure
$d\nu(\cdot;qs^{-1}t^{-1},qts^{-1};-qst|q^2)$, and even uniformly
for $x$ in compact subsets of the support of the measure.  
\qed\enddemo 

%%%%%%%%%%%%%%%%%%%%%%%%%%%%%%%%%%%%%%%%%%%%%%%%%%%%%%%%%%%%%%%%%
\subhead \the\sectionnumber.3 The trace of 
$f(\pt(\rst))Q^2$\endsubhead 
We calculate the trace of $f(\pt(\rst))Q^2$ in this subsection 
and we integrate the result over $[0,2\pi]$. This gives the 
Haar functional on $f(\rst)$, see
Remark \thtag{2100}, as an explicit
Askey-Wilson type integral with unbounded support
of the form \eqtag{590}. 

Since $f(\pt(\rst))$ preserves the decomposition
$\H=V^\th(t)\oplus W^\th(t)$ arising from Proposition 
\thtag{410} we have that 
$$
\text{Tr}\vert_{\H }f(\pt(\rst))Q^2 =  
\text{Tr}\vert_{V^\th(t)} f(\pt(\rst))Q^2
+ \text{Tr}\vert_{W^\th(t)} f(\pt(\rst))Q^2 
$$
under suitable conditions on $f$, cf. Proposition
\thtag{560} and Proposition \thtag{5110}. 
In order to sum these two expressions
using Propositions \thtag{560} and \thtag{5110} we first 
have to sum the kernel $R_u$ introduced
in Proposition \thtag{5110}. The summation formula needed is
stated in the following lemma, which has been proved
by Mizan Rahman. 
The proof is presented in Appendix B. 

\proclaim{Lemma \thname{5150}}{\rm (Mizan Rahman)} We have 
for $|s|, |t|>1$, satisfying $st\not\in \pm q^{-\N}$,
$$
\multline
R_{q^{-2}}(\cos\psi;s,t|q^2) = 
{{(q^2,q^2,-q^2t^2s^2,-t^{-2}s^{-2},q{t\over s}e^{i\psi},
q{t\over
s}e^{-i\psi};q^2)_\infty}\over{(s^{-2},q^2s^{-2},-q^2,
-1,qste^{i\psi},qste^{-i\psi};q^2)_\infty}} 
- \\ {{(q^2t^2,q^3{t\over s}e^{i\psi},q^3{t\over s}e^{-i\psi},
{q\over{st}}e^{i\psi},{q\over{st}}e^{-i\psi}, q^3ste^{i\psi},
q^3ste^{-i\psi},-qste^{i\psi},-qste^{-i\psi};q^2)}
\over{(-1,-q^2,-q^2,-1, 
q^2e^{2i\psi},q^2e^{-2i\psi},qste^{i\psi},
qste^{-i\psi},q^2s^{-2};q^2)_\infty}}
\\ \times {{(-{q\over{st}}
e^{i\psi}, -{q\over{st}}e^{-i\psi};q^2)_\infty}\over{
(q^2s^{-2},q^4t^2;q^2)_\infty}}
\, {}_8W_7 (q^2t^2;q^2,q{t\over s}e^{i\psi},
q{t\over s}e^{-i\psi}, qste^{i\psi},qste^{-i\psi};q^2,q^2)
\endmultline
$$
and it remains valid for the discrete mass points
of the measure in Theorem \thtag{5100}, i.e. for
$e^{i\psi}= -stq^{1-2k}$, $k\in\Z$, $|-stq^{1-2k}|>1$. 
Explicitly,
$$
R_{q^{-2}}(\mu(-stq^{1-2k});s,t|q^2) = 
{{(q^2,q^2,-q^2t^2s^2,-t^{-2}s^{-2},
-q^{2-2k}t^2,-q^{2k}s^{-2};q^2)_\infty}\over
{(s^{-2},q^2s^{-2},-q^2,-1,-q^{2-2k}s^2t^2,-q^{2k};q^2)_\infty}}.
$$
\endproclaim

Rahman's lemma gives the explicit evaluation of the 
Poisson kernel corresponding to the little $q$-Jacobi function 
in one specific point $u=q^{-2}$. Note that
the explicit expression for $R_{q^{-2}}(\mu(-stq^{1-2k});s,t|q^2)$
follows from the fact that the second term in the
general expression vanishes since the factor in front of the
${}_8W_7$-series is zero and the ${}_8W_7$-series is 
non-singular for this value.

The condition $st\not\in \pm q^{-\N}$
ensures that the right hand
side in Lemma \thtag{5150} does not have simple 
poles for $\psi=0$, or $\psi=\pi$. In the subsequent application
of Lemma \thtag{5150} we multiply the result by the weight function
as in Proposition \thtag{5110}, which cancels the poles. 

\proclaim{Theorem \thname{5160}} Let $f$ be a
continuous, compactly supported function on the
spectrum of $\pt(\rst)$ and assume that
$f\bigl(\pt(\rst)\vert_{V^\th(t)}\bigr)PQ^2$ is 
of trace class on $V^\th(t)$ and that
$f\bigl(\pt(\rst)\vert_{W^\th(t)}\bigr)
(1-P)Q^2$ is of trace class on $W^\th(t)$.
Let $|s|\geq |t|>q^{-1}$ and
$s^2t^{\pm 2}\not\in q^{2\Z}$, then
$$
\align
{1\over{2\pi}} \int_0^{2\pi} &\text{\rm Tr}\vert_{\H} 
\bigl(f(\pt(\rst))Q^2\bigr)\, d\th = C 
\int_\R f(x)\,  d\nu (x;qs/t,qt/s;-qst|q^2), \\ 
C &= {{(q^2,q^2,-s^2,-q^2s^{-2},-t^2,-q^2t^{-2};q^2)_\infty}
\over{(t^2-1)(s^2-1)}}  
\endalign
$$
where the measure is defined in \eqtag{590}. 
\endproclaim

\demo{Remark} Note that the right hand side is symmetric in $s$ and
$t$ and invariant under $(s,t)\mapsto (-s,-t)$.
Since $\pi_{e^{i\th}}\circ\psi =
\pi_{e^{-i\th}}$ on $A_q(SU(1,1))$ and $\psi(\rst)=\rho_{-t,-s}$
with $\psi$ defined in Remark \thtag{3101} 
we see that the left hand side is 
also invariant under $(s,t)\mapsto (-t,-s)$. So the condition
$|s|\geq |t|$ is not essential. 
\enddemo

\demo{Proof} Propositions \thtag{5100} and 
\thtag{550} imply that the discrete spectrum of 
$f(\pt(\rst))\vert_{W^\th(t)}$ 
and  the discrete spectrum of 
$f(\pt(\rst))\vert_{V^\th(t)}$ 
do not overlap, but the continuous spectrum is the same 
in both cases. 
We consider the continuous and discrete
spectrum separately.

Let us consider the absolutely continuous part
on $[-1,1]$ first. Using Propositions \thtag{560} and \thtag{5110}
we have to consider
$$
\multline
{{(1-q^2/t^2)}\over{(t^2-1)}}
{{(1-e^{\pm 2i\psi})}\over{(1-{{qs}\over t}e^{\pm i\psi})
(1-{q\over{st}}e^{\pm i\psi})}}
\, {}_8W_7({{q^2}\over{t^2}};q^2,{{qs}\over t}e^{\pm i\psi},
{q\over{st}}e^{\pm i\psi};q^2,q^2) \\ +
{{(q^2s^{-2},-1,-q^2;q^2)_\infty^2}\over{(1-t^{-2})}}
{{ R_{q^{-2}}(\cos\psi;s,t|q^2)
(e^{\pm2i\psi};q^2)_\infty}\over{
({{qe^{\pm i\psi}}\over{st}},{{qte^{\pm i\psi}}\over{s}},
-qste^{\pm i\psi},-{{qe^{\pm i\psi}}\over{st}};q^2)_\infty}},
\endmultline
\tag\eq{5170}
$$
where we have also used \eqtag{590} and the $\pm$-signs means
that we have to take two terms, one with $+$ and one with $-$.
Now we  
can use Lemma \thtag{5150} to write $R_{q^{-2}}$ as a sum
of a very-well poised ${}_8W_7$-series and an explicit
term of infinite $q$-shifted factorials.  The two very-well
poised ${}_8W_7$-series can be summed using 
Bailey's summation formula 
see \cite{\GaspR, (2.11.7)}. In this case we write 
Bailey's formula as, cf.
\cite{\KoelV, p.~413},
$$
\multline
{{(1-ab) (1-e^{\pm 2i\psi})}\over{(1-ae^{\pm i\psi})
(1-be^{\pm i\psi})}} \, {}_8W_7(ab;q,ae^{\pm i\psi},
be^{\pm i\psi};q,q) \\ - {q\over{ab}}
{{(1-q^2/ab) (1-e^{\pm 2i\psi})}\over{
(1-qe^{\pm i\psi}/a) (1-qe^{\pm i\psi}/b)}}
\, {}_8W_7({{q^2}\over{ab}};q,{q\over a}e^{\pm i\psi},
{q\over b}e^{\pm i\psi};q,q) \\ = 
{{ (ab,q/ab, aq/b, bq/a,q,q;q)_\infty
(e^{\pm 2i\psi};q)_\infty}\over
{(ae^{\pm i\psi},be^{\pm i\psi}, qe^{\pm i\psi}/a,
qe^{\pm i\psi}/b;q)_\infty}}.
\endmultline
$$
Using Bailey's summation shows that \eqtag{5170}
equals
$$
\multline
{{(q^2,q^2,-q^2t^2s^2,-t^{-2}s^{-2},-1,-q^2,e^{\pm 2i\psi};
q^2)_\infty}\over{(1-s^{-2})(1-t^{-2})(qste^{\pm i\psi},
{q\over{st}}e^{\pm i\psi}, -qste^{\pm i\psi},
-{q\over{st}}e^{\pm i\psi};q^2)_\infty}} \\
- {{(q^2t^2,q^2t^{-2},q^2s^{-2},q^2s^2,q^2,q^2,
e^{\pm 2i\psi};q^2)_\infty}\over{ (q{t\over s}e^{\pm i\psi},
qste^{\pm i\psi}, q{s\over t}e^{\pm i\psi}, {q\over{st}}e^{\pm
i\psi} ;q^2)_\infty}} = %\\
{{(q^2,q^2,e^{\pm 2i\psi};q^2)_\infty}\over{
(1-t^{-2})(1-s^{-2})
(qste^{\pm i\psi}, {q\over{st}}e^{\pm i\psi};q^2)_\infty}}\\
\times
\Bigl(
{{(-q^2t^2s^2,-t^{-2}s^{-2},-1,-q^2;q^2)_\infty
}\over{(-qste^{\pm i\psi},
-{q\over{st}}e^{\pm i\psi};q^2)_\infty}} -
{{(t^{-2},q^2t^2,q^2s^2,s^{-2};q^2)_\infty}\over{
(q{s\over t}e^{\pm i\psi}, q{t\over s}e^{\pm i\psi};q^2)_\infty}}
\Bigr).
\endmultline
$$
Use the notation
$\Theta(a)=(a,q^2/a;q^2)_\infty$ for a theta-product in
base $q^2$ and
$S(a,b,c,d)=\Theta(a)\Theta(b)\Theta(c)\Theta(d)$ for
the product of four theta products in base $q^2$. The 
following identity for theta-products, 
$$
S(x\la,x/\la,\mu\nu,\mu/\nu)-S(x\nu,x/\nu,\la\mu,\mu/\la) =
{{\mu}\over{\la}} S(x\mu,x/\mu,\la\nu,\la/\nu),
\tag\eq{5175}
$$
see \cite{\GaspR, Ex.~2.16}, can be used to rewrite
the term in parentheses as
$$
\multline
{{S(-1,-t^{-2}s^{-2},qse^{i\psi}/t, qte^{i\psi}/s)
- S(t^{-2},s^{-2}, -qste^{i\psi}, -qe^{i\psi}/st)}\over
{S(-qste^{i\psi}, -qste^{-i\psi},qse^{i\psi}/t,
qse^{-i\psi}/t)}} = \\
-{q\over{st}}e^{i\psi} {{S(qe^{i\psi}/st,
e^{-i\psi}/qst,-s^2,-t^2)}\over
{S(-qste^{i\psi}, -qste^{-i\psi},qse^{i\psi}/t,
qse^{-i\psi}/t)}} =\\
{1\over{s^2t^2}}
{{(qe^{\pm i\psi}/st,qste^{\pm i\psi},-s^2,-q^2s^{-2},
-t^2,-q^2t^{-2};q^2)_\infty}\over{(qse^{\pm i\psi}/t,
qte^{\pm i\psi}/s,
-qe^{\pm i\psi}/st, -qste^{\pm i\psi};q^2)_\infty}}
\endmultline
$$
by taking $x=1/st$, $\mu=qe^{i\psi}$, $\la=-st$, $\nu=s/t$.
Plugging this in for the term in
parentheses we have evaluated \eqtag{5170} explicitly as
$$
{{(q^2,q^2,-s^2,-q^2s^{-2},-t^2,-q^2t^{-2};q^2)_\infty}
\over{(t^2-1)(s^2-1)}}  
{{(e^{\pm 2i\psi};q^2)_\infty}\over{
(-{q\over{st}}e^{\pm i\psi}, -qste^{\pm i\psi},
q{s\over t}e^{\pm i\psi},q{t\over s}e^{\pm i\psi};q^2)_\infty}}.
$$
This proves the statement concerning the absolutely continuous
part. 

It remains to check the discrete mass points. Since 
$|s|\geq |t|$, we only have an infinite set of
discrete mass points from Proposition \thtag{5110} and possibly
a finite set of discrete mass points from Proposition \thtag{560}. 
In case discrete 
mass points arise from Proposition \thtag{560} we have to verify 
$$
{-1\over{1-q^2}} 
{{w_k(qs/t;q/st,qt/s,qst|q^2)}
\over{h_0(qs/t,q/st,qt/s,qst|q^2)}}
= C\, w_k(qs/t;qt/s,-q/st,-qst|q^2)
$$
and this is a straightforward calculation using 
\eqtag{545} and
the value for $C$. 
For the infinite set of discrete mass points arising from 
Proposition \thtag{5110} we have by Lemma \thtag{5150} and
\eqtag{595} 
$$
\multline
{{(q^2s^{-2},-1,-q^2;q^2)_\infty^2}\over{1-t^{-2}}} 
R_{q^{-2}}(\mu(-stq^{1-2k});s,t|q^2) \\ \times 
\text{Res}_{z=-stq^{1-2k}} 
\Bigl( z^{-1}\, w(z;{q\over{st}}, q{t\over
s}, -{q\over{st}},-qst|q^2)\Bigr)  = 
{{q^{2(k-1)}(s^2t^2q^{2-4k}-1)}\over{(t^2-1)(s^2-1)}}
\endmultline
$$
using \eqtag{237}. From \eqtag{595} and \eqtag{237}
we also obtain 
$$
C\, \text{Res}_{z=-stq^{1-2k}} z^{-1}\, 
w(z;q{s\over t}, q{t\over s},
-{q\over{st}},-qst|q^2) = 
{{q^{2(k-1)}(s^2t^2q^{2-4k}-1)}\over{(t^2-1)(s^2-1)}},
$$
so that we have the desired result for the infinite set of
discrete mass points. 
\qed\enddemo

It follows directly from \eqtag{505} and \eqtag{565} 
that the unitary operator $T_t(e^{i\th})$ defined by
$w^\th_p(t)\mapsto e^{-2ip\th} w^0_p(t)$ and 
$v^\th_p(t)\mapsto e^{-2ip\th} v^0_p(t)$ satisfies
$T_t(e^{i\th})\pt(\rst)T_t(e^{i\th})^\ast=\pi_1(\rst)$. 
Note that $T_t(e^{i\th})$ is unitary by
Proposition \thtag{410}. So using \cite{\DixmvNA, Ch. II.2, \S 6}
we find that for a bounded continuous function $f$ 
$$
f(\pi(\rst)) = {1\over{2\pi}} \int_0^{2\pi} 
f(\pt(\rst))\, d\th  = T_t^\ast\bigl(\id \otimes f(\pi_1(\rst))
\bigr) T_t^\ast,
\tag\eq{5176}
$$
$T_t={1\over{2\pi}}\int_0^{2\pi} T_t(e^{i\th})\, d\th$, 
using $L^2(\T;\H)\cong L^2(\T)\otimes\H$ as tensor product of
Hilbert spaces, cf. \eqtag{415}. As before, $T_t$ commutes with
multiplication by a function from $L^2(\T)$, so that 
$f(\pi(\rst))$ is decomposable. So 
we can apply the Haar functional to it, see
Remark \thtag{2100}. 

\proclaim{Corollary \thname{5177}} Let $|s|\geq |t|> q^{-1}$
and $s^2t^{\pm 2}\not\in q^{2\Z}$. Let $f$ be a
continuous, compactly supported function, such
that $f(\pi_1(\rst))Q^2$ is of trace class on $\H$, 
then $f(\pi(\rst))$ is a 
decomposable operator from ${\Cal B}(L^2(\T;\H))$ and
$$
h\bigl( f(\pi(\rst)\bigr) = 
{{(q^2,q^2,-s^2,-q^2s^{-2},-t^2,-q^2t^{-2};q^2)_\infty}
\over{(t^2-1)(s^2-1)}} 
\int_\R f(x)\,  d\nu (x;qs/t,qt/s;-qst|q^2)
$$
with the measure defined in \eqtag{590}.
\endproclaim

\demo{Proof} Since $f(\pi_1(\rst))Q^2$ is of trace class on
$\H$, we have $f(\pi_1(\rst))PQ^2$ 
as trace class operator on $V^0(t)$ and 
$f(\pi_1(\rst))(1-P)Q^2$ as trace class operator on $W^0(t)$.
Then $f(\pt(\rst))PQ^2$ and $f(\pt(\rst))(1-P)Q^2$
are trace class operators on $V^\th(t)$ and $W^\th(t)$.
Now apply Theorem \thtag{5160} and Remark \thtag{2100}.
\qed\enddemo
 
%%%%%%%%%%%%%%%%%%%%%%%%%%%%%%%%%%%%%%%%%%%%%%%%%%%%%%%%%%%%%%%%%
\subhead \the\sectionnumber.4 The Haar functional\endsubhead
In this subsection we give the measure for the Haar functional
on a specific bi-$K$-type of the Cartan decomposition of
Theorem \thtag{3140}. This gives an explicit measure 
space of Askey-Wilson type. The Haar functional 
on bi-$K$-invariant elements is obtained in Corollary \thtag{5177}. 

In order to describe the Haar functional on the non-trivial
$K$-types of the Cartan decomposition of Theorem \thtag{3140}
we need to generalise the measure $d\nu(\cdot;a,b;d|q)$.
Define, cf. \eqtag{540},
$$
W_r(z;a,b,c;d|q) =
{{(z^2,z^{-2},qz/d,q/zd;q)_\infty}\over{
(rdz,rd/z,qz/rd,q/rdz,az,a/z,bz,b/z,cz,c/z;q)_\infty}}
$$
and observe that it differs from the Askey-Wilson weight
function by a quotient of theta functions;
$W_r(z;a,b,c;d|q) = \psi_r(z)w(z;a,b,c,d|q)$ with
$$
\psi_r(z) = {{(dz,q/dz,d/z,qz/d;q)_\infty}\over{
(rdz,q/rdz,rd/z,qz/rd;q)_\infty}}=\hat\psi_r(\mu(z)).
$$
The corresponding measure is defined in terms of the
Askey-Wilson measure of \eqtag{530} by
$$
\multline
\int_\R f(x)\, d\nu_r(x;a,b,c;d|q) = h_0(a,b,c,d|q) 
\int_\R f(x)\hat\psi_r(x)\, dm(x;a,b,c,d|q)\\
+ \sum_{k\in\Z, |rdq^{-k}|>1} f(\mu(q^{-k}rd))
\text{Res}_{z=rdq^{-k}} {{W_r(z;a,b,c;d|q)}\over z},
\endmultline
\tag\eq{5190} 
$$
cf. \eqtag{590}. Note that any possible discrete mass points
of the Askey-Wilson measure at $\mu(dq^k)$, $k\in\Zp$, are
annihilated by $\hat\psi_r(\mu(dq^k))=\psi_r(dq^k)=0$. So the
support of the measure defined in \eqtag{5190} is given by
$[-1,1]$, where the measure is absolutely continuous, plus
a finite discrete set of points of the form $\mu(eq^k)$,
$k\in\Zp$ such that $|eq^k|>1$ for $e=a$, $b$ or $c$ and an
infinite discrete set $\mu(rdq^k)$, $k\in\Z$, with $|rdq^k|>1$.
Note also that for $r=-1$, 
$c=q/d$ we obtain the definition 
\eqtag{590} of $d\nu(\cdot;a,b;d|q)$ as a special case.
Note that the measure $d\nu_r(\cdot;a,b,c;d|q)$ is symmetric
in $a$, $b$, and $c$. The measure
$d\nu_r(\cdot;a,b,c;d|q)$ is positive
if $r<0$, $0<b\leq a<d/q$, $0<c\leq a<d/q$, 
$bd\geq q$, $cd\geq q$, $ab<1$, $ac<1$, 
where we assume that $a$ is the largest of
the parameters $a$, $b$ and $c$.
For the general discussion of this measure 
we refer to \cite{\KoelSAW}.

For an element 
$\xi^{(i)}=\Ga^{(i)}_{l,m}(s,t)f(\rst)$ 
corresponding to the Cartan decomposition of Theorem \thtag{3140},
we define the corresponding quadratic form
$$
\langle \xi^{(i)}, \xi^{(i)}\rangle =
h\bigl( \bar f(\pi(\rst))
\pi(\Ga^{(i)}_{l,m}(s,t)^\ast \Ga^{(i)}_{l,m}(s,t)) 
f(\pi(\rst))\bigr).
$$
By \eqtag{5176} and \eqtag{3137} we regard the operator in parentheses
as a decomposable operator for suitable $f$ and we assume it satisfies
the conditions of Remark \thtag{2100}, cf. \S 4.4. 

\proclaim{Theorem \thname{5200}} Let 
$f$ be a continuous, compactly supported function on the
spectrum of $\pi_1(\rst)$ 
and $l\in\hf\Zp$, $m\in\{-l,-l+1,\ldots,l\}$,
$|s|\geq |t|>q^{-1}$, $s^2t^{\pm 2}\notin q^{2\Z}$. Assume that
$$
\bar f(\pi_1(\rst))
\pi_1(\Ga^{(i)}_{l,m}(s,t)^\ast \Ga^{(i)}_{l,m}(s,t)) 
f(\pi_1(\rst)) Q^2
$$ 
is of trace class on $\H$.
Then we have for positive constants $C_i$ independent of $f$, 
$$
\align
\langle \xi^{(1)}, \xi^{(1)}\rangle &=
C_1 \int_\R |f(x)|^2 \, d\nu_{-s^2t^2}(x;{s\over t} q^{1+2l-2m}, 
q{t\over s}, q^{1+2l+2m}st; {q\over{st}}|q^2), \\
\langle \xi^{(2)}, \xi^{(2)}\rangle &=
C_2 \int_\R |f(x)|^2 \, d\nu_{-1}(x;{s\over t} q^{1+2l+2m}, 
q{t\over s}, q^{1+2l-2m}/st; qst|q^2), \\
\langle \xi^{(3)}, \xi^{(3)}\rangle &=
C_3 \int_\R |f(x)|^2 \, d\nu_{-s^2t^2}(x;q{s\over t}, 
{t\over s}q^{1+2l-2m}, q^{1+2l+2m}st; {q\over{st}}|q^2), \\
\langle \xi^{(4)}, \xi^{(4)}\rangle &=
C_4 \int_\R |f(x)|^2 \, d\nu_{-1}(x;q{s\over t}, 
{t\over s}q^{1+2l+2m}, q^{1+2l-2m}/st; qst|q^2).
\endalign
$$
\endproclaim

\demo{Proof} This is a direct consequence of
Corollary \thtag{5177} and \eqtag{3137},
cf. the proof of Theorem \thtag{470}.
\qed\enddemo

Observe that the quadratic forms in Theorem \thtag{5200}
are not always positive definite, cf. Remark \thtag{472}.

The content of Remark \thtag{475} applies
here as well up to some minor changes.

%%%%%%%%%%%%%%%%%%%%%%%%%%%%%%%%%%%%%%%%%%%%%%%%%%%%%%%%%%%%%%%%%%%%
%N E W   S E C T I O N%
%%%%%%%%%%%%%%%%%%%%%%%%%%%%%%%%%%%%%%%%%%%%%%%%%%%%%%%%%%%%%%%%%%%%
\head \newsec Spherical Fourier transforms\endhead

In this section we give a formal interpretation of the Askey-Wilson
function transform as studied in \cite{\KoelSAW} as a Fourier
transform on the quantum $SU(1,1)$ group. Parts of these results
only hold at a formal level, so this section mainly serves as the
motivation for the study of the Askey-Wilson function transform. In
this section we derive which symmetric operators on the function 
spaces of Theorem \thtag{5200} have to be studied 
and what are the natural
eigenfunctions of this operator to be considered. For the spherical
case we calculate the spherical functions
and the related action of the Casimir operator. For the Fourier
transforms related to the other parts of the Cartan decomposition we
only sketch parts of the formal arguments. 

%%%%%%%%%%%%%%%%%%%%%%%%%%%%%%%%%%%%%%%%%%%%%%%%%%%%%%%%%%%%%%%%%%%%
\subhead \the\sectionnumber.1\ Unitary representations
of $\U$\endsubhead 
The irreducible unitary representations, i.e.
$\ast$-representations, of $\U$, are known,
see Burban and Klimyk \cite{\BurbK}, Masuda et al.
\cite{\MasuMNNSU}, Vaksman and Korogodski\u\i\ \cite{\VaksK}. 
We are only interested in the admissible representations, i.e. 
we require that the eigenvalues of $A$ are contained in 
$q^{\hf \Z}$, that the corresponding eigenspaces are
finite-dimensional, and that the direct sum of the eigenspaces
is equal to the representation space.  
We now recall the classification, see 
Masuda et al. \cite{\MasuMNNSU} and 
Burban and Klimyk \cite{\BurbK} for a more general
situation. 
The irreducible admissible 
unitary representations act in $\ell^2(\Zp)$ or in $\H$, 
and we use $\{ e_n\}$ with $n\in\Zp$ or $\Z$ for the standard
orthonormal basis of $\ell^2(\Zp)$ or $\H$. There are, apart
from the trivial representation, 
five types of representations; positive discrete series,
negative discrete series, principal unitary series, 
complementary series and strange series. 
The representations are in terms of unbounded operators on
$\H$ or on $\ell^2(\Zp)$ with common domain the finite linear
combinations of the standard basis vectors $e_k$, cf. \S 2.1.
We also give the action of
the Casimir operator $\Om$, see \eqtag{315}, in each of the
irreducible admissible representations. 
The Casimir operator is central, so it acts by a scalar 
$[\la+\hf]^2$ for some $\la\in\C$, where 
$[a]=(q^a-q^{-a})/(q-q^{-1})$ is the $q$-number. 
The eigenvalues of $A$ are contained in $q^{\ep+\Z}$ for
$\ep=0$ or $\ep=\hf$. 

In the following list of irreducible admissible unitary
representations of $\U$ we give the action of the generators
of $\U$ on the orthonormal basis $\{e_k\}$, and the $\la\in\C$
corresponding to the action of $\Om$ and $\ep\in\{0,\hf\}$
corrresponding to the set $q^{\ep+\Z}$ in which $A$ takes its
eigenvalues. We remark that the scalar $[\la+\hf]^2$ and the 
the eigenvalues of $A$ determine the irreducible admissible
representation of $\U$ up to equivalence.  

\subsubhead Positive discrete series\endsubsubhead
The representation space is $\ell^2(\Zp)$.
Let $k\in\hf\N$, and $\la=-k$, 
and define the action of the generators by
$$
\align 
A\cdot e_n&=q^{k+n}e_n,\qquad D\cdot e_n=q^{-k-n}e_n,\\
(q^{-1}-q)B\cdot e_n&=q^{-\hf-k-n}
\sqrt{(1-q^{2n+2})(1-q^{4k+2n})}\,e_{n+1}\\
(q^{-1}-q)C\cdot e_n&=-q^{\hf-k-n}
\sqrt{(1-q^{2n})(1-q^{4k+2n-2})}\,e_{n-1}
\endalign
$$
with the convention $e_{-1}=0$.
We denote this representation by $T_k^+$. Now
$\ep=\hf$ if $k\in\hf+\N$ and $\ep=0$ if $k\in\N$.

\subsubhead Negative discrete series\endsubsubhead
The representation space is $\ell^2(\Zp)$.
The negative discrete
series representation is $T_k^-=T_k^+\circ \phi$, 
where $\phi:\U\to\U$ is the $\ast$-algebra involution defined by
$\phi(A)=D$, $\phi(B)=C$. The parameters $\la$ and $\ep$ are
the same as for the positive discrete series.

\subsubhead Principal series\endsubsubhead 
The representation space is $\H$.
Let $\la=-\hf+ib$ with $0\leq b\leq-\frac{\pi}{2\ln q}$ and
$\ep\in \{0,\hf\}$ and
assume  $(\la,\ep)\not=(-\hf,\hf)$. The action of the 
generators is defined by 
$$
\align
A\cdot e_n&=q^{n+\ep}e_n,\qquad D\cdot e_n=q^{-n-\ep}e_n,\\
(q^{-1}-q)B\cdot e_n&=q^{-\hf-n-\ep-ib}(1-q^{1+2n+2\ep+2ib})
\,e_{n+1}, \\
(q^{-1}-q)C\cdot e_n&=-q^{\hf-n-\ep+ib}
(1-q^{-1+2n+2\ep-2ib})\,e_{n-1}.
\endalign
$$
We denote the representation by $T_{\la,\ep}^{P}$. In case
$(\la,\ep)=(-\hf,\hf)$ this still defines an admissible 
unitary representation. It splits as the direct sum 
$T_{-\hf,\hf}^P= T_{\hf}^+\oplus T_{\hf}^-$ of a
positive and negative discrete series representation
by restricting to the invariant 
subspaces $\text{span}\{ e_n\mid n\geq 0\}$ and 
to $\text{span}\{ e_n\mid n <0\}$. 

\subsubhead Complementary series\endsubsubhead  
The representation space is $\H$.
Let $\ep=0$ and $-\hf<\la<0$. The action of the 
generators is defined by 
$$
\align
A\cdot e_n&=q^ne_n,\qquad D\cdot e_n=q^{-n}e_n,\\
(q^{-1}-q)B\cdot e_n&= q^{-n-\hf} \sqrt{(1-q^{2\la+2n+2})
(1-q^{2n-2\la})}\,e_{n+1},\\
(q^{-1}-q)C\cdot e_n&=-q^{-n+\hf} \sqrt{(1-q^{2\la+2n})
(1-q^{2n-2\la-2})}\,e_{n-1}.
\endalign
$$
We denote this representation by $T_{\la,0}^{C}$. 

\subsubhead Strange series\endsubsubhead
The representation space is $\H$. Let $\ep\in \{0,\hf\}$, 
and put $\la=-\hf - \frac{i\pi}{2\ln q} +a$, $a>0$.
The action of the generators is defined by 
$$
\align
A\cdot e_n&=q^{n+\ep}e_n,\qquad D\cdot e_n=q^{-n-\ep}e_n,\\
(q^{-1}-q)B\cdot e_n&= q^{-n-\ep-\hf} \sqrt{(1+q^{2n+2\ep+1+2a})
(1+q^{2n+2\ep-2a+1})}\,e_{n+1},\\
(q^{-1}-q)C\cdot e_n&=-q^{-n-\ep+\hf} \sqrt{(1+q^{2n+2\ep-1+2a})
(1+q^{2n+2\ep-2a-1})}\,e_{n-1}.
\endalign
$$
We denote this representation by  $T_{\la,\ep}^{S}$.

\topinsert
\hskip1truecm
{\epsfxsize 5.25truein \epsfysize 3truein \epsfbox{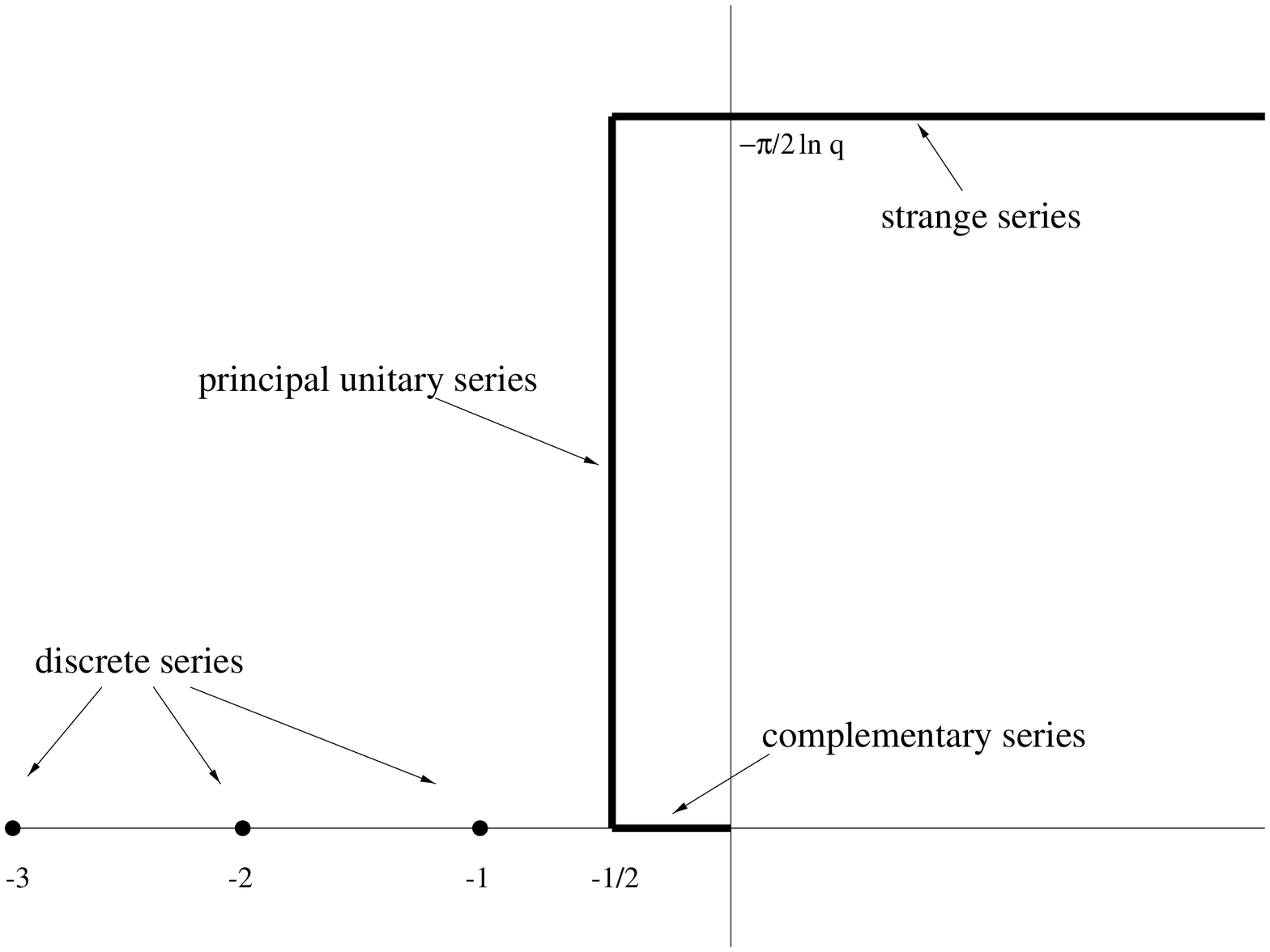}}
\botcaption{Figure 6.1}
Representation label $\la$ in the complex plane 
for $\ep=0$.
\endcaption
\endinsert 

\demo{Remark \thname{610}}
The matrix elements of the irreducible 
admissible unitary representations 
in terms of the standard basis $\{ e_k\}$ 
satisfy \eqtag{3120} for $(s,t)=(\infty,\infty)$, and
can be written in the form
$\Ga^{(i)}_{l,m}(\infty,\infty)g(\rho_{\infty,\infty})$, 
where $g$ is a power series and $\rho_{\infty,\infty}=\ga^\ast\ga$,
cf. the Cartan decomposition of Theorem \thtag{3140}. 
The corresponding power series 
have been calculated explicitly in \cite{\MasuMNNSU},  
see also \cite{\VaksK}, 
using the explicit duality between $\U$ and $A_q(SU(1,1))$.
The power series is a ${}_2\vp_1$-series, and can be interpreted
as a little $q$-Jacobi function. In the next subsections
we compute explicitly matrix coefficients which behave as a
character under the left $AY_t$, 
respectively right $Y_sA$-action, and we identify them 
with big $q$-Jacobi functions and with Askey-Wilson functions. 
\enddemo

%%%%%%%%%%%%%%%%%%%%%%%%%%%%%%%%%%%%%%%%%%%%%%%%%%%%%%%%%%%%%%%%%%%%
\subhead \the\sectionnumber.2\ $K$-fixed vectors\endsubhead
In this subsection we look for (generalised) eigenvectors
of $Y_sA \in\U$ for the eigenvalue zero. For the 
discrete series representations this involves the
Al-Salam and Chihara polynomials. For the other series
this involves transforms with a ${}_2\vp_1$-series as
kernel, which can be considered as Al-Salam and 
Chihara functions. 
For each of the series of representations of $\U$
we use the notation of \S 6.1.

\subsubhead Irreducible representations in $\ell^2(\Zp)$
\endsubsubhead
For the positive
discrete series representation the spectrum
of $Y_sA$ has been calculated 
in \cite{\KoelVdJ, \S 4} using the Al-Salam and Chihara
polynomials, see \S 5.1. From \cite{\KoelVdJ, Prop.~4.1}
we conclude that zero is a (generalised) 
eigenvalue of $T^+_k(Y_sA)$ 
if and only if $\mu(s)$ is in the support of the
orthogonality measure $dm(\cdot;q^{2k}s,q^{2k}/s|q^2)$, 
where we use the notation of \S 5.1. Since $|s|>1$ we have 
$|\mu(s)|>1$, so this can only happen if there 
exists a discrete mass point. As $k>0$ implies 
$|q^{2k}/s|<1$, we have to have $q^{2k+2n}s=s$ for
some $n\in\Zp$, which is not possible. 
We conclude that in the positive discrete series we do not
have eigenvectors of $Y_sA$ for the eigenvalue zero. 

It follows immediately from \S 6.1 and \eqtag{360} that
for the spectrum of $T_k^-(Y_sA)$ 
we have to study the recurrence relation
$$
\gathered
((q-q^{-1})T^-_k(Y_sA) -s-s^{-1}) \cdot e_n = 
a_n\, e_{n-1} + b_n\, e_n +a_{n+1}\, e_{n+1},\qquad n\in\Zp, \\
a_n = q^{1-2k-2n}\sqrt{(1-q^{2n})(1-q^{4k+2n-2})},
\quad b_n=-q^{-2k-2n}(s+s^{-1}).
\endgathered
\tag\eq{620}
$$
In order to determine the spectrum 
we have to study the corresponding orthonormal 
polynomials $x\, p_n(x) = 
a_n\, p_{n-1}(x) + b_n\, p_n(x) +a_{n+1}\, p_{n+1}(x)$. 
We let
$$
p_n(x) = (-1)^n q^{n(1+2k)} {{(q^2;q^2)_n^\hf}\over{
(q^{4k};q^2)_n^\hf}} \, P_n(x),
$$
so that $P_n(x)$ satisfies the recurrence relation
$$
(1-q^{2n+2})P_{n+1}(x) = (-q^{-2k}(s+s^{-1}) - xq^{2n})\,
P_n(x) - (q^{-4k}-q^{2n-2})\, P_{n-1}(x).
$$
This is precisely the form of the Al-Salam and Chihara
polynomials in base $q^{-2}>1$ as studied by
Askey and Ismail \cite{\AskeI, \S 3.12, \S 3.13}. 
It follows from \cite{\AskeI, Thm.~3.2} that 
for $|s|\geq q^{-1}$ the associated moment problem
is determinate, and in that case the support of
the orthogonality measure is $\{ 2\mu(-sq^{-2p-2k})\mid
p\in\Zp\}$, \cite{\AskeI, (3.80)-(3.82)}. 
Part of this statement can also be found
in \cite{\AlSaC, p.26}. So we see that $T^-_k(Y_sA)$
has an eigenvalue zero if and only if 
$2\mu(-s)$ is in the support of the orthogonality
measure. Since $k>0$ and $s\in\R$, $|s|\geq q^{-1}$, 
we see that this is impossible. 
We conclude that in the negative discrete series we do not
have eigenvectors of $Y_sA$ for the eigenvalue zero. 

We formalise this into the following lemma.

\proclaim{Lemma \thname{625}} Let $s\in\R$.
For $|s|\geq 1$, the operator $T^+_k(Y_sA)$ is self-adjoint,
and zero is not contained in its spectrum. 
For $|s|\geq q^{-1}$ the operator $T^-_k(Y_sA)$ is essentially
self-adjoint, and zero is not contained in its spectrum.  
\endproclaim 

\subsubhead Irreducible representations in $\H$
\endsubsubhead
For the spectrum of $Y_sA$ 
in the case of the principal unitary, complementary
and strange series representations  
we have to study the recurrence relation
$$
((q^{-1}-q)T^\bullet_{\la,\ep}(Y_sA)  +s+s^{-1}) \cdot e_n = 
a_n^\bullet\, e_{n-1} + b_n\, e_n + 
\overline{a_{n+1}^\bullet}\, e_{n+1}, 
\tag\eq{630}
$$
with $b_n=q^{2n+2\ep}(s+s^{-1})$ and 
$\bullet\in\{ P,C,S\}$ for the principal unitary,
complementary or strange series. The various values for
$a_n^\bullet$, $\bullet\in\{ P,C,S\}$, follow immediately
from \eqtag{360} and \S 6.1;
$$
\aligned
a_n^P &= q^{ib}(1-q^{2n-1+2\ep-2ib}), 
\qquad \ep\in \{0,\hf\}, 
\quad 0\leq b\leq \frac{-\pi}{2\ln q}, \quad (b,\ep)\not=(0,\hf),\\
a_n^C& = \sqrt{(1-q^{2n+2\la})(1-q^{2n-2-2\la})},
\qquad -\hf<\la<0, \\
a_n^S& = \sqrt{(1+q^{2n+2\ep-1+2a})(1+q^{2n+2\ep-1-2a})},
\qquad \ep\in \{0,\hf\}, \quad a>0.
\endaligned
\tag\eq{640}
$$
Define a new orthonormal basis $\{f^\bullet_k\}_{k\in\Z}$ of
$\H$ by $f^\bullet_n=e^{i\phi^\bullet_n}e_{-n}$, 
with $\phi^\bullet_n$ a sequence of real numbers satisfying 
$\phi^\bullet_n=\phi^\bullet_{n+1}-\arg (a^\bullet_{-n})$, then
$$
((q^{-1}-q)T^\bullet_{\la,\ep}(Y_sA)  +s+s^{-1}) \cdot 
f^\bullet_n = 
a_n\, f^\bullet_{n+1} + b_n\, f^\bullet_n + a_{n-1}\,
f^\bullet_{n-1},
$$
with $a_n$, $b_n$ as in Lemma~A.2 %harde referentie 
in base $q^2$ and $c$, $d$, $z$ replaced by $q^2s^{-2}$, 
$q^{2+2\la}s^{-1}$, $q^{-2\ep-2\la}$,
where $\la=-\hf+ib$, $0\leq b\leq -\pi/2\ln q$, $\ep\in\{0,\hf\}$, 
$(b,\ep)\not=(0,\hf)$ for $\bullet=P$,
$\ep=0$, $-\hf<\la<0$ for $\bullet=C$, and
$\la=-\hf-i\pi/2\ln q + a$, $a>0$, $\ep\in\{0,\hf\}$
for $\bullet=S$. This recurrence relation is related to the
second order $q$-difference equation for the ${}_2\vp_1$-series, 
see Appendix A. In case of the principal unitary series, the
parameters satisfy the conditions of case (1) of 
Lemma~A.2. %harde referentie 
For the complementary series, the
parameters satisfy the conditions of case (2) of 
Lemma~A.2. %harde referentie 
and for the strange series, the
parameters satisfy the conditions of case (3) of 
Lemma~A.2. %harde referentie 
Now Theorem~A.5 implies the
following result, since all conditions are met in the respective
cases for $|s|\geq q^{-1}$.

\proclaim{Proposition \thname{650}} 
Assume $|s|\geq q^{-1}$ so that $T^\bullet_{\la,\ep}(Y_sA)$
is essentially self-adjoint. With the
notation of \S 6.1 we have that zero is in the discrete
spectrum of $T^\bullet_{\la,0}(Y_sA)$, 
$\bullet\in\{ P,C,S\}$, and zero is not in the 
spectrum of $T^\bullet_{\la,\hf}(Y_sA)$. The eigenvector
$v_s^\bullet$ of $T^\bullet_{\la,0}(Y_sA)$ 
for the eigenvalue zero is given by $v^\bullet_s =$
$$
 \sum_{n=-\infty}^\infty
e^{i(\psi_{-n}+\phi^\bullet_{-n})} |q^{2-2\la}s|^n 
\sqrt{ {{(q^{-2\la+2n};q^2)_\infty}\over{
(q^{2\bar \la+2+2n};q^2)_\infty}}}
\, {}_2\vp_1\left( {{q^{2+2\la}s^{-2}, q^{2+2\la}}\atop
{q^2s^{-2}}};q^2, q^{-2n-2\la}\right)\, e_n
\tag\eq{651}
$$
where $\psi_{k+1}=\psi_k+\arg
(q^{2+2\la}s^{-1}(1-q^{-2k+2\bar\la}))$. 
Here we use the analytic continuation of the ${}_2\vp_1$-series
as described in \S A.2
\endproclaim

\demo{Remark} Since we have $Y_\infty=A-D$, we obtain directly from
\S 6.1 that $T^\bullet_{\lambda,\ep}(A^2-1)$ has a vector in its
kernel if and only if $\ep=0$ and $\bullet\in \{P,C,S\}$, and in
that case $e_0$ spans the kernel. We can formally obtain this from
Proposition \thtag{650} by taking termwise limits $s\to\infty$. 
\enddemo

%%%%%%%%%%%%%%%%%%%%%%%%%%%%%%%%%%%%%%%%%%%%%%%%%%%%%%%%%%%%%%%%%%%%
\subhead \the\sectionnumber.3\ 
Zonal spherical functions\endsubhead 
In this subsection we will now give a 
formal derivation of the zonal spherical functions
that may occur in the spherical Fourier transform on the
quantum $SU(1,1)$ group. 

With the eigenvectors of Proposition \thtag{650} at hand
we can consider the linear functional 
$$
f^\bullet_\la \colon \U\to\C, 
\qquad X\mapsto  \langle T^\bullet_{\la,0}(XA) v_t^\bullet,
v_s^\bullet \rangle_{\H},
\tag\eq{660}
$$
where we take $s,t\in\R$ with $|s|,|t|\geq q^{-1}$. 
Then we formally have for $\bullet\in\{ P,C,S\}$ 
$$
\aligned
(Y_t.f^\bullet_\la)(X) &= \langle T^\bullet_{\la,0}
(XY_tA)v_t^\bullet, v_s^\bullet \rangle_{\H} = 0, \\
(f^\bullet_\la.Y_s)(X) &= \langle T^\bullet_{\la,0}
(Y_sXA)v_t^\bullet, v_s^\bullet \rangle_{\H} 
= \langle T^\bullet_{\la,0}
(DXA)v_t^\bullet, T^\bullet_{\la,0}(Y_sA)v_s^\bullet \rangle_{\H} 
= 0
\endaligned
\tag\eq{670}
$$
by Proposition \thtag{650} and the fact that $T^\bullet_{\la,0}$
is unitary and $(Y_sA)^\ast=Y_sA$. 
Note that \eqtag{660} and \eqtag{670} 
can be made rigorous for the
limit case $s\to\infty$, so that $Y_sA$ has to be replaced
by $A^2-1$, which has only an eigenvector for the eigenvalue
zero in the principal unitary series, complementary series and
the strange series for $\ep=0$. In these cases $e_0$ is
the eigenvector, and the analogue of 
$f^\bullet_\la \colon \U\to\C$ for this case is
given by $X\mapsto  \langle v_t^\bullet,
T^\bullet_{\la,0}(AX^\ast)e_0 \rangle_{\H}$, which is well-defined
for every $X\in\U$ since $T^\bullet_{\la,0}(AX^\ast)e_0$ has only
finitely many terms. 

The matrix elements $T^\bullet_{\la,0;n,m}\colon 
X\mapsto \langle T^\bullet_{\la,0}(X)e_m,e_n\rangle$  
have been calculated by 
Masuda et al. \cite{\MasuMNNSU}, see Remark \thtag{610}.
We can formally write 
$$
f^\bullet_\la= \sum_{n,m=-\infty}^\infty 
\langle v^\bullet_t,e_m\rangle \langle e_n,v^\bullet_s\rangle
q^m\, T^\bullet_{\la,0;n,m}.
\tag\eq{680}
$$
Note again that for the case $s\to\infty$ this can be
made rigorous, since the double sum reduces to a single
sum and pairing with an arbitrary element $X\in\U$ gives
a finite sum in this case. 

Because of \eqtag{670} we consider $f^\bullet_\la$ 
as the zonal spherical function. Because of the Cartan
decomposition of Theorem \thtag{3140} we formally have that 
the expression in \eqtag{680} is a function in 
$\rst$, 
$$
f^\bullet_\la= \sum_{n,m=-\infty}^\infty 
\langle v^\bullet_t,e_m\rangle \langle e_n,v^\bullet_s\rangle
q^m\, T^\bullet_{\la,0;n,m}= 
\phi_\la(\rst).
$$
In order to determine $\phi_\la$ we evaluate
at $A^\nu$, $\nu\in\Z$, with $A^{-1}=D$. Since 
$A$ is group-like, i.e. $\De(A)=A\otimes A$, we have that
pairing with $A$ is a homomorphism. Since
$\rst(A^\nu) = \mu(q^{\nu})$ and
$T^\bullet_{\la,0;n,m}(A^\nu) = 
\langle T^\bullet_{\la,0}(A^\nu)e_m,e_n\rangle = 
q^{n\nu}\de_{n,m}$ we see that we can determine $\phi_\la$ from
$$
\phi_\la(\mu(q^\nu)) = 
\sum_{n=-\infty}^\infty 
\langle v^\bullet_t,e_n\rangle \langle e_n,v^\bullet_s\rangle
q^{n(\nu+1)}.
$$

Now we can use the following summation formula. 
This lemma has been proved by
Mizan Rahman, and the proof is given in Appendix B. 

\proclaim{Lemma \thname{690}}{\rm (Mizan Rahman)} 
Let $|s|, |t|\geq 1$, assume $st>0$. For 
$\la$ corresponding to the principal unitary series,
complementary series and strange series, i.e.
$\la=-\hf+ib$, $0\leq b\leq -\pi/2\ln q$, or
$-\hf <\la<0$, or $\la=-\hf +a -i\pi/2\ln q$, $a>0$, 
and for $z$ in the annulus 
$|q/st|<|z|<|st/q|$ we have  
$$
\multline
\sum_{n=-\infty}^\infty 
\langle v^\bullet_t,e_n\rangle \langle e_n,v^\bullet_s\rangle
q^nz^n = {{(q^2,q^2/s^2t^2,
q^{1-2\la}/zst,q^{1-2\la}z/st;q^2)_\infty}\over
{(q^{-2\la},q^{2-2\la}/s^2t^2,qz/st,q/zst;q^2)_\infty}} \\ \times\, 
{}_8W_7(q^{-2\la}/s^2t^2;q^{-2\la}/s^2, q^{-2\la}/t^2,
q^{-2\la}, qz/st,q/zst;q^2,q^{2+2\la}).
\endmultline
$$
\endproclaim

{}From Lemma \thtag{690} we formally conclude that the 
spherical elements of \eqtag{680} can be
expressed in terms of a very-well-poised
${}_8W_7$-series;
$$
\multline
f^\bullet_\la = \phi_\la(\rst) = 
{{(q^2, q^2/s^2t^2,q^{1-2\la}/zst,q^{1-2\la}z/st;q^2)_\infty}\over
{(q^{-2\la},q^{2-2\la}/s^2t^2,qz/st,q/zst;q^2)_\infty}} \\\times \, 
{}_8W_7(q^{-2\la}/s^2t^2;q^{-2\la}/s^2, q^{-2\la}/t^2,
q^{-2\la}, qz/st,q/zst;q^2,q^{2+2\la})\Big\vert_{
\mu(z)=\rst}.
\endmultline
\tag\eq{6100}
$$
Note that the right hand side is symmetric in $z$ and $z^{-1}$,
so that we can make this specialisation. 

\demo{Remark \thname{6105}} Using the limit transition $-2\rst/qs\to \rit$
as $s\to\infty$, see \eqtag{3105}, 
we formally obtain the limit case $s\to\infty$
of the spherical function in \eqtag{6100} as
$$
{{(q^2,-q^{2-2\la}\rit/t;q^2)_\infty}\over{(q^{-2\la},
-q^2\rit/t;q^2)_\infty}} \, 
{}_3\vp_2 \left( {{q^{-2\la}t^{-2}, q^{-2\la},-q^2\rit t^{-1}}
\atop{q^2t^{-2}, -q^{2-2\la}\rit t^{-1}}};q^2, q^{2+2\la}\right),
$$
which is, up to a scalar, a special case of the function 
considered in \cite{\KoelSbig, (3.8)} and is the big $q$-Legendre
function. Next letting $t\to\infty$ using
$t^{-1}\rit\to\rho_{\infty,\infty}=\ga^\ast\ga$, cf. \eqtag{3105}, 
we see that the spherical function in the case $s,t\to\infty$ is
$$
\multline
{{(q^2,-q^{2-2\la}\rho_{\infty,\infty};q^2)_\infty}\over{(q^{-2\la},
-q^2\rho_{\infty,\infty};q^2)_\infty}} \, 
{}_2\vp_1 \left( {{q^{-2\la},-q^2\rho_{\infty,\infty}}
\atop{-q^{2-2\la}\rho_{\infty,\infty}}};q^2, q^{2+2\la}\right) \\
= {{(q^2,q^2;q^2)_\infty}\over{(q^{-2\la}, q^{2+2\la};q^2)_\infty}} \,
{}_2\vp_1 \left( {{q^{2+2\la},q^{-2\la}}\atop{q^2}};q^2, -q^2
\rho_{\infty,\infty}\right)
\endmultline
$$
by \cite{\GaspR, (1.4.5)}. This gives back the spherical
function, the little $q$-Legendre function, 
as studied by Kakehi, Masuda and Ueno \cite{\KakeMU}
and Vaksman and Korogodski\u\i\ \cite{\VaksK}. So the 
function $\phi_\la$ of \eqtag{6100} is a 2-parameter extension 
of the little $q$-Legendre function. 
\enddemo 

%%%%%%%%%%%%%%%%%%%%%%%%%%%%%%%%%%%%%%%%%%%%%%%%%%%%%%%%%%%%%%%%%%%%
\subhead \the\sectionnumber.4\ 
The action of the Casimir element\endsubhead 
Since the Casimir element $\Om$ acts in any of the irreducible
unitary representations of \S 6.1 by the
constant $[\la+\hf]^2$, we see from \eqtag{680} that
we formally have that the spherical function is an
eigenfunction of the action of the Casimir operator; 
$\Om.f^\bullet_\la=[\la+\hf]^2f^\bullet_\la$.

On the other hand, observe that the $(s,t)$-spherical
elements as defined in Proposition \thtag{3110} are invariant
under the action of the Casimir operator, since 
$\Om$ is in the centre of $\U$. So we can restrict its
action to the subalgebra of $(s,t)$-spherical elements,
or the subalgebra generated by $\rst$. For this
we have to calculate the radial part of $\Om$, and this
is stated in the following lemma. The proof is
the same as Koornwinder's proof of \cite{\KoorZSE, Lemma 5.1},
so we skip the proof. 

\proclaim{Lemma \thname{6110}} 
Put 
$$
\psi(z)= {{(1-qstz)(1-qsz/t)(1-qzt/s)(1-qz/st)}
\over{(1-z^2)(1-q^2z^2)}},
$$
then
$$
q(q^{-1}-q)^2A^{\nu}\Omega \equiv 
\psi(q^{\nu})\bigl(A^{\nu+2}-A^{\nu}\bigr)
+\psi(q^{-\nu})\bigl(A^{\nu-2}-A^{\nu}\bigr)
+(1-q)^2A^{\nu}
$$
modulo $\U Y_t+Y_s\U$.
\endproclaim

As for the polynomial case discussed by Koornwinder \cite{\KoorZSE},
we derive from this equation that the action of the
Casimir operator on the subalgebra of $(s,t)$-spherical
elements is given by the Askey-Wilson $q$-difference
operator \cite{\AskeW};
$$
\multline
q(q^{-1}-q)^2\Om.(f(\rst))= 
\Bigl( \psi(q^\nu)\bigl(f(\mu(q^{\nu+2}))-f(\mu(q^{\nu}))\bigr)
+ \\\psi(q^{-\nu})\bigl(f(\mu(q^{\nu-2}))-f(\mu(q^{\nu}))\bigr) 
+(1-q)^2f(\mu(q^{\nu}))\Bigr)\Big\vert_{\mu(q^\nu)=\rst}.
\endmultline
\tag\eq{6120}
$$
Combining \eqtag{6120} with the scalar action of $\Om$ in
the irreducible representations we formally find that the 
spherical function $\phi_\la(\mu(z))$ is an eigenfunction
of 
$$
\gathered
L\phi_\la(\mu(z)) = \bigl(-1-q^2 +q(q^{2\la+1}+q^{-2\la-1})
\bigr) \phi_\la(\mu(z)), \\
L = \psi(z)(T_{q^2}-1) +\psi(z^{-1})(T_{q^{-2}}-1), 
\qquad (T_q f)(z)= f(qz).
\endgathered
\tag\eq{6130}
$$
This is only a formal derivation, due to the fact that the
series \eqtag{680} is only a formal expression. 
Note that the eigenvalues in \eqtag{6130}
are real for $\la$ corresponding to the principal
unitary series, complementary series and strange 
series. The
function $\phi_\la(\mu(z))$ given in 
\eqtag{6100} is indeed an eigenfunction of the Askey-Wilson 
$q$-difference equation as in \eqtag{6130}, see
Ismail and Rahman \cite{\IsmaR}, Suslov \cite{\SuslJPA}, 
\cite{\SuslPP}. So we call $\phi_\la$ of \eqtag{6100} an
Askey-Wilson function. 

\demo{Remark \thname{6135}} For the limit case $s\to\infty$ we 
obtain the same eigenvalue equation as in \eqtag{6130} but now 
with the operator
$$
\gathered
L = A(z)(T_{q^2}-1) +B(z)(T_{q^{-2}}-1),\\
A(z) = q^2(1+{1\over{q^2tz}})(1+{t\over{q^2z}}), \quad
B(z) = (1+{1\over{tz}})(1+{t\over z}).
\endgathered
\tag\eq{6137}
$$
Then it is known \cite{\GuptIM} 
that the spherical function given in Remark 
\thtag{6105} is indeed a solution to the eigenvalue equation.
See \cite{\KoelSbig} for more information. For the limit case 
$s,t\to\infty$ we find the same eigenvalue equation \eqtag{6137},
but with now 
$A(z)=q^2(1+q^{-2}z^{-1})$ and $B(z)=1+z^{-1}$.
The little $q$-Legendre function as in Remark \thtag{6105} is a solution
of the eigenvalue equation as follows from \thetag{A.8}, see
also \cite{\Kake}, \cite{\KakeMU}, \cite{\VaksK}. 
\enddemo

\proclaim{Proposition \thname{6140}} The action of the Casimir
operator on the space of 
$(s,t)$-spherical elements is symmetric, i.e. 
$$
\int_\R (Lf)(x) \bar g(x)\,  d\nu (x;qs/t,qt/s;-qst|q^2)
= \int_\R f(x) \overline{(Lg)(x)}\,  d\nu (x;qs/t,qt/s;-qst|q^2)
$$
for continuous, compactly supported functions
$f$ and $g$ such that the functions $F(z)=f(\mu(z))$ and
$G(z)=g(\mu(z))$ have an analytic continuation to a neighbourhood 
of  $\{ z\in\C\mid q^2\leq |z|\leq q^{-2}\}$.
\endproclaim

So we interpret this as $h\bigl( g(\rst)^\ast \Om.f(\rst)\bigr)
= h\bigl( (\Om.g(\rst))^\ast f(\rst)\bigr)$ using
Corollary \thtag{5177}. 

\demo{Proof} This is a calculation using Cauchy's theorem
and shifting sums, see \cite{\KoelSAW} for details. 
\qed\enddemo

Proposition \thtag{6140} remains valid for the 
limit case $s\to\infty$
with the same proof, see \cite{\KoelSbig}. Taking furthermore
$t\to\infty$ leads to the situation considered by Kakehi, Masuda 
and Ueno \cite{\KakeMU}, see also \cite{\Kake}, \cite{\VaksK}. 

%%%%%%%%%%%%%%%%%%%%%%%%%%%%%%%%%%%%%%%%%%%%%%%%%%%%%%%%%%%%%%%%%%%%
\subhead \the\sectionnumber.5\ The spherical Fourier transform
\endsubhead
Suppose that $s\geq t\geq 1$, and define
$$
\multline
\Phi_{\mu(q^{1+2\la})}(\mu(x)) =
{{(q^{-2\la}, q^{2+2\la};q^2)_\infty}\over
{(q^2,q^2,q^{-2\la}s^{-2},q^{2+2\la}s^{-2},
q^2t^{-2};q^2)_\infty}}
\phi_\la(\mu(x)) = \\
{{(q^{3+2\la}x^{\pm 1}/st;q^2)_\infty}\over{
(q^{4+2\la}t^{-2},q^{2+2\la}s^{-2},q^2s^{-2},
qx^{\pm 1}/st;q^2)_\infty}} \\ \times 
\, {}_8W_7(q^{2+2\la}/t^2;q{s\over t}x^{\pm 1},q^{2+2\la},
q^{2+2\la}, q^{2+2\la}/t^2;q^2,q^{-2\la}s^{-2})
\endmultline
\tag\eq{6150} 
$$
by an application of \cite{\GaspR, (III.24)}.
Here $\phi_\la$ is defined in \eqtag{6100}. For \eqtag{6150}
to be well-defined we need that 
$\phi_\la$ is invariant under interchanging 
$q^{1+2\la}$ and $q^{-1-2\la}$, or changing $\la$
into $-1-\la$. This is not obvious from \eqtag{6100}, 
but it can be obtained from Bailey's transformation for
a very-well-poised ${}_8\vp_7$-series \cite{\GaspR, (2.10.1)},
or directly from the proof of Lemma \thtag{690} as given in
Appendix B. The quantum group theoretic interpretation
of the invariance is that
the principal unitary, complementary and  strange series
representations are all obtained from the 
so-called principal series
representations which are equivalent for $\la$ and $-1-\la$,
see Burban and Klimyk \cite{\BurbK}, Masuda et al.
\cite{\MasuMNNSU}.  

We now define the spherical Fourier transform of
a $(s,t)$-spherical element $\xi=f(\rst)$, with $f$
continuous and compactly supported on the 
spectrum of $\pi_1(\rst)$, by
$$
\bigl( {\Cal F}\xi\bigr)(\si) =
\int_\R f(x) \overline{\Phi_\si(x)}\,
d\nu(x;q{t\over s}, q{s\over t};-qst|q^2)
= \bigl( {\Cal F}f)(\si),
\tag\eq{6155}
$$
which is, up to constant, formally equal to
$h\bigl( (\phi_\la(\rst))^\ast\xi\bigr)$ with
$\si=\mu(q^{1+2\la})$. The spherical Fourier
transform \eqtag{6155} is a special case of the
Askey-Wilson function transform as studied in \cite{\KoelSAW}.
There an inversion formula is obtained, which reduces to
the following theorem in this situation.

\proclaim{Theorem \thname{6160}} Assume $s\geq t\geq 1$.
The spherical Fourier transform
of \eqtag{6155} is inverted by
$$
\gather
f(x) = C \int_\R \bigl({\Cal F} f\bigr)(z)\, \Phi_z(x)\,
d\nu_{-q^{-4}s^{-2}}(z;q,q,qt^{-2};qs^2|q^2), \\
C = (qs)^{-1}
(q^2,q^2,q^2t^{-2},q^2t^{-2},q^2s^{-2};q^2)_\infty^2
\Theta(-q^2)^2 \Theta(-t^{-2})\Theta(-s^{-2}),
\endgather
$$
as an identity in $L^2(\R, d\nu(\cdot;qs/t,qt/s;-qst|q^2))$.
The notation $\Theta(a)=(a,q^2/a;q^2)_\infty$ for a
(normalised) theta product is used.
\endproclaim

\demo{Remark \thname{6170}} {\rm (i)} We refer
to \cite{\KoelSAW} for complete proofs and the appropriate
generalisation of Theorem \thtag{6160}.
Note that the spherical Fourier transform
is self-dual for the case $s=t=1$.
In compliance with the situation for the
compact quantum $SU(2)$ group case, we could call the
spherical functions for the case $s=t=1$ the 
continuous $q$-Legendre functions, cf. 
\cite{\AskeW}, \cite{\KoelSIAM}, 
\cite{\KoelAAM}, \cite{\KoelFIC}. For the $SU(1,1)$ group
the spherical Fourier transform is given by the Legendre
function transform, which is also known as the Mehler-Fock
transform, see \cite{\KoorJF}, \cite{\Vile, Ch. VI}, 
\cite{\VileK, Ch. 7}. 
So the transform \eqtag{6150} and its inverse
of Theorem \thtag{6160} is a two-parameter $q$-analogue of the
Legendre function (or Mehler-Fock) transform.

{\rm (ii)} We see that the support of the Plancherel measure
of the spherical Fourier transform 
is $[-1,1]$, which corresponds to all of the 
principal unitary series representations, plus  
the discrete set $\{ \mu(-q^{1-2k})\vert k\in\N\}$, which
corresponds to the strange series representations 
with $\la = -1 + k - i\pi/2\ln q$, $k\in\N$,  (and $\ep=0$). 
Note that the support is independent of $s$ and $t$. 
Indeed, the existence of a non-trivial kernel of $Y_sA$ in an
admissible irreducible unitary representation of $\U$ is
independent of $s$.

{\rm (iii)} For the limiting cases we obtain the big
$q$-Legendre function transform, which is studied and inverted in
\cite{\KoelSbig}, for $s\to\infty$, and the little $q$-Legendre
function transform, which is studied and inverted in \cite{\KakeMU},
\cite{\VaksK}, \cite{\Kake}, Appendix A, for $s,t\to\infty$.
In all these cases the support of the Plancherel measure is as in
part (ii) of this Remark. 
\enddemo

%%%%%%%%%%%%%%%%%%%%%%%%%%%%%%%%%%%%%%%%%%%%%%%%%%%%%%%%%%%%%%%%%%%%
\subhead \the\sectionnumber.6\ 
Other $K$-types\endsubhead 
In the previous subsections we have interpreted in a formal way a
special (2 continuous parameters) case of the Askey-Wilson function
transform as the spherical Fourier transform on the quantum
$SU(1,1)$ group. This is connected to the $(s,t)$-spherical part of
the Cartan decomposition in Theorem \thtag{3140}. It is also possible to
associate a Fourier transform related to the non-trivial $K$-types
in the Cartan decomposition of Theorem \thtag{3140}, and this allows us to
interpret a 4-parameter (2 continuous, 2 discrete) case of the
Askey-Wilson function transform on the quantum $SU(1,1)$ group.
Since the derivation lives on the same formal level we only shortly
discuss this more general case, and we refer to \cite{\KoelSAW} for
the precise analytic proof of the Askey-Wilson function transform.
For the limiting cases we refer to \cite{\KoelSbig} for the big
$q$-Jacobi function transform and to \cite{\Kake}, or Appendix A,
for the little $q$-Jacobi function transform. 
We stress that the formal results of this subsection have served as
the motivation for the analytic definition of the general 
Askey-Wilson, respectively big $q$-Jacobi, function transform in
\cite{\KoelSAW}, respectively \cite{\KoelSbig}. 

First we consider the action of the Casimir operator $\Om$. Since
$\Om$ is central, it preserves the Cartan decomposition. So we have,
cf. Theorem \thtag{3140},
$$
\Om. \bigl( \Ga^{(p)}_{i,j}(s,t)\, f(\rst)\bigr) = 
\Ga^{(p)}_{i,j}(s,t)\, (Lf)(\rst)
\tag\eq{6180}
$$
for some linear operator $L$. In the rest of this subsection
we take $p=2$, the other cases can be treated similarly.
In order to determine $L$ we proceed
by determining $A^\nu\Om$ modulo $\U(AY_t-\la_{-i}(t)) +
(Y_sA-\la_j(s))\U$ with $\la_j(t)$ as in Lemma \thtag{370}. This is done
as in Lemma \thtag{6110} using Koornwinder's method, and we find that it
is a linear combination of $A^{\nu+2}$, $A^\nu$ and $A^{\nu-2}$ with
explicit rational coefficients in $q^\nu$. Next we evaluate
\eqtag{6180} in $A^\nu$. Since $A^\nu$ is group-like in $\U$,
i.e. $\De(A^\nu)= A^\nu\otimes A^\nu$, this is a homomorphism. So
$$
\multline
\Ga^{(2)}_{i,j}(s,t)(A^\nu)\, (Lf)(\mu(q^\nu)) =
\bigl( \Ga^{(2)}_{i,j}(s,t)\, (Lf)(\rst)\bigr) (A^\nu) =
\bigl( \Ga^{(2)}_{i,j}(s,t)\, f(\rst)\bigr) (A^\nu\Om) = \\
\psi^+(q^k) \Ga^{(2)}_{i,j}(s,t)(A^{\nu+2}) \, f(\mu(q^{\nu+2}))+
\psi^0(q^k) \Ga^{(p)}_{i,j}(s,t)(A^\nu) \, f(\mu(q^\nu)) \\ +
\psi^-(q^k) \Ga^{(p)}_{i,j}(s,t)(A^{\nu-2}) \, f(\mu(q^{\nu-2}))
\endmultline
$$
for certain explicit 
rational functions $\psi^+,\psi^0,\psi^-$. Using
\eqtag{3135} and the homomorphism property we can calculate 
$\Ga^{(2)}_{i,j}(s,t)(A^\nu)$ explicitly in terms of finite
$q$-shifted factorials. For $j\in\{-i,1-i,\ldots,i\}$  we have
$\Ga^{(2)}_{i,j}(s,t)(A^\nu) = Cq^{-\nu i}
(q^{1+\nu}/st;q^2)_{i-j}(q^{1+\nu}s/t;q^2)_{i+j}$ for a
non-zero constant $C$ independent of $\nu$. In this way we
can determine $L$ in
terms if an Askey-Wilson difference operator. We find
$$
\gathered
q^{2i+1}(q-q^{-1})^2 \, L = 
\psi(z) (T_{q^2}-1) + \psi(z^{-1})(T_{q^{-2}}-1) + (1-q^{2i+1})^2,
\\ \psi(z) = 
{{(1-qtz/s)(1-q^{1+2i-2j}z/st)(1-qstz)(1-q^{1+2i+2j}sz/t)}\over
{(1-z^2)(1-q^2z^2)}}.
\endgathered
\tag\eq{6185}
$$
With respect to the measure
$d\nu_{-1}(\cdot;q^{1+2i+2j}s/t,qt/s,q^{1+2i-2j}/st;qst|q^2)$,
the operator $L$ is formally symmetric,
cf. Theorem \thtag{5200}, Proposition \thtag{6140}, and see  
\cite{\KoelSAW} for the general case.
See also \cite{\KoelAAM, \S 7} for the compact case.

To find the appropriate eigenfunctions of $L$ we have to determine
to which of the irreducible admissible unitary representations of
$\U$ as in \S 6.1 we formally can associate an element in the
corresponding part of the Cartan decomposition. So we have to
determine for which of the representations there exists eigenvectors
of $Y_sA$ and $Y_tA$ for the eigenvalues $\la_i(s)$ and $\la_j(t)$
as defined in Lemma \thtag{370} with $i,j\in\hf\Z$.
This is done in Appendix A, cf. \S 6.2, where
essentially the complete spectral analysis of $Y_sA$ in any of the
irreducible admissible unitary representations is described. Now for
the principal unitary, complementary and strange series we have an
eigenvector of $T^\bullet_{\la,\ep}(Y_sA)$ 
for the eigenvalue $\la_i(s)$
for every $\bullet\in\{P,C,S\}$ and $\la$ with $\ep\equiv
 i\mod\Z$. In order to have $\la_i(s)$ in the discrete
spectrum of $T^\bullet_{\la,\ep}(Y_sA)$ we need
$|sq^{2i}|>1$. In
the discrete series, $T^\pm_k(Y_sA)$ has an eigenvector for the
eigenvalue $\la_i(s)$ for only finitely many values of $k$.
Moreover, we need $k\equiv i\mod\Z$ and for $i<0$
the eigenvalue can occur
only in the negative discrete series and for $i>0$ it can
occur only in
the positive discrete series. Here we assume $|s|\geq q^{-1}$
so that we are dealing with essentially self-adjoint operators.
Let us denote such an eigenvector, if it exists, by $v_s(i)$. 
Assuming that the irreducible admissible unitary representation
$T^\bullet_{\la,\ep}$ or $T^\pm_k$
contains both $v_s(j)$ and $v_t(-i)$ we formally see that
$f_\la(X)= \langle
T^\bullet_{\la,\ep}(XA)v_t(-i),v_s(j)\rangle$ satisfies
\eqtag{3120} with $\la=\la_{-i}(t)$, $\mu=\la_j(s)$.
In case $j\in\{-i,1-i,\ldots,i\}$ we formally obtain
$$
f_\la = 
\sum_{n,m} \langle v_t(-i),e_m\rangle\langle e_n,v_s(j)\rangle
q^{m+\ep} T^\bullet_{\la,\ep;n,m} =
\Ga^{(2)}_{i,j}(s,t)\, \phi_{\la,\ep}^{(i,j)}(\rst)
$$
for some function $\phi_{\la,\ep}^{(i,j)}$.
Here $n,m$ run through $\Z$ if
$T^\bullet_{\la,\ep}$ is in the principal unitary,
complementary or strange series
representations and through $\Zp$ if 
$T^\bullet_{\la,\ep}=T^\pm_k$
is in the discrete series
representation. Evaluating in $A^\nu$ we can obtain
$\phi_{\la,\ep}^{(i,j)}$ from
$$
\sum_{n} \langle v_t(-i),e_n\rangle\langle e_n,v_s(j)\rangle
q^{(n+\ep)(1+\nu)} = C q^{-\nu i}
(q^{1+\nu}/st;q^2)_{i-j}(q^{1+\nu}s/t;q^2)_{i+j}
\, \phi_{\la,\ep}^{(i,j)}(\mu(q^\nu)).
$$
{}From this we can, in a similar way as for Lemma \thtag{690} determine 
$\phi_{\la,\ep}^{(i,j)}$ explicitly for the principal unitary,
complementary and strange series representations for the
diagonal case $i=-j$. An extension of Rahman's method in Appendix
B can be used to sum the other cases. (This is pointed out to
us by Hjalmar Rosengren.) In
case of the positive discrete series the sum runs through
$\Zp$ and the coefficients of $v_t(j)$ are Al-Salam and Chihara
polynomials, see \eqtag{515} and \cite{\KoelVdJ}.
The sum can then be evaluated using
the Poisson kernel for the Al-Salam and Chihara polynomials obtained
by Askey, Rahman and Suslov \cite{\AskeRS, (14.8)}. 
In the case $j\in\{-i,1-i,\ldots,i\}$ we find, up to a scalar
independent of $z$, in the case of the positive discrete series
$T^+_k$,
$$
\multline
\phi_{\la,\ep}^{(i,j)}(\mu(z)) =
z^{i-\la} {{(q^2z^2,q^{1+2j-2\la}stz,q^{1+2j-2\la}sz/t,
q^{1-2i-2\la}stz,q^{1+2i-2\la}sz/t;q^2)_\infty}\over
{(q^{2-2\la+2j}s^2z^2,q^{1+2j-2\la}stz,qsz/t,qz/st,
q^{1-2i-2j}tz/s;q^2)_\infty}} \\
\times \, {}_8W_7(q^{2j-2\la}s^2z^2;qstz,qsz/t,
q^{2j-2\la}s^2,q^{1+2j-2i}stz,q^{1+2i+2j}sz/t;q^2,
q^{-2\la-2j}s^{-2}),
\endmultline
$$
where $\la$ is equal to $-k$, see \S 6.1.
After application of \cite{\GaspR, (2.10.1)} we can relate
the right hand side with the asymptotically free solution
of $Lf(z)=[\hf+\la]^2 f(z)$ for $z\to 0$ as considered in
\cite{\KoelSAW}. Using the connection coefficient formula
\cite{\GaspR, (2.11.1)}, see \cite{\KoelSAW}, we can show
that the right hand side is indeed invariant under $z$ to
$z^{-1}$, and that it coincides, up to a constant, with the
Askey-Wilson function, i.e. the spherical function for the
Askey-Wilson function transform, since one of the connection
coefficients vanishes.
By comparing with \cite{\IsmaR}, \cite{\SuslJPA},
\cite{\SuslPP}, we see that these functions are indeed
solutions to the Askey-Wilson difference operator of
\eqtag{6185}.

Next we formally associate the corresponding Fourier
transform to the diagonal case by
$$
f\mapsto \hat f (\mu(q^{1+2\la})) = 
h\Bigl( \bigl( \Ga^{(2)}_{i,-i}(s,t)\phi_{\la,\ep}^{(i,-i)}\bigr)^\ast
\Ga^{(2)}_{i,-i}(s,t) f(\rst)\Bigr).
$$
The explicit expression of $h$, see Theorem \thtag{5200}, and of
$\phi_{\la,\ep}^{(i,-i)}$ can be used to
formally invert this transform by a
spectral analysis of the operator $L$ related to the Casimir element
$\Om$, assuming that the measure for the case $p=2$ in Theorem
\thtag{5200} is positive.
See \cite{\KoelSAW} for the rigorous analytic derivations and
for the explicit inversion formulas.
{}From this result we see that the Plancherel formula is supported
on the principal unitary series and the same discrete subset of the
strange series, cf. Remark \thtag{6170}(ii), plus on the discrete series
representations that allow a map $f_\la$ as before, i.e. 
for those discrete series representations that contain the appropriate
eigenvectors of $Y_sA$ and $Y_tA$. 
For all other
cases we can proceed in a similar fashion. 

The limit case $s\to\infty$ gives a 3-parameter family of big
$q$-Jacobi function transforms in this way, see \cite{\KoelSbig} for
the analytic proofs. Taking moreover $t\to\infty$ brings us back to
the case studied by Kakehi \cite{\Kake}. 

%%%%%%%%%%%%%%%%%%%%%%%%%%%%%%%%%%%%%%%%%%%%%%%%%%%%%%%%%%%%%%%%%%%%
%N E W   S E C T I O N%
%%%%%%%%%%%%%%%%%%%%%%%%%%%%%%%%%%%%%%%%%%%%%%%%%%%%%%%%%%%%%%%%%%%%
\head \newappA Spectral analysis of a doubly infinite
Jacobi matrix\endhead

In this subsection we give the spectral analysis of a doubly 
infinite
Jacobi matrix that arises from the second order $q$-difference
equation for the basic hypergeometric series ${}_2\vp_1$. In a way
the results can be viewed as the spectral analysis 
of a $q$-integral operator on $(0,\infty)$ 
with a basic hypergeometric series as kernel. 
The result covers in particular the little $q$-Jacobi function
transform as studied by Kakehi \cite{\Kake}, 
see also \cite{\KakeMU},
\cite{\VaksK}. The method of proof is similar to the one used in
\cite{\Kake}, so we are brief. The result is more general. 

%%%%%%%%%%%%%%%%%%%%%%%%%%%%%%%%%%%%%%%%%%%%%%%%%%%%%%%%%%%%%%%%%%%%
\subhead A.1\ Generalities\endsubhead
In this subsection we collect some 
generalities on the study of the symmetric
operator on the Hilbert space $\H$ defined by 
$$
L\, e_k = a_k \, e_{k+1} + b_k\, e_k + 
\overline{a_{k-1}}\, e_{k-1}, \qquad
a_k\not= 0, \ b_k\in\R,
\tag\Aeq{710}
$$
where 
$\{ e_k\}_{k\in\Z}$ is the standard orthonormal basis of $\H$. 
By replacing $e_k$ by $e^{i\psi_k}e_k$ with $\psi_k=\psi_{k+1} -
\arg a_k$ we see that we may assume that $a_k>0$, 
which we assume in this subsection from now on. 
We use the standard terminology and results as in Dunford and
Schwartz \cite{\DunfS, Ch.~XII}, see also 
Berezanski\u\i\ \cite{\Bere}, Kakehi \cite{\Kake}, 
Kakehi et al. \cite{\KakeMU}, Koelink and Stokman \cite{\KoelSbig}, 
Masson and Repka \cite{\MassR}, Rudin 
\cite{\Rudi}, Simon \cite{\Simo}.

The domain ${\Cal D}$ 
of $L$ is the dense subspace ${\Cal D}(\Z)$ of finite linear 
combinations of the basis elements $e_k$, then $L$ is a
densely defined symmetric operator.
Let $L^\ast$ with
domain ${\Cal D}^\ast$ be the adjoint and $L^{\ast\ast}$
with domain ${\Cal D}^{\ast\ast}$ the closure of $L$. 
The deficiency indices are equal since $L$ commutes with
complex conjugation and they are less than or equal to $2$, so that
$L$ has self-adjoint extensions. 

For any two vectors
$u=\sum_{k=-\infty}^\infty u(k)e_k$ and 
$v=\sum_{k=-\infty}^\infty v(k)e_k$ 
we define the Wronskian by
$$
[u,v](k) = a_k\bigl(u(k+1)v(k)-v(k+1)u(k)\bigr).
\tag\Aeq{717}
$$
Note that 
the Wronskian $[u,v](k)$ is independent of $k$ if
$Lu=xu$ and $Lv=xv$ for $x\in\C$. In this case we have that
$u$ and $v$ are linearly independent solutions if and only
if $[u,v]\not= 0$.

Associated with the
operator $L$ we have two Jacobi matrices $J^+$ and $J^-$
acting on $\ell^2(\Zp)$ with orthonormal basis $\{f_k\}_{k\in\Zp}$, which are 
given by
$$
\aligned
J^+\, f_k &= \cases
a_k \, f_{k+1} + b_k\, f_k + a_{k-1}\, f_{k-1},& 
\text{for $k\geq 1$,} \\
a_0 \, f_1 + b_0\, f_0, & \text{for $k=0$,}\endcases \\
J^-\, f_k &= \cases
a_{-k-1}\, f_{k+1} + b_{-k}\, f_k + a_{-k}\, f_{k-1},& 
\text{for $k\geq 1$,} \\
a_{-1} \, f_1 + b_0\, f_0, & \text{for $k=0$,}\endcases
\endaligned
\tag\Aeq{720}
$$
initially defined on ${\Cal D}(\Zp)$. Then $J^{\pm}$ are densely
defined symmetric operators with deficiency indices $(0,0)$
or $(1,1)$ corresponding to whether the associated Hamburger
moment problems is determinate or indeterminate, see 
Akhiezer \cite{\Akhi}, Berezanski\u\i\ \cite{\Bere}, 
Simon \cite{\Simo}.
Moreover, by \cite{\MassR, Thm.~2.1} the deficiency indices of
$L$ are obtained by summing the deficiency indices of $J^+$
and $J^-$.

{}From now on we assume that $a_k$ is bounded as $k\to -\infty$.
Then
$\lim_{m\to -\infty} [u,\bar v](m)=0$ for $u,v\in{\Cal D}^\ast$.
By \cite{\Bere, Thm.~1.3, p.~504} it follows that $J^-$ is
self-adjoint, hence the space $S_x^- = 
\{ u\mid Lu=xu,\, \sum_{k=-\infty}^N |u(k)|^2<\infty\, 
\text{for some $N\in\Z$}\}$ is one-dimensional for 
$x\in\C\backslash\R$ by \cite{\Akhi, \S 1.3}. Let us
say that $\Phi_x$ spans $S_x^-$ for $x\in\C\backslash\R$.
The similarly defined space $S_x^+ = 
\{ u\mid Lu=xu,\, \sum_{k=N}^\infty |u(k)|^2<\infty\, 
\text{for some $N\in\Z$}\}$ is either one-dimensional or
two-dimensional according to whether $J^+$ has deficiency
indices $(0,0)$ or $(1,1)$, see \cite{\Akhi, Ch.~1}.

For the purposes of this appendix, it suffices to 
consider the case that $J^+$ has deficiency indices 
$(0,0)$, which we will assume from now on. In particular 
$L$ is essentially self-adjoint, i.e. 
${\Cal D}^{\ast}={\Cal D}^{\ast\ast}$.
The closure $L^{\ast\ast}$ of $L$ satisfies the 
same formula \eqtag{710}, so we denote it also by $L$. 
We thus have that $S_x^+$ is one-dimensional, say spanned 
by $\phi_x$, and we have that $[\phi_x,\Phi_x]\not= 0$. 
Indeed, $\Phi_x$ cannot be in $S_x^+$ since  otherwise 
$Lu=iu$ would have a non-trivial solution in $\H$.

We also have to deal with possible non-real solutions
of \eqtag{710}. Note that if $\psi_x$ is a solution
of $L\psi_x=x\psi_x$, then so is $\overline{\psi_{\bar x}}$
defined by $\overline{\psi_{\bar x}}(k)= 
\overline{\psi_{\bar x}(k)}$, since we assume that the 
coefficients
$a_k$ and $b_k$ are real. Observe in particular, that 
$\overline{\phi_{\bar x}}$
and $\overline{\Phi_{\bar x}}$ are multiples of $\phi_x$ 
and $\Phi_x$, respectively, since the subspaces $S_x^{\pm}$ 
are one-dimensional.
 
Having the solutions $\Phi_x$ and $\phi_x$ at hand we can define
the Green kernel for $x\in\C\backslash\R$ by 
$$
G_x(k,l) = {1\over{[\Phi_x, \overline{\phi_{\bar x}}]}}
\cases 
\Phi_x(k)\overline{\phi_{\bar x}(l)}, 
&\text{for $k\leq l$,} \\
\Phi_x(l)\overline{\phi_{\bar x}(k)},
&\text{for $k\geq l$}\endcases
\tag\Aeq{725}
$$
and the operator $(G_xu)(k) = \langle u, \overline{G_x(k,\cdot)}
\rangle = \sum_{l=-\infty}^\infty u(l) G_x(k,l)$ for
$u\in\H$. Note that this is well defined, since 
$G_x(k,\cdot)\in\H$ for all $k$. 
Special cases of the following proposition are proved in 
\cite{\KakeMU}, \cite{\Kake} and \cite{\CiccKK}. 

\proclaim{Proposition \Athname{7100}} Let $L$ with 
domain ${\Cal D}$
be essentially self-adjoint, then the resolvent of 
the closure of $L$ 
is given by
$\bigl( (x-L)^{-1}u\bigr)(k) = 
\sum_{l=-\infty}^\infty u(l) G_x(k,l)$.
\endproclaim

Since $L$ with domain ${\Cal D}^{\ast\ast}$ is self-adjoint we have
the spectral decomposition, $L=\int_\R t\, dE(t)$, for
a unique projection valued measure $E$ on $\R$. This means that for
any vectors $u\in{\Cal D}^{\ast\ast}$, $v\in\H$ we 
have a complex measure $E_{u,v}$ on $\R$ such that 
$\langle Lu,v\rangle =\int_\R tdE_{u,v}(t)$, 
where $E_{u,v}(B)=\langle E(B)u,v\rangle$ for any Borel subset
$B\subset \R$, see \cite{\Rudi, Thm. 13.30}. 
The measure can be obtained from the resolvent by the
inversion formula, see \cite{\DunfS, Thm. XII.2.10}, 
$$
E_{u,v}\bigl( (x_1,x_2)\bigr) = 
\lim_{\de\downarrow 0} \lim_{\ep\downarrow 0} {1\over{2\pi i}}
\int_{x_1+\de}^{x_2-\de} \langle (x-i\ep-L)^{-1}u,v\rangle
- \langle (x+i\ep-L)^{-1}u,v\rangle\, dx,
\tag\Aeq{7110}
$$
where $x_1<x_2$. Combined with Proposition \thtag{7100}, 
we see that the Wronskian is crucial for the
structure of the spectral measure. In particular, if 
$\langle (x-L)^{-1}u,v\rangle$ is meromorphic in a subset of
$\C$ we find that $E_{u,v}$ has discrete mass points at
the real poles, and for a real pole $x_0$ we can rewrite
\eqtag{7110} as
$$
E_{u,v}\bigl( \{x_0\}\bigr) = 
{1\over{2\pi i}}
\oint_{(x_0)} \langle (x-L)^{-1}u,v\rangle\, dx 
\tag\Aeq{7120}
$$
where the contour is taken in the subset where
$\langle (x-L)^{-1}u,v\rangle$ is meromorphic and 
such that it encircles only the pole $x_0$.

Finally, observe that from the explicit formula for
the Green kernel \eqtag{725} we get for $x\in\R$ and $\ep>0$,
$$
\langle (x\pm i\ep-L)^{-1}u,v\rangle = \sum_{k\leq l}
{{\Phi_{x\pm i\ep}(k) \overline{\phi_{x\mp i\ep}(l)}}\over{
[\Phi_{x\pm i\ep},\overline{\phi_{x\mp i\ep}}]}} 
\bigl( u(l)\overline{v(k)} + u(k)\overline{v(l)}\bigr) 
(1-\hf\de_{k,l}).
\tag\Aeq{7125}
$$

%%%%%%%%%%%%%%%%%%%%%%%%%%%%%%%%%%%%%%%%%%%%%%%%%%%%%%%%%%%%%%%%%%%%
\subhead A.2\ The $q$-hypergeometric difference
equation\endsubhead
Consider the second-order $q$-diffe\-rence equation, see
\cite{\GaspR, exerc.~1.13}, in the following form
$$
(y+y^{-1})\, f(k) = (d-{{cq^{-k}}\over{dz}})\, f(k+1) 
+q^{-k}{{c+q}\over{dz}}\, f(k)
+ (d^{-1}-{{q^{1-k}}\over{dz}})\, f(k-1).
\tag\Aeq{7130}
$$
We assume that $d$ and $z$ are non-zero, and as usual we
take $0<q<1$. 
For the difference equation we have the following solutions
in terms of basic hypergeometric series;
$$
f_{\mu(y)}(k) = 
{}_2\vp_1\left( {{dy,d/y}\atop c};q,zq^k\right), 
\qquad c\not\in q^{-\Zp},\quad \mu(y)=\hf(y+y^{-1}), 
\tag\Aeq{7135}
$$
which is symmetric in $y$ and $y^{-1}$ and 
$$
F_y(k) = (dy)^{-k} \, {}_2\vp_1\left({{dy,qdy/c}\atop{qy^2}};q,
{{q^{1-k}c}\over{d^2z}}\right), \qquad y^2\not\in q^{-\N},
$$
so that we also have $F_{y^{-1}}(k)$ as a solution
to \eqtag{7130}. Here we use Jackson's transformation formula
\cite{\GaspR, (1.5.4)} to give a meaning to 
$f_{\mu(y)}(k)$ and $F_y(k)$
in case that $|zq^k|\geq 1$ and $|q^{1-k}c/d^2z|\geq 1$, 
respectively, for $z,d^2z/c\notin q^\Z$. 

These solutions are related by the expansion
$$
\aligned
f_{\mu(y)}(k) 
&= c(y) F_y(k) + c(y^{-1})F_{y^{-1}}(k), \\
c(x) &= {{(c/dy,d/y,dzy,q/dzy;q)_\infty}\over
{(y^{-2},c,z,q/z;q)_\infty}},
\endaligned
\tag\Aeq{7140}
$$
for $d,c,z\not= 0$, $|\arg (-z)|<\pi$, $c\not\in q^{-\Zp}$, 
$y^2\not\in q^\Z$, see \cite{\GaspR, (4.3.2)} and use the
theta-product identity \eqtag{237}. 

Next we consider the associated operator 
$$
\xi_k \mapsto  (d-{{cq^{-k}}\over{dz}})\, \xi_{k+1} 
+q^{-k}{{c+q}\over{dz}}\, \xi_k
+ (d^{-1}-{{q^{1-k}}\over{dz}})\, \xi_{k-1},
\tag\Aeq{7150}
$$
where we now assume that $\{\xi_k\}_{k\in\Z}$ is an
orthogonal basis of $\H$. We can now ask for what
values of the parameters $d$, $c$ and $z$ we can 
rewrite the operator as a symmetric operator of the form
as in \eqtag{710}. Inserting the orthonormal basis 
$e_k = \xi_k/\|\xi_k\|$
in \eqtag{7150} shows that we have to have 
that $(c+q)/dz\in\R$ and 
$$
{{ \|\xi_{k+1}\|^2}\over{\|\xi_k\|^2}} = 
{{\bar d^{-1} - q^{-k}/\bar d\bar z}\over{d-cq^{-k}/dz}}
= {1\over{|d|^2}} {{1-q^{-k}/\bar z}\over{1-cq^{-k}/d^2z}}.
\tag\Aeq{7160}
$$
Hence the right hand side of
\eqtag{7160} must be positive for all $k\in\Z$. Note that we
assume that the numerator and the denominator are non-zero
for all $k\in\Z$, in order not to reduce to the 
Jacobi matrix case.
So we assume $z, c/d^2z \not\in q^\Z$. 
On the other hand, if the 
right hand side of \eqtag{7160} is positive for all
$k\in\Z$ we can define $\|\xi_k\|$ recursively from
\eqtag{7160} and we find a symmetric operator of the
form \eqtag{710} assuming that $(c+q)/dz\in\R$. 
So positivity of the right hand
side of \eqtag{7160} and $(c+q)/dz\in\R$
are necessary and sufficient
for the mapping in \eqtag{7150} to be symmetric. 
In the following lemma we give an explicit 
description of the parameter domain which satisfy these 
conditions. The proof is similar to
the determination of unitary structures on irreducible 
principal series representations of
$\U$, see \cite{\MasuMNNSU, part II, \S 2}. 

\proclaim{Lemma \Athname{7170}} Assume $z, c/d^2z\not\in q^\Z$. 
The right hand side of the $q$-hypergeometric difference
equation \eqtag{7130} can be written as a symmetric operator
on $\H$ of the form \eqtag{710} if and only if $(c+q)/dz\in\R$
and one of
the following conditions holds:
{\rm (1)} $\bar z c=d^2 z$, or {\rm (2)} 
$z>0$, $c\not=d^2$ and
$zq^{k_0+1}<c/d^2<zq^{k_0}$, where $k_0\in \Z$ is such 
that $1<q^{k_0}z<q^{-1}$,
or {\rm (3)} $z<0$, $c\not=d^2$ and $c/d^2 >0$. 
In these cases the parameters
of \eqtag{710} are given by $b_k =q^{-k}(c+q)/dz$ and 
$$
a_k = \sqrt{ 
(1-{{q^{-k}}\over z})(1-{{cq^{-k}}\over{d^2z}})},
$$
after multiplying the basis $\{e_k\}$ with suitable phase factors.
\endproclaim

\demo{Remark}
For later purposes, we furthermore assume that $c, dz\in \R$.   
For the cases (2) and (3) this implies $c>0$, $d\in \R$, 
while for case (1) this implies that $c>0$ and 
$c=|d|^2$, since $|z|^2=(dz)^2$. Note that we may 
assume $k_0=0$ by replacing $k$ by $k+k_0$ in \eqtag{7130},
and replacing $z$ by $zq^{k_0}$.
\enddemo

We now consider the cases described in Lemma \thtag{7170}. The
symmetric operators are given by
$2Le_k=a_ke_{k+1}+b_ke_k+a_{k-1}e_{k-1}$
with $a_k$ and $b_k$ as in Lemma \thtag{7170}. Put 
$$
w(k) = e^{i\psi_k}|d|^{k} \sqrt{{{(cq^{1-k}/d^2z;q)_\infty}\over
{(q^{1-k}/\bar z;q)_\infty}}},
\tag\Aeq{7177}
$$
where $\psi_k\in\R$ are such that 
$\psi_{k+1}-\psi_k=\arg\bigl(d(1-q^{-k}/\bar z)\bigr)=
\arg\bigl(d(1-cq^{-k}/d^2z)\bigr)$ for all $k$. 
Then $u=wf=\sum_{k\in\Z} w(k) f(k)e_k$ is a solution to
$Lu=\mu(y)u$ if $f(k)$ is a solution to
the hypergeometric $q$-difference equation \eqtag{7130}.
Observe furthermore that 
for $k\rightarrow -\infty$ we have 
$|\arg(1-q^{-k}/\bar z)|={\Cal O}(q^{-k})$,
so that $\psi_{k+1}-\psi_k\,\, (\text{mod}\, 2\pi) 
\rightarrow \arg(d)$ as $k\rightarrow -\infty$.

\proclaim{Lemma \Athname{7179}} 
Let $c, dz\in \R$ and assume that the parameters 
satisfy the conditions as described in Lemma \thtag{7170}. 
Then the operator $L$ with domain ${\Cal D}(\Z)$ is essentially 
self-adjoint for $0<c\leq q^2$.
\endproclaim

\demo{Proof}
The $a_k$ are bounded for $k\rightarrow -\infty$, so it suffices 
to show that the Jacobi matrix $J^+$ associated to $L$ is 
essentially self-adjoint, see the previous subsection.
By \cite{\Bere,  Ch. VII, \S 1, Thm.~1.4, Cor.} we have that 
$J^+$ is essentially self-adjoint if
$a_k + a_{k-1}\pm b_k$ is bounded from above as $k\to\infty$
for a choice of the sign. 
Use $a_k=q^{-k}\sqrt{c/d^2z^2}-\hf(z+d^2z/c)+{\Cal O}(q^k)$,
$k\to\infty$, then the boundedness condition is satisfied if the
coefficient of $q^{-k}$ in $a_k + a_{k-1}\pm b_k$ is non-positive.
Since $c>0$, $dz\in\R$, this is
the case when $(1+q)\sqrt{c}\leq c+q$.
For $0<c\leq q^2$ the inequality holds. 
\qed\enddemo

{}From now on, we will assume throughout this appendix 
that $c,dz\in \R$, $0<c\leq q^2$, and that the parameters 
satisfy the conditions as described in Lemma \thtag{7170}.
Let $S_{x}^{\pm}$ be the eigenspaces of $L$ corresponding 
to the eigenvalue $x$ as defined in the previous subsection, 
and $[\cdot,\cdot]$ the Wronskian associated to $L$.

\proclaim{Lemma \Athname{7180}} The solution $wF_y$ spans 
$S^-_{\mu(y)}$
for $y\in\C$, $|y|<1$, and $wf_{\mu(y)}$ spans $S^+_{\mu(y)}$ for
$\mu(y)\in \C\backslash\R$. Furthermore, 
$$
[\overline{wf_{\overline{\mu(y)}}},
wF_y] = \hf\overline{c(\bar y^{-1})}(y-y^{-1})
$$
when $\mu(y)\in \C\backslash \R$, where $c(y)$ is defined 
in \eqtag{7140}. 
\endproclaim

\demo{Proof} Since $F_y(k)= {\Cal O}\bigl((dy)^{-k}\bigr)$
as $k\to-\infty$, the first statement follows from
\eqtag{7177}. 
Since $f_{\mu(y)}(k)={\Cal O}(1)$ as $k\to\infty$ and, by
\eqtag{237},
$$
w(k)=e^{i\psi_k}|d|^k\Bigl({{c\bar z} \over {d^2z}}\Bigr)^{k/2}
\sqrt{{{(\bar z q^k, d^2z/c, cq/d^2z;q)_\infty}\over{
(d^2zq^k/c,\bar z, q/\bar z;q)_\infty}}} \Rightarrow
|w(k)|={\Cal O}(c^{\hf k}), \ k\to\infty, 
$$
we have
$wf_{\mu(y)}\in S^+_{\mu(y)}$ for $|c|<1$. By Lemma 
\thtag{7179} and the generalities of the previous subsection
it follows that $wf_{\mu(y)}$ spans the one-dimensional
space $S^+_{\mu(y)}$.

It remains to calculate the Wronskian.
By \eqtag{7140} and the fact that $\overline{wF_{\overline{y}}}$
is a constant multiple of $wF_y$, see \S 7.1, we have
$$
[\overline{wf_{\overline{\mu(y)}}}, wF_y]=
\overline{c(\overline{y}^{-1})}
[\overline{wF_{\overline{y}^{-1}}}, wF_y].
$$
The lemma follows now from 
$$
\align
&[wF_y, \overline{wF_{\overline{y}^{-1}}}] =
\lim_{k\to-\infty} [wF_y, \overline{wF_{\overline{y}^{-1}}}](k) = \\
&\lim_{k\to -\infty}{{a_k}\over{2}}|w(k)w(k+1)|
\bigl(e^{i(\psi_{k+1}-\psi_k)}F_y(k+1)
\overline{F_{\overline{y}^{-1}}(k)}-
e^{i(\psi_k- \psi_{k+1})}F_y(k)
\overline{F_{\overline{y}^{-1}}(k+1)}\bigr)\\
&=\lim_{k\to -\infty}\frac{1}{2}|d|^{2k+1}
\bigl(e^{i\arg(d)}(yd)^{-k-1}{\overline{(\overline{y}^{-1}d)}}^{-k}
- e^{-i\arg(d)}(yd)^{-k}{\overline{(\overline{y}^{-1}d)}}^{-k-1}
\bigr)\\ &=\frac{1}{2}(y^{-1}-y). \qed
\endalign
$$
\enddemo 

We define for $x\in \C\backslash \R$, $\phi_x=wf_x$ and 
$\Phi_x=wF_y$, where $y$ is the unique element in the open unit 
disk such that $x=\mu(y)$. By Lemma \thtag{7180} and 
Proposition \thtag{7100},
we can give an expression of the resolvent $(x-L)^{-1}$ 
in terms of the
two functions $\overline{\phi_{\overline{x}}}$ and $\Phi_x$. 
In order to use \eqtag{7110} for the computation of the 
spectral measure of $L$, we have to calculate the limits 
as $\ep\downarrow 0$ in \eqtag{7125}.
Note that $\phi_{x\pm i\ep}\to wf_{x}$ as $\ep\downarrow 0$
for $x\in\R$.  For the asymptotic solution $\Phi_x$ we have 
to be more careful in computing the limit. 
For $x\in \R$ satisfying $|x|>1$ we have 
$\Phi_{x\pm i\ep} \rightarrow wF_y$ as 
$\ep\downarrow 0$, where 
$y\in (-1,1)\backslash \{0\}$ is such that $\mu(y)=x$. 
If $t\in [-1,1]$, then we put $x=\cos\chi=\mu(e^{i\chi})$ 
with $\chi\in [0,\pi]$, and then
$\Phi_{x-i\ep}\rightarrow wF_{e^{i\chi}}$ 
and $\Phi_{x+i\ep}\rightarrow wF_{e^{-i\chi}}$ 
as $\ep\downarrow 0$.  

Let us for the moment assume that the zeros of the
$c$-function of \eqtag{7140} are simple and do not
coincide with its poles. 
For the case $|x|< 1$, $x=\cos\chi=\mu(e^{i\chi})$ and 
$u,v\in {\Cal D}(\Z)$ we consider the limit
$$
\align 
&\lim_{\ep\downarrow 0} \langle (x-i\ep-L)^{-1}u,v\rangle
- \langle (x+i\ep-L)^{-1}u,v\rangle  
\tag\Aeq{7184}
\\ =& 
2\sum_{k\leq l}
\Bigl({{w(k)F_{e^{i\chi}}(k) \overline{w(l)f_{\cos\chi}(l)}}
\over{\overline{c(e^{i\chi})}(e^{-i\chi}-e^{i\chi})}} 
- {{w(k)F_{e^{-i\chi}}(k) \overline{w(l)f_{\cos\chi}(l)}}
\over{\overline{c(e^{-i\chi})}
(e^{i\chi}-e^{-i\chi})}}\Bigr)\\
&\qquad\qquad\qquad \times 
\bigl(u(l)\overline{v(k)}+u(k)\overline{v(l)}\bigr)
(1-\frac{1}{2}\delta_{k,l}).
\endalign
$$
Observe that the term within the big brackets can 
be written in the following two ways,
$$
\align
{{w(k)f_{\cos\chi}(k)\overline{w(l)f_{\cos\chi}(l)}}\over
{|c(e^{i\chi})|^2(e^{-i\chi}-e^{i\chi})}} &=
{{w(k)F_{e^{i\chi}}(k) \overline{w(l)f_{\cos\chi}(l)}}
\over{\overline{c(e^{i\chi})}(e^{-i\chi}-e^{i\chi})}} 
- {{w(k)F_{e^{-i\chi}}(k) \overline{w(l)f_{\cos\chi}(l)}}
\over{\overline{c(e^{-i\chi})}
(e^{i\chi}-e^{-i\chi})}}\\
&=\hf \overline{\Bigl({w(k)f_{\cos\chi}(k)w(l)f_{\cos\chi}(l)
\over {c(e^{i\chi})c(e^{-i\chi})[wF_{e^{-i\chi}}, 
wF_{e^{i\chi}}]}}\Bigr)}.\tag\Aeq{7186}
\endalign
$$
Here the first equality follows from the connection 
coefficient formula \eqtag{7140} and the fact that 
$|c(e^{i\chi})|=|c(e^{-i\chi})|$ by the conditions 
on the parameters. The second equality follows again 
by the connection coefficient formula \eqtag{7140}, 
and the fact that for $y\in \C\backslash \{0\}$ with 
$\mu(y)\in \C\backslash \R$, 
$$
{wF_y \over {[wF_y, \overline{wF_{\bar y^{-1}}}]}} =
{{\overline{wF_{\bar y}}} \over 
{[\overline{wF_{\bar y}}, 
\overline{wF_{\bar y^{-1}}}]}} \Rightarrow
{{2wF_y}\over{y^{-1}-y}} = 
{{\overline{wF_{\bar y}}} \over 
{\overline{ [wF_{\bar y}, wF_{\bar y^{-1}}]}}},  
$$ 
since $\overline{wF_{\bar y}}$ is a constant multiple of 
$wF_y$ and using the last step of the proof of Lemma \thtag{7180}.
{}From the second identity for the term in big brackets in 
\eqtag{7184} we see that it is symmetric in $k$ and $l$, 
so we can symmetrise the sum
in \eqtag{7184}. Using then the first identity for the term 
in big brackets in
\eqtag{7184} we obtain
$$
\lim_{\ep\downarrow 0} \langle (x-i\ep-L)^{-1}u,v\rangle
- \langle (x+i\ep-L)^{-1}u,v\rangle = 
2\sum_{k,l=-\infty}^\infty 
{{\overline{w(k)f_{\cos\chi}(k)}u(k) 
w(l)f_{\cos\chi}(l)\overline{v(l)}}\over{
|c(e^{i\chi})|^2(e^{-i\chi}-e^{i\chi})}}.
$$
Hence, with $dx = \frac{i}{2}(e^{i\chi}-e^{-i\chi})d\chi$, we
obtain for $0\leq \chi_1<\chi_2\leq \pi$ and $u,v\in {\Cal D}(\Z)$,
$$
E_{u,v}\bigl( (\cos\chi_2,\cos\chi_1)\bigr) = {1\over{2\pi}}
\int_{\chi_1}^{\chi_2} 
\bigl({\Cal F}u\bigr)(\cos\chi) 
\overline{\bigl({\Cal F}v\bigr)(\cos\chi)}
{{d\chi}\over{|c(e^{i\chi})|^2}},
$$
where 
$$
\bigl( {\Cal F}u\bigr)(x) = \langle u ,wf_x\rangle = 
\sum_{k=-\infty}^\infty u(k)\overline{w(k)f_x(k)}
\tag\Aeq{7185}
$$
for $u\in {\Cal D}(\Z)$ is the corresponding Fourier transform. 

Next we consider the case $|x|>1$, $x\in \R$, 
then we have from \eqtag{7125} and Lemma \thtag{7180} that
$$
\lim_{\ep\downarrow 0} 
\langle (x\pm i\ep-L)^{-1}u,v\rangle = 2\sum_{k\leq l}
{{w(k)F_y(k) \overline{w(l)f_{\mu(y)}(l)}}\over{
\overline{c( y^{-1})}(y^{-1}-y)}} 
\bigl( u(l)\overline{v(k)} + u(k)\overline{v(l)}\bigr) 
(1-\hf\de_{k,l})
$$
where $u,v\in {\Cal D}(\Z)$ and where
$y\in (-1,1)\backslash\{0\}$ is such that $x=\mu(y)$, 
provided that $y^{-1}$ is not a zero of $c(\cdot)$. 
It follows by the bounded convergence theorem that 
$E_{u,v}((x_1,x_2))=0$ when $(x_1,x_2)\cap [-1,1]=\emptyset$ 
and $(x_1,x_2)$ does not contain $x_0=\mu(y_0)$ with $y_0
\in (-1,1)$ a zero of the map $y\mapsto c(y^{-1})$.

Suppose now that $(x_1,x_2)\cap [-1,1]=\emptyset$  
and that $(x_1,x_2)$ contains exactly one point $x_0=\mu(y_0)$ 
with $c(y_0)=0$, where $y_0\in \R$ is such that $|y_0|>1$. 
Suppose furthermore that $y_0$ is a simple zero of $c(\cdot)$,
and that $c(y_0^{-1})\not=0$. Then it follows from 
\eqtag{7120} after the change of variable $x=\mu(y)$, that
$$
\align
&\langle E(\{x_0\})u,v\rangle = \langle E((x_1,x_2))u,v\rangle =\\
&\sum_{k\leq l}\text{Res}_{y=y_0^{-1}}\Bigl({ -1 \over 
\overline{c(y^{-1})}y}
\Bigr)w(k)F_{y_0^{-1}}(k)\overline{w(l)f_{x_0}(l)}
(u(l)\overline{v(k)}+u(k)\overline{v(l)})
(1-\frac{1}{2}\delta_{k,l}).
\endalign
$$
Now using the connection coefficient formula \eqtag{7140} and 
the fact that $c(y_0)=0$, we have 
$w(k)F_{y_0^{-1}}(k)=c(y_0^{-1})^{-1}w(k)f_{x_0}(k)$. 
Since $\overline{w(l)f_{x_0}(l)}w(k)f_{x_0}(k)$ is symmetric
in $k$ and $l$, cf. \eqtag{7186}, we can symmetrise to find  
$$
\langle E(\{x_0\})u,v\rangle =\sum_{k,l=-\infty}^{\infty}
\text{Res}_{y=y_0}
\Bigl({1 \over \overline{c(y)}c(y^{-1})y}\Bigr)
\overline{w(k)f_{x_0}(k)}u(k)
w(l)f_{x_0}(l)\overline{v(l)}.
$$

Observe that $(z,q/z;q)_\infty c(y)$ is real for $y\in\R$ and
that all zeros of the $c$-function outside the
unit disk are real. For parameters satisfying condition 
(2) or (3) of Lemma \thtag{7170} and $c>0$, $dz\in\R$, 
this is obvious. For parameters satisfying condition (1)
of lemma \thtag{7170} and $c>0$ and $dz\in\R$, 
this follows from the fact that $|d|=|c/d|<1$ since
$|d|^2=c\leq q^2$. It follows now easily that, 
for generic parameters, the support of the resolution of 
the identity $E$ of $L$ is given by $[-1,1]$, 
which is exactly the continuous spectrum of $L$, 
together with the discrete set $\{ x_0=\mu(y_0)\mid 
y_0\in \R\backslash [-1,1],\ c(y_0)=0 \}$, 
which is exactly the point spectrum of $L$, 
cf. \cite{\KakeMU} and \cite{\KoelSbig}.
These remarks prove a large
part of the following theorem, see \cite{\KakeMU}, 
\cite{\KoelSbig} for more details. 

\proclaim{Theorem \Athname{7190}} 
Consider $d,z$ as non-zero complex parameters such that $dz\in \R$.
Suppose that $0<c\leq q^2$, and that $z, c/d^2z\not\in q^{\Z}$.
Assume furthermore that the parameters satisfy one of 
the following three conditions: {\rm (1)} 
$\bar z c=d^2z$, or {\rm (2)}
$z>0$, $c\not=d^2$ and $zq^{k_0+1}<c/d^2<zq^{k_0}$, 
where $k_0\in \Z$ is such that $1<q^{k_0}z<q^{-1}$, or {\rm (3)}
$z<0$, $c\not=d^2$ and $c/d^2>0$. 
Consider the
following unbounded operators on $\H$ defined initially on the 
domain ${\Cal D}$ of
finite linear combinations of the 
orthonormal basis vectors $\{e_k\}_{k\in\Z}$;
$$
\gather
2L\, e_k = a_k \, e_{k+1} + b_k \, e_k 
+a_{k-1}\, e_{k-1},\\ 
b_k =q^{-k}(c+q)/dz\in \R, \qquad
a_k = \sqrt{ 
(1-{{q^{-k}}\over z})(1-{{cq^{-k}}\over{d^2z}})} >0.
\endgather
$$
Then $L$ is essentially self-adjoint, and
the closure $L^{\ast\ast}$ of the operator $L$ is
given by the same formula on ${\Cal D}^{\ast\ast}$. 
The spectral decomposition 
$L = \int_\R x\, dE(x)$ is given by
$$
\langle L u, v\rangle  = |(c,z,q/z;q)_\infty|^2 
\int_\R x \bigl({\Cal F}u\bigr) (x) \overline{\bigl( 
{\Cal F}v\bigr)(x)} \, d\nu(x;c/d, d;q/dz|q), 
\quad u\in {\Cal D}^{\ast\ast},\, v\in\H, 
$$
where the measure $d\nu (\cdot;a,b;d|q)$ is defined 
in \eqtag{590}, \eqtag{530}, and where 
the Fourier transform ${\Cal F}: \H \rightarrow L^2(\R, 
d\nu(\cdot;c/d,d;q/dz|q))$ is the unique continuous linear
map which coincides with the formulas \eqtag{7185}, 
\eqtag{7177} and \eqtag{7135} on ${\Cal D}$. 
\endproclaim

\demo{Remark} {\rm (i)} 
This theorem extends the result by Kakehi \cite{\Kake}
to a much larger parameter set. Kakehi's result corresponds
to case (3) with $z=-1$ and $c$ and $d$ in a discrete subset.
The proof is essentially the same. 

\noindent
{\rm (ii)} The Fourier transform ${\Cal F}\colon \H \to
L^2(\R; d\nu(\cdot;c/d, d;dz|q))$ is in fact 
an isometric isomorphism after scaling it by $(c,z,q/z;q)_\infty$.

\noindent
{\rm (iii)} It can be shown that the closure of 
$L$ has deficiency indices $(1,1)$ if one replaces the 
condition $0<c\leq q^2$ by $q^2<c<1$.
Indeed, since $|c|<1$ we have 
$wf_{x}\in S_{x}^+$ for $x\in \C\backslash \R$.
On the other hand,
$$
g_{x}(k) = q^kc^{-k} \, {}_2\vp_1\left( 
{{qdy/c, qd/cy}\atop{q^2/c}};q, zq^k\right), \qquad x=\mu(y), 
$$
is also a solution of the $q$-hypergeometric difference
equation \eqtag{7130}. 
Since $|w(k) g_{x}(k)|={\Cal O}(q^k|c|^{-k/2})$ 
as $k\to\infty$, 
we find a new $\ell^2$-solution of $Lf=xf$ as 
$k\to\infty$ for $1>c>q^2$,
which is linear independent of $wf_{x}$.
So $S^+_{x}$, $x\in\C\backslash\R$, is two-dimensional
for $q^2<c<1$, which implies that $L$ has deficiency 
indices $(1,1)$.
\enddemo

%%%%%%%%%%%%%%%%%%%%%%%%%%%%%%%%%%%%%%%%%%%%%%%%%%%%%%%%%%%%%%%%%%%%
%N E W   S E C T I O N%
%%%%%%%%%%%%%%%%%%%%%%%%%%%%%%%%%%%%%%%%%%%%%%%%%%%%%%%%%%%%%%%%%%%%
\head \newappB Summation formulas
\\ by Mizan Rahman\endhead

In this Appendix the proofs of Lemma \thtag{5150} 
and Lemma \thtag{690} are 
given. In both cases it involves an expression for the
Poisson kernel of the little $q$-Jacobi functions. The structure of
the proof is similar in both cases.
The proof of Lemma \thtag{5150} splits into two cases;
one for the absolutely
continuous part and one for the infinite set of discrete mass
points. This is treated in the first two subsections. 
The proof of Lemma \thtag{690}, treated in the last subsection, 
is similar to, but simpler than, the
proof of Lemma \thtag{5150} for the absolutely continuous case.

%%%%%%%%%%%%%%%%%%%%%%%%%%%%%%%%%%%%%%%%%%%%%%%%%%%%%%%%%%%%%%%%%%%%
\subhead B.1\ Proof of Lemma \thtag{5150} for the
absolutely continuous part\endsubhead
The idea of the proof is to write the product of two little
$q$-Jacobi functions as an infinite sum of Askey-Wilson
polynomials, and next to use an integral representation for
the Askey-Wilson polynomials. Interchanging summation and
integration gives a summable series as the integrand. The
resulting integral can then be evaluated, and after some
series manipulation we arrive at the desired result.
We give the proof in several steps. Recall that our
basic assumption is that the real parameters
$s$ and $t$ satisfy $|t|>1$, $|s|>1$.

First use \cite{\GaspR, (1.4.6)} to write the little
$q$-Jacobi function of \S 5.2 as
$$
\phi_n(x;s,t|q^2) = {{(-q^{2n}t^{-2};q^2)_\infty}\over{
(-q^{2n};q^2)_\infty}}
\, {}_2\vp_1\left( {{qtz/s,qt/sz}\atop{q^2s^{-2}}};
q^2,-q^{2n}t^{-2}\right), x=\mu(z)=\hf(z+z^{-1}),
$$
where we use the analytic continuation of the ${}_2\vp_1$-series
as in \cite{\GaspR, Ch.~4}. Using the theta-product identity
\eqtag{237} we see that we have to evaluate
$$
\multline
 \sum_{n=-\infty}^\infty \Bigl( {{uq^2}\over{s^2}}
\Bigr)^n \, {}_2\vp_1\left( {{qte^{i\th}/s,qte^{-i\th}/s}
\atop{q^2s^{-2}}}
;q^2,-q^{2n}t^{-2}\right)
\, {}_2\vp_1\left( {{qe^{i\th}/ts,qe^{-i\th}/ts}\atop{q^2s^{-2}}}
;q^2,-q^{2n}\right) \\
= Q_u(\cos\th) =
{{(-1,-q^2;q^2)_\infty}\over{(-q^2t^2,-t^{-2};q^2)_\infty}}
R_u(\cos\th;s,t|q^2). 
\endmultline
\tag\Beq{810}
$$

Recall the definition of the Askey-Wilson polynomials,
see \cite{\AskeW}, \cite{\GaspR, \S 7.5},
$$
p_m(x;a,b,c,d|q) =
{}_4\vp_3\left( {{q^{-m}, q^{m-1}abcd, ax,a/x}\atop{ab,ac,ad}};
q,q\right).
\tag\Beq{815}
$$
We can take the first step, which allows us
to separate the summation variable $n$ from the product of
the two ${}_2\vp_1$-series in \eqtag{810}.

\proclaim{Lemma \Bthname{820}}
For $|w|<1$ we have
$$
\multline
{}_2\vp_1\left( {{qte^{i\th}/s,qte^{-i\th}/s}\atop{q^2s^{-2}}}
;q^2,wt^{-2}\right)
\, {}_2\vp_1\left( {{qe^{i\th}/st,qe^{-i\th}/st}
\atop{q^2s^{-2}}};q^2,w \right) = \\ 
\sum_{m=0}^\infty
{{ (q^2s^{-2}e^{2i\th};q^2)_m}\over{(q^2;q^2)_m}}
\bigl( wq^{-1}se^{-i\th}t^{-1}\bigr)^m
p_m(t;{{qe^{i\th}}\over s},
{{qe^{i\th}}\over s},{{qe^{-i\th}}\over s},
{{qe^{-i\th}}\over s}|q^2).
\endmultline
$$
\endproclaim

\demo{Proof} Since $|w|<1$ and $|t|>1$, we have
$|wt^{-2}|<1$. Use the series representation of the two
${}_2\vp_1$-series to write the left hand side as an
absolutely convergent 
double sum. Next split off the power of $w$
in order to write the left hand side as
$$
\sum_{m=0}^\infty w^m
\sum_{k=0}^m {{(qe^{i\th}/st, qe^{-i\th}/st;q^2)_{m-k}}
\over{(q^2s^{-2},q^2;q^2)_{m-k}}}
{{(qte^{i\th}/s,qte^{-i\th}/s;q^2)_k}\over{
(q^2s^{-2},q^2;q^2)_k}}  t^{-2k}.
$$
Using elementary relations for the $q$-shifted factorials,
see \cite{\GaspR, \S 1.2}, we can rewrite this as
$$
\sum_{m=0}^\infty w^m 
{{(qe^{i\th}/st, qe^{-i\th}/st;q^2)_m}
\over{(q^2s^{-2},q^2;q^2)_m}}
\, {}_4\vp_3\left( {{q^{-2m}, s^2q^{-2m}, qte^{i\th}/s,
qte^{-i\th}/s}\atop{ q^2s^{-2}, q^{1-2m}ste^{i\th},
q^{1-2m}ste^{-i\th}}};q^2,q^2\right).
$$
Since the ${}_4\vp_3$-series is terminating and balanced
we can transform it using Sears's transformation
\cite{\GaspR, (2.10.4)} with $a$, $d$ replaced by
$qte^{i\th}/s$, $q^2s^{-2}$. Then the ${}_4\vp_3$-series
can be written as an Askey-Wilson polynomial, and
keeping track of the constants proves the lemma.
\qed\enddemo

Our next step is to use an integral representation for
the Askey-Wilson polynomial in Lemma \thtag{820}. There is a number
of ($q$-)integrals for the Askey-Wilson polynomial
available. 

\proclaim{Lemma \Bthname{830}} We have the integral
representation for the Askey-Wilson polynomial;
$$
\align
&p_m(t;{{qe^{i\th}}\over s},
{{qe^{i\th}}\over s},{{qe^{-i\th}}\over s},
{{qe^{-i\th}}\over s}|q^2)
=  A {1\over{2\pi}} \int_{-\pi}^\pi
{{( {{q\si}\over s}e^{-i(\th+\psi)};q^2)_m}\over{
({{q^3}\over{\si s^3}}e^{i(\th+\psi)};q^2)_m}} \Bigl(
{{qe^{i(\th+\psi)}}\over{\si s}}\Bigr)^m \\
& \qquad\qquad\times
{{({{kt}\over\si}e^{i\psi}, {{q^2\si}\over{kt}}e^{-i\psi},
k\si te^{-i\psi}, {{q^2}\over{k\si t}}e^{i\psi},
{{q^3}\over{\si s^3}} e^{i(\th+\psi)};q^2)_\infty}\over
{({q\over{\si s}}e^{i(\th+\psi)},{q\over{\si s}}e^{i(\th+\psi)},
{q\over{\si s}}e^{i(\psi-\th)},\si te^{-i\psi},
{\si\over t}e^{-i\psi};q^2)_\infty}}
\, d\psi, \\
&A = {{(q^2,qte^{i\th}/s,qe^{i\th}/st, qte^{-i\th}/s,
qe^{-i\th}/st, qte^{i\th}/s,qe^{i\th}/st;q^2)_\infty}\over
{(k,q^2/k,t^2k,q^2/t^2k,q^2e^{2i\th}/s^2,q^2s^{-2}, q^2s^{-2};
q^2)_\infty}}, 
\endalign
$$
where $\si$ and $k$ are free parameters such that there
are no zeros in the denominator. 
\endproclaim

\demo{Proof} We start with the integral representation for
a very-well-poised ${}_8\vp_7$-series given in
\cite{\GaspR, Exerc.~4.4, p.~122}, which can be proved by
a residue calculation. We use
\cite{\GaspR, Exerc.~4.4, p.~122}
in base $q^2$ and with the parameters $a=g=qe^{i\th}/\si s$,
$b=qe^{-i\th}/\si s$, $c=\si t$, $d=\si/t$,
$f=sq^{1-2m}e^{i\th}/\si$ and $h=sqe^{i\th}/\si$ where
$\si$ and $k$ are free parameters. This gives an integral
representation for a terminating very-well-poised
${}_8\vp_7$-series;
$$
\multline
{1\over{2\pi}} \int_{-\pi}^\pi
{{({{kt}\over\si}e^{i\psi}, {{q^2\si}\over{kt}}e^{-i\psi},
k\si te^{-i\psi}, {{q^2}\over{k\si t}}e^{i\psi},
{{q^{3+2m}}\over{\si s^3}} e^{i(\th+\psi)},
{s\over\si}q^{1-2m}e^{i(\th+\psi)};q^2)_\infty}\over
{({q\over{\si s}}e^{i(\th+\psi)},{q\over{\si s}}e^{i(\th+\psi)},
{q\over{\si s}}e^{i(\psi-\th)},\si te^{-i\psi},
{\si\over t}e^{-i\psi},{{sq}\over\si}e^{i(\th+\psi)}
;q^2)_\infty}}\, d\psi = \\
{{(k,\frac{q^2}{k},t^2k,\frac{q^2}{t^2k},\frac{q^2}{s^2}e^{2i\th},
\frac{q^2}{s^2},
ste^{i\th}q^{1-2m}, \frac{s}{t}e^{i\th}q^{1-2m},
q^2e^{2i\th},\frac{q^{2+2m}}{s^2};q^2)_\infty}\over
{(q^2,q\frac{t}{s}e^{i\th},\frac{q}{st}e^{i\th}, 
q\frac{t}{s}e^{-i\th},
\frac{q}{st}e^{-i\th}, q\frac{t}{s}e^{i\th},
\frac{q}{st}e^{i\th},qste^{i\th},
q\frac{s}{t}e^{i\th},q^{2-2m}e^{2i\th};q^2)_\infty}}
\\ \times
\, {}_8W_7(q^{-2m}e^{2i\th};
s^2q^{-2m}, s^2q^{-2m}e^{2i\th},
qte^{i\th}/s,qe^{i\th}/st,
q^{-2m};q^2, q^{2+2m}s^{-2}).
\endmultline
$$
Now we can apply Watson's formula \cite{\GaspR, (2.5.1)} with
$d=qte^{i\th}/s$, $e=qe^{i\th}/st$
to transform the terminating very-well-poised 
${}_8\vp_7$-series 
into a terminating balanced ${}_4\vp_3$-series. 
This shows that the ${}_8W_7$-series equals
$$
{{(q^{2-2m}e^{2i\th},s^2q^{-2m};q^2)_m}\over{(q^{1-2m}se^{i\th}/t,
q^{1-2m}ste^{i\th};q^2)_m}} \,
p_m(t;{{qe^{i\th}}\over s},
{{qe^{i\th}}\over s},{{qe^{-i\th}}\over s},
{{qe^{-i\th}}\over s}|q^2),
$$
the desired Askey-Wilson
polynomial. Collecting the results proves the lemma.
\qed\enddemo

In the proof the freedom to choose $k$ wisely is crucial,
but the $\si$-dependence is not essential. 
Combining Lemmas \thtag{820} and \thtag{830} gives
the following expression
$$ 
\aligned
& {}_2\vp_1\left( {{qte^{i\th}/s,qte^{-i\th}/s}\atop{q^2s^{-2}}}
;q^2,wt^{-2}\right)
\, {}_2\vp_1\left( {{qe^{i\th}/st,qe^{-i\th}/st}
\atop{q^2s^{-2}}};q^2,w \right)
= \\ &\qquad A {1\over{2\pi}} \int_{-\pi}^\pi
B(e^{i\psi})\,
{}_2\vp_1\left( {{q^2s^{-2}e^{2i\th},
 {{q\si}\over s}e^{-i(\th+\psi)}}\atop{
{{q^3}\over{\si s^3}}e^{i(\th+\psi)} }};q^2, 
{{we^{i\psi}}\over{\si t}}\right)\, d\psi , \\
& B(e^{i\psi}) =
{{({{kt}\over\si}e^{i\psi}, {{q^2\si}\over{kt}}e^{-i\psi},
k\si te^{-i\psi}, {{q^2}\over{k\si t}}e^{i\psi},
{{q^3}\over{\si s^3}} e^{i(\th+\psi)};q^2)_\infty}\over
{({q\over{\si s}}e^{i(\th+\psi)},{q\over{\si s}}e^{i(\th+\psi)},
{q\over{\si s}}e^{i(\psi-\th)},\si te^{-i\psi},
{\si\over t}e^{-i\psi};q^2)_\infty}}
\endaligned
\tag\Beq{840}
$$
and $A$ as in Lemma \thtag{830}. For $|w/\si t|<1$ 
interchanging integration and 
summation is justified. 
Note that \eqtag{840} also gives the analytic extension
of the left hand side to $w\in\C\backslash[1,\infty)$ using
the analytic continuation of the ${}_2\vp_1$-series in the
integrand, e.g. using \cite{\GaspR, (1.4.4)}, 
$$
\multline
{}_2\vp_1\left( {{q^2s^{-2}e^{2i\th},
 {{q\si}\over s}e^{-i(\th+\psi)}}\atop{
{{q^3}\over{\si s^3}}e^{i(\th+\psi)} }};q^2, 
{{we^{i\psi}}\over{\si t}}\right) = \\ 
{{(q^2e^{2i\th}/s^2,wqe^{-i\th}/st;q^2)_\infty}\over{
(q^3e^{i(\th+\psi)}/\si s^3, we^{i\psi}/\si t;q^2)_\infty}}
\, {}_2\vp_1\left(  {{qe^{i(\psi-\th)}/\si s, we^{i\psi}/\si t}
\atop{wqe^{-i\th}/st}};q^2, q^2e^{2i\th}s^{-2}\right).
\endmultline
\tag\Beq{850}
$$

We use \eqtag{840} with \eqtag{850} 
in \eqtag{810} and we interchange
summation and integration, which is easily justified. 
Then we have to evaluate
a sum where now the summand consists of one ${}_2\vp_1$-series.
This is done in the following lemma.

\proclaim{Lemma \Bthname{860}} For 
$\max (1, |s\si /q|)<|u|<s^2q^{-2}$ we have
$$
\multline
\sum_{n=-\infty}^\infty (uq^2s^{-2})^n
{}_2\vp_1\left( {{q^2s^{-2}e^{2i\th},
 {{q\si}\over s}e^{-i(\th+\psi)}}\atop{
{{q^3}\over{\si s^3}}e^{i(\th+\psi)} }};q^2, 
{{-q^{2n}e^{i\psi}}\over{\si t}}\right) = \\ 
{{(q^2,q^2e^{2i\th}/s^2, qe^{i(\th+\psi)}/\si su,
q\si e^{-i(\th+\psi)}/s, -\si ts^2e^{-i\psi}/u,
-uq^2e^{i\psi}/\si ts^2;q^2)_\infty}\over{
(e^{2i\th}/u, uq^2/s^2,\si se^{-i(\psi+\th)}/qu,
q^3 e^{i(\th+\psi)}/\si s^3, -e^{i\psi}/\si t,
-q^2\si te^{-i\psi};q^2)_\infty}}
\endmultline
$$
\endproclaim

\demo{Proof} Use the analytic continuation of \eqtag{850}
and interchange summation to write the left hand side as
$$
\multline
{{(q^2e^{2i\th}s^{-2};q^2)_\infty}\over{
(q^3e^{i(\th+\psi)}/\si s^3;q^2)_\infty}}
\sum_{m=0}^\infty {{(qe^{i(\psi-\th)}/\si s;q^2)_m}\over
{(q^2;q^2)_m}} \Bigl( {{q^2e^{2i\th}}\over{s^2}}\Bigr)^m \\
\times \sum_{n=-\infty}^\infty
{{(-q^{1+2n+2m}e^{-i\th}/st;q^2)_\infty}\over{
(-q^{2n+2m}e^{i\psi}/\si t;q^2)_\infty}} 
\Bigl( {{uq^2}\over{s^2}} \Bigr)^n.
\endmultline
$$
The inner sum can be evaluated by Ramanujan's
${}_1\psi_1$-summation formula \cite{\GaspR, (5.2.1)}
for $|s\si/q|<|u|<|s^2q^{-2}|$. The dependence on $m$
of the result is easy using the theta-product identity
\eqtag{237}. Explicitly, the inner sum equals
$$
\bigl( {{s^2}\over{uq^2}}\bigr)^m
{{(q^2,q\si e^{-i(\th+\psi)}/s, -uq^2e^{i\psi}/\si ts^2,
-e^{-i\psi}\si ts^2/u;q^2)_\infty}\over{
(uq^2/s^2,s\si e^{-i(\psi+\th)}/qu,-e^{i\psi}/\si t,
-q^2\si te^{-i\psi};q^2)_\infty}}.
$$
Then the inner sum to be evaluated reduces to
$$
\sum_{m=0}^\infty {{(qe^{i(\psi-\th)}/\si s;q^2)_m}\over
{(q^2;q^2)_m}} \bigl( {{e^{2i\th}}\over{u}}\bigr)^m  =
{{(qe^{i(\psi+\th)}/\si su;q^2)_\infty}\over{
(e^{2i\th}/u;q^2)_\infty}}
$$
by \cite{\GaspR, (1.3.2)} for $|u|>1$.
Collecting the intermediate results gives the
lemma.
\qed\enddemo

Combining \eqtag{840} and Lemma \thtag{860} gives an integral
representation for $Q_u(\cos\th)$ on $[-\pi,\pi]$
where the integrand consists of a quotient of eight
infinite $q$-shifted factorials in the numerator and
denominator. If we specialise one of the
the free parameters, $k=-q^2$,
this reduces to six infinite products in
the nominator and denominator. Explicitly,
$$
\aligned
Q_u(\cos\th) &= {B\over{2\pi}}\int_{-\pi}^\pi
{{(\frac{-q^2te^{i\psi}}{\si}, \frac{-\si e^{-i\psi}}{t}, 
\frac{qe^{i(\psi+\th)}}{\si su},
\frac{q\si e^{-i(\th+\psi)}}{s}, \frac{-\si ts^2e^{-i\psi}}{u},
\frac{-uq^2e^{i\psi}}{\si t s^2};q^2)_\infty}\over
{(\frac{qe^{i(\psi+\th)}}{\si s},
\frac{qe^{i(\psi+\th)}}{\si s}, \frac{qe^{i(\psi-\th)}}{\si s},
\si te^{-i\psi},\frac{\si}{t} e^{-i\psi},
\frac{\si se^{-i(\th+\psi)}}{qu};q^2)_\infty}}\, d\psi, \\
B &= {{(q^2,q^2,qte^{i\th}/s, qe^{i\th}/st, qte^{-i\th}/s,
qe^{-i\th}/st,qte^{i\th}/s, qe^{i\th}/st;q^2)_\infty}\over
{(-1,-q^2,-q^2t^2,-t^{-2},e^{2i\th}/u,q^2s^{-2},q^2s^{-2},
uq^2s^{-2};q^2)_\infty}}
\endaligned
\tag\Beq{870}
$$
The integral in \eqtag{870} is of the type considered
in \cite{\GaspR, \S\S 4.9-10} in base $q^2$ meaning that
we can evaluate the integral by residue calculus.
Indeed we may apply \cite{\GaspR, (4.10.8)} with the
specialisation $A=B=C=D=3$, $m=0$ and the twelve
other parameters given by
$a_1=-q^2t/\si$, $a_2=qe^{i\th}/\si su$, $a_3=-uq^2/\si ts^2$,
$b_1=-\si/t$, $b_2=q\si e^{-i\th}/s$, $b_3= -\si ts^2/u$,
$c_2=c_1= qe^{i\th}/\si s$, $c_3= qe^{-i\th}/\si s$,
$d_1=\si t$, $d_2=\si/t$, $d_3=\si se^{-i\th}/qu$.
Then the condition \cite{\GaspR, (4.10.2)} is
trivially satisfied, so that we may apply \cite{\GaspR, (4.10.8)}
to write the integral \eqtag{870} as a sum of three
${}_6\vp_5$-series. Due to $a_1b_1=q^2$, $a_3b_3=q^2$, these
${}_6\vp_5$-series reduce to ${}_4\vp_3$-series. So we have written
the integral in \eqtag{870} as a sum of three ${}_4\vp_3$-series
in this way. If we further specialise $u=q^{-2}$ two of these
${}_4\vp_3$-series reduce to ${}_3\vp_2$-series. The result is 
independent of $\si$ and it reads 
$$
\gathered
{1\over{2\pi}}\int_{-\pi}^\pi
{{(\frac{-q^2te^{i\psi}}{\si}, \frac{-\si e^{-i\psi}}{t}, 
\frac{q^3e^{i(\psi+\th)}}{\si s},
\frac{q\si e^{-i(\th+\psi)}}{s}, -\si ts^2q^2e^{-i\psi},
\frac{-e^{i\psi}}{\si t s^2};q^2)_\infty}\over
{(\frac{qe^{i(\psi+\th)}}{\si s},
\frac{qe^{i(\psi+\th)}}{\si s}, \frac{qe^{i(\psi-\th)}}{\si s},
\si te^{-i\psi},\frac{\si e^{-i\psi}}{t},
\si sqe^{-i(\th+\psi)};q^2)_\infty}}\, d\psi = \\
{{(-q^2t^2,\frac{q^3te^{i\th}}{s},
-s^{-2},-t^{-2},\frac{qe^{-i\th}}{st},-q^2s^2;q^2)_\infty}\over
{(q^2,\frac{qte^{i\th}}{s},\frac{qte^{i\th}}{s},
\frac{qte^{-i\th}}{s},t^{-2},
\frac{qse^{-i\th}}{t};q^2)_\infty}}
\, {}_3\vp_2\left( {{\frac{qte^{i\th}}{s}, 
\frac{qte^{-i\th}}{s},qste^{i\th}}
\atop{q^2t^2, \frac{q^3te^{i\th}}{s}}};q^2,q^2\right) \\
+ {{(-q^2,\frac{q^3e^{i\th}}{st},-t^{-2}s^{-2},-1,
q\frac{t}{s}e^{-i\th},-q^2s^2t^2;q^2)_\infty}\over{
(q^2,\frac{q}{st}e^{i\th}, \frac{q}{st}e^{i\th}, 
\frac{q}{st}e^{-i\th},t^2, qtse^{-i\th};
q^2)_\infty}}\, {}_3\vp_2\left( {{\frac{q}{st}e^{i\th}, 
\frac{q}{st}e^{-i\th},
q\frac{s}{t}e^{i\th}}\atop{\frac{q^3}{st}e^{i\th},
\frac{q^2}{t^2}}};q^2,q^2\right) \\
+ {{(-q^3ste^{-i\th},q^4,\frac{-qe^{-i\th}}{st},
\frac{-e^{i\th}}{stq}, s^{-2},
-stqe^{i\th};q^2)_\infty}\over{(q^2,q^2,q^2,q^2e^{-2i\th},
\frac{te^{i\th}}{qs},\frac{e^{i\th}}{stq};q^2)_\infty}}
\, {}_4\vp_3\left( {{q^2,q^2,q^2e^{-2i\th},q^2s^2}\atop{
q^4,q^3\frac{s}{t}e^{-i\th},q^3ste^{-i\th}}};q^2,q^2\right).
\endgathered
\tag\Beq{880}
$$

Combining \eqtag{880} with \eqtag{870} and \eqtag{810}
gives an expression for $R_{q^{-2}}(\cos\th;s,t|q^2)$
in terms of two ${}_3\vp_2$-series and one ${}_4\vp_3$-series.
In order to bring the very-well-poised ${}_8\vp_7$-series into
play we use Bailey's extension of Watson's transformation
formula \cite{\GaspR, (2.10.10)} in base $q^2$ with parameters
$a=q^2t^2$, $b=qste^{-i\th}$, $c=q^2$, $d=qte^{i\th}/s$,
$e=qte^{-i\th}/s$, $f=qste^{i\th}$. This gives the possibility
to write the ${}_4\vp_3$-series in \eqtag{880} as the sum of
a very-well-poised series as in Lemma \thtag{5150} and a
${}_3\vp_2$-series with the same parameters as the 
first ${}_3\vp_2$-series
on the right hand side of \eqtag{880};
$$
\multline
{{(-q^3ste^{-i\th},q^4,\frac{-qe^{-i\th}}{st},
\frac{-e^{i\th}}{stq}, s^{-2},
-stqe^{i\th};q^2)_\infty}\over{(q^2,q^2,q^2,q^2e^{-2i\th},
\frac{t}{sq}e^{i\th},\frac{e^{i\th}}{stq};q^2)_\infty}}
\, {}_4\vp_3\left( {{q^2,q^2,q^2e^{-2i\th},q^2s^2}\atop{
q^4,q^3\frac{s}{t}e^{-i\th},q^3ste^{-i\th}}};q^2,q^2\right) = \\
{{(q^2t^2,s^{-2},q^3ste^{i\th}, q^3te^{i\th}/s,
q^3te^{-i\th}/s,-stq^3e^{-i\th}, -e^{i\th}/qst,
-qste^{i\th},-qe^{-i\th}/st;q^2)_\infty}\over
{(q^2,q^2,q^4t^2,q^2e^{-2i\th},qte^{i\th}/s,qte^{-i\th}/s,
qste^{i\th},qte^{i\th}/s,e^{i\th}/qst;q^2)_\infty}} \\
\times {}_8W_7(q^2t^2;qte^{i\th}/s,qte^{-i\th}/s,qste^{i\th},
qste^{-i\th},q^2;q^2,q^2) \\
-{{(q^2t^2,q^3te^{i\th}/s,q^2s^2,s^{-2},-e^{i\th}/qst,
-stq^3e^{-i\th},-qste^{i\th},-qe^{-i\th}/st;q^2)_\infty}\over
{(q^2,q^3ste^{-i\th},sqe^{-i\th}/t,
qte^{i\th}/s,qte^{i\th}/s,qte^{-i\th}/s,qste^{i\th},
e^{i\th}/qst;q^2)_\infty}}
\\ \times \, {}_3\vp_2\left( 
{{qte^{i\th}/s, qte^{-i\th}/s,qste^{i\th}}\atop
{q^2t^2,q^3te^{i\th}/s}};q^2,q^2\right).
\endmultline
\tag\Beq{890}
$$
If we now use \eqtag{890} in \eqtag{880} and \eqtag{870}
we obtain the ${}_8W_7$-series on the right hand side of
Lemma \thtag{5150} with the factor in front using straightforward
manipulations of $q$-shifted factorials. 

It remains to show that the remaining terms can be summed
explicitly. For this we first consider the factor in front
of the first ${}_3\vp_2$-series in \eqtag{880} after having
plugged in \eqtag{890} for the ${}_4\vp_3$-series. This factor is
$$
\multline
{{(q^3te^{i\th}/s;q^2)_\infty}\over{(q^2,qte^{i\th}/s,
qte^{i\th}/s,qte^{-i\th}/s,qse^{-i\th}/t,qste^{i\th},
q^3ste^{-i\th},e^{i\th}/qst,t^{-2};q^2)_\infty}} \\
\Bigl( (-q^2t^2,-t^{-2},-s^{-2},-q^2s^2,\frac{qe^{-i\th}}{st},
qste^{i\th},q^3ste^{-i\th},\frac{e^{i\th}}{qst};q^2)_\infty
- \\ (q^2t^2,t^{-2},s^{-2},q^2s^2,\frac{-qe^{-i\th}}{st},
-qste^{i\th},-q^3ste^{-i\th},\frac{-e^{i\th}}{qst};q^2)_\infty
\Bigr) \\
=
{{(q^3te^{i\th}/s,qse^{i\th}/t,-1,-q^2,-q^2s^2t^2,
-s^{-2}t^{-2};q^2)_\infty}\over{
(qte^{i\th}/s,qste^{i\th},qste^{-i\th},qe^{i\th}/st,
q^2,t^{-2};q^2)_\infty}}
\endmultline
$$
where we used the theta-product identity 
\eqtag{5175} with $\la=qst$, $\mu=-e^{i\th}$,
$x=-qt/s$ and $\nu=-qst$, 
see \cite{\GaspR, Exerc. 2.16}. 
Having used this identity 
we see that the resulting sum of two
${}_3\vp_2$-series in \eqtag{880} after having 
applied \eqtag{890}
can be summed by the non-terminating version of
the Saalsch\"utz summation formula 
\cite{\GaspR, (2.10.12)} with $e=q^2t^2$, $f=tq^3e^{i\th}/s$ in the
form
$$
\multline
{{(t^{-2}, q{t\over s}e^{i\th}, q{t\over s}e^{-i\th},
qste^{i\th},{q^3\over{st}}e^{i\th};q^2)_\infty}\over{
(t^2,{q\over{st}}e^{i\th},{q\over{st}}e^{-i\th},q{s\over t}e^{i\th},
q^3{t\over s}e^{i\th};q^2)_\infty}}
{}_3\vp_2\left( {{\frac{q}{st}e^{i\th}, 
\frac{q}{st}e^{-i\th},
q\frac{s}{t}e^{i\th}}\atop{\frac{q^3}{st}e^{i\th},
\frac{q^2}{t^2}}};q^2,q^2\right) + \\
{}_3\vp_2\left( {{\frac{qte^{i\th}}{s}, 
\frac{qte^{-i\th}}{s},qste^{i\th}}
\atop{q^2t^2, \frac{q^3te^{i\th}}{s}}};q^2,q^2\right) 
= {{(t^{-2},q^2,q^2e^{2i\th},q^2/s^2;q^2)_\infty}\over{
({q\over{st}}e^{i\th}, {q\over{st}}e^{-i\th}, q{s\over t}e^{i\th},
q^3{t\over s}e^{i\th};q^2)_\infty}}.
\endmultline
$$ 
This gives the term with $q$-shifted factorials in
Lemma \thtag{5150}. \qed 

%%%%%%%%%%%%%%%%%%%%%%%%%%%%%%%%%%%%%%%%%%%%%%%%%%%%%%%%%%%%%%%%%%%%
\subhead B.2\ Proof of Lemma \thtag{5150} for the
infinite set of discrete mass points\endsubhead
Since the radius of convergence for $u$ of $R_u(x;s,t|q^2)$
depends on $x$,  we cannot use the result obtained
for $x\in [-1,1]$ in the previous subsection to obtain
the value for $x=\mu(q^{1-2k}st)$, cf. \cite{\KoelV, \S 6}. 
In order to prove Lemma \thtag{5150} for the infinite set of
discrete mass points we take $x=\mu(q^{1-2k}st)$
in the definition of $R_u(x;s,t|q^2)$ in Proposition \thtag{5110}. 
For this argument we can use \cite{\GaspR, (1.4.4), (1.4.5)}
to rewrite the little $q$-Jacobi function as 
$$
\align
 {}_2\vp_1& \left( {{ -q^{2-2k},
-q^{2k}s^{-2}t^{-2}}\atop{q^2s^{-2}}};q^2, -q^{2n}\right)\\  = &
{{(-q^{2k}s^{-2}t^{-2},q^{2+2n-2k};q^2)_\infty}\over
{(q^2s^{-2},-q^{2n};q^2)_\infty}} \ {}_2\vp_1
\left( {{-q^{2-2k}t^2,-q^{2n}}\atop{q^{2+2n-2k}}};q^2,
-{{q^{2k}}\over{s^2t^2}}\right) \\
=& {{(-q^{2k}s^{-2},q^{2+2n-2k};q^2)_\infty}\over{
(q^2s^{-2},-q^{2n};q^2)_\infty}} \ {}_2\vp_1\left(
{{-q^{2n}t^{-2},-q^{2-2k}}\atop{q^{2+2n-2k}}};q^2, -{{q^{2k}}
\over{s^2}}\right),
\endalign
$$
This shows the $q$-Bessel coefficient behaviour 
of the little $q$-Jacobi function at these
discrete mass points. In case the absolute value of the argument
is greater than one we can use Jackson's transformation of
a ${}_2\vp_1$-series to a ${}_2\vp_2$-series, see 
\cite{\GaspR, (1.5.4)}, to give the analytic extension 
which respects the $q$-Bessel coefficient behaviour. 
Hence, we find
$$
\multline
R_u(\mu(-stq^{1-2k});s,t|q^2) 
= {{(-t^{-2},-q^2t^2,-q^{2k}s^{-2},
-q^{2k}s^{-2}t^{-2};q^2)_\infty}\over{
(-1,-q^2,q^2s^{-2},q^2s^{-2};q^2)_\infty}} \\
\sum_{n=-\infty}^\infty
\bigl( {{uq^2}\over{s^2}}\bigr)^n 
{{(q^{2+2n-2k};q^2)_\infty}\over{(-q^{2n};q^2)_\infty}}
\ {}_2\vp_1
\left( {{-q^{2-2k}t^2,-q^{2n}}\atop{q^{2+2n-2k}}};q^2,
-{{q^{2k}}\over{s^2t^2}}\right)  \\ \times
{{(q^{2+2n-2k};q^2)_\infty}\over{(-q^{2n}t^{-2};q^2)_\infty}}
{}_2\vp_1\left(
{{-q^{2n}t^{-2},-q^{2-2k}}\atop{q^{2+2n-2k}}};q^2, -{{q^{2k}}
\over{s^2}}\right)
\endmultline
\tag\Beq{8100}
$$

Recall the identity, see 
\cite{\KoelITSF, Prop.~2.2, with $n=0$, $w=tx^{-1}y^{-1}$},
$$
\multline
\sum_{m=-\infty}^\infty w^m
{{(q^{m+1};q)_\infty}\over{(aq^m;q)_\infty}}
\, {}_2\vp_1\left( {{aq^m,b}\atop{q^{m+1}}};q,-x\right)
{{(q^{m+1};q)_\infty}\over{(cq^m;q)_\infty}}
\, {}_2\vp_1\left( {{cq^m,d}\atop{q^{m+1}}};q,-y\right)  = \\
{{(q,q;q)_\infty}\over{(a,c;q)_\infty}}
\sum_{p=0}^\infty (-x)^p {{(a,b;q)_p}\over{(q,q;q)_p}}
\, {}_2\vp_1\left( {{q^{-p},d}\atop{q^{1-p}a^{-1}}};q,
-{{qy}\over{aw}}\right)
\, {}_2\vp_1\left( {{q^{-p},c}\atop{q^{1-p}b^{-1}}};q,
-{{qw}\over{bx}}\right).
\endmultline
$$
After specialising $a=-y/w$, $c=-x/w$ we can use 
the $q$-Chu-Vandermonde sums 
\cite{\GaspR, (1.5.2), (1.5.3)} to sum the two terminating
${}_2\vp_1$-series in the summand on the right hand side. 
This gives 
$$
\multline
\sum_{m=-\infty}^\infty 
{{w^m(q^{m+1};q)_\infty}\over{(-yq^m/w;q)_\infty}}
\, {}_2\vp_1\left( {{-\frac{yq^m}{w},b}\atop{q^{m+1}}};q,-x\right)
{{(q^{m+1};q)_\infty}\over{(-xq^m/w;q)_\infty}}
\, {}_2\vp_1\left( {{-\frac{xq^m}{w},d}
\atop{q^{m+1}}};q,-y\right)  \\ =
{{(q,q;q)_\infty}\over{(-y/w,-x/w;q)_\infty}}
\, {}_2\vp_1\left( {{-yd/w,-bx/w}\atop q};q,w\right), 
\endmultline
$$
valid for $1>|w|>|xy|$. Indeed, using \eqtag{238} we see that
$$
{{(q^{m+1};q)_\infty}\over{(-xq^m/w;q)_\infty}}
\, {}_2\vp_1\left( {{-\frac{xq^m}{w},d}
\atop{q^{m+1}}};q,-y\right) = 
\cases {\Cal O}(1), & m\to\infty , \\
{\Cal O}((-x)^{-m}), & m\to -\infty. \endcases
$$
Use this identity in base $q^2$ with  $b= -q^{2-2k}t^2$,
$d= -q^{2-2k}$, $x= q^{2k}s^{-2}t^{-2}$,
$y= q^{2k}s^{-2}$, and $w=s^{-2}$, so we specialise
$u=q^{-2}$. Note that $s^{-2}<1$ since $|s|>1$, 
and $|xy|= q^{4k}s^{-4}t^{-2}< q^2s^{-2}$, since we have to
evaluate at the discrete mass point $\mu(-q^{1-2k}st)$
so that $|q^{1-2k}st|>1$. So we may use this identity 
to see that the sum in \eqtag{8100} for $u=q^{-2}$ 
equals
$$
s^{-2k} {{(q^2,q^2;q^2)_\infty}\over{
(-q^{2k},-q^{2k}t^{-2};q^2)_\infty}}
\, {}_2\vp_1\left( {{q^2,q^2}\atop
{q^2}};q^2, s^{-2}\right) = {{s^{-2k}}\over{1-s^{-2}}} 
{{(q^2,q^2;q^2)_\infty}\over{
(-q^{2k},-q^{2k}t^{-2};q^2)_\infty}}. 
$$
Using this gives 
$$
R_{q^{-2}}(\mu(-stq^{1-2k});s,t|q^2) = 
 {{(-t^{-2},-q^2t^2,-q^{2k}s^{-2},
-q^{2k}s^{-2}t^{-2},q^2,q^2;q^2)_\infty}\over{
(-1,-q^2,q^2s^{-2},q^2s^{-2},-q^{2k},
-q^{2k}t^{-2};q^2)_\infty}} {{s^{-2k}}\over{1-s^{-2}}}.
$$
Note that specialising $e^{i\th}=-q^{1-2k}st$ in
Lemma \thtag{5150} results in the same answer after
manipulating theta-products, since the ${}_8W_7$-series is
not singular for this value and the factor in front of the
${}_8W_7$-series is zero. Taking into account the first term gives
the result. \qed

%%%%%%%%%%%%%%%%%%%%%%%%%%%%%%%%%%%%%%%%%%%%%%%%%%%%%%%%%%%%%%%%%%%%
\subhead B.3\ Proof of Lemma \thtag{690}\endsubhead 
The proof of Lemma \thtag{690} is similar to the proof
given in the first subsection of this appendix, but it is
simpler. It again uses the Askey-Wilson polynomials and
a corresponding $q$-integral representation. However,
we have to distinguish between the principal unitary
series on the one hand and the complementary and
the strange series on the other hand in some
derivations. 

We first observe that $|\langle v^\bullet_t,e_n\rangle|$
behaves as $|s/q|^n$ as $n\to-\infty$ and as $|s|^{-n}$
as $n\to\infty$. This follows from the results from
Appendix A. It follows that the 
doubly infinite sum of Lemma \thtag{690} is absolutely convergent
in the annulus $|q/st|<|z|<|st/q|$. 

The analogue of Lemma \thtag{820} is the
following. 

\proclaim{Lemma \Bthname{8110}} With the notation of
\S 6.2 and \S 6.3 and with the assumptions
of Lemma \thtag{690} and assuming $n\leq 0$ we have 
$$
\langle v^\bullet_t,e_n\rangle \langle e_n,v^\bullet_s\rangle = 
(st)^n q^{-2n}\sum_{m=0}^\infty 
q^{-2nm}\, p_m(q^{1+2\la};q,qs^{-2},qt^{-2},q|q^2)
$$
for $\bullet\in\{ P,C,S\}$ using the notation \eqtag{815}
for the Askey-Wilson polynomials.
\endproclaim

\demo{Proof}
The proof is similar to the proof of Lemma \thtag{820}, but we
have to choose the right form of the ${}_2\vp_1$-series
in order to have the ${}_4\vp_3$-series balanced. 
We start with, cf. \eqtag{651},
$$
\multline
\langle v^\bullet_t,e_n\rangle \langle e_n,v^\bullet_s\rangle
= (st)^n q^{-4n(1+\Re \la)} \sqrt{
{{(q^{-2\la+2n},q^{-2\bar\la +2n};q^2)_\infty}\over{
(q^{2\la+2n+2},q^{2\bar\la +2n+2};q^2)_\infty}} } \\
\times {}_2\vp_1\left( {{q^{2+2\la}t^{-2},q^{2+2\la}}\atop
{q^2t^{-2}}};q^2,q^{-2n-2\la}\right)
{}_2\vp_1\left( {{q^{2+2\bar\la}s^{-2},q^{2+2\bar\la}}\atop
{q^2s^{-2}}};q^2,q^{-2n-2\bar\la}\right).
\endmultline
\tag\Beq{8115}
$$
Note that for $n\leq 0$ both ${}_2\vp_1$-series are
absolutely convergent, since $|q^{2+2\la}|\leq q$
for $\la$ as in Lemma \thtag{690}. As in the proof of 
Lemma \thtag{820} we rewrite the product of the two ${}_2\vp_1$-series
as
$$
\sum_{m=0}^\infty q^{-2m(n+\bar\la)} 
{{(q^{2+2\bar\la}s^{-2},q^{2+2\bar\la};q^2)_m}\over
{(q^2,q^2s^{-2};q^2)_m}} 
\, {}_4\vp_3\left( {{ q^{2+2\la}t^{-2},q^{2+2\la},q^{-2m},
q^{-2m}s^2}\atop{q^2t^{-2},q^{-2m-2\bar\la}s^2, q^{-2m-2\bar\la}}};
q^2, q^{-4\Re\la}\right)
$$
and the ${}_4\vp_3$-series is balanced if $\Re\la=-\hf$, i.e.
for $\la$ corresponding to the principal unitary series. We can 
apply Sears's transformation \cite{\GaspR, (2.10.4)} with
$a$ and $d$ specialised to $q^{2+2\la}$ and $q^2t^{-2}$ to
rewrite this as, using $\Re\la=-\hf$, 
$$
\sum_{m=0}^\infty q^{-2mn}
{}_4\vp_3\left( {{ q^{-2m}, q^{2+2m}t^{-2}s^{-2},
q^{2+2\la}, q^{-\la}}\atop{q^2t^{-2},q^2s^{-2}, q^2}};
q^2, q^2\right)
\tag\Beq{8117}
$$
and the ${}_4\vp_3$-series is the Askey-Wilson polynomial
as in the lemma. Note that the square root of $q$-shifted
factorials in \eqtag{8115} reduces to $1$ for $\Re\la=-\hf$.
This proves the lemma for the principal unitary series. 

For the complementary series and the strange series 
we use Heine's transformation formula \cite{\GaspR, (1.4.6)}; 
$$
\multline
{}_2\vp_1\left( {{q^{2+2\bar\la}s^{-2},q^{2+2\bar\la}}\atop
{q^2s^{-2}}};q^2,q^{-2n-2\bar\la}\right) = \\
{{(q^{2+2\bar\la-2n};q^2)_\infty}\over
{(q^{-2n-2\bar\la};q^2)_\infty}}
\, {}_2\vp_1\left( {{q^{-2\bar\la}, q^{-2\bar\la}s^{-2}}\atop
{q^2s^{-2}}};q^2, q^{2+2\bar\la-2n}\right).
\endmultline
$$
The theta product identity \eqtag{237} 
can now be used to rewrite \eqtag{8115} as
$$
\multline
\langle v^\bullet_t,e_n\rangle \langle e_n,v^\bullet_s\rangle
= 
(ts)^n q^{-2n(1-\la+\bar\la)} \sqrt{
{{(q^{-2\la+2n},q^{2\bar\la+2n+2};q^2)_\infty}\over
{(q^{2\la+2n+2},q^{-2\bar\la+2n};q^2)_\infty}} }\\
\times {}_2\vp_1\left( {{q^{2+2\la}t^{-2},q^{2+2\la}}\atop
{q^2t^{-2}}};q^2,q^{-2n-2\la}\right)
{}_2\vp_1\left( {{q^{-2\bar\la},q^{-2\bar\la}s^{-2}}\atop
{q^2s^{-2}}};q^2,q^{2+2\bar\la-2n}\right).
\endmultline
\tag\Beq{8120}
$$
The product of the two ${}_2\vp_1$-series can be written as
$$
\sum_{m=0}^\infty q^{-2m(n+\la)}
{{(q^{2+2\la}t^{-2}, q^{2+2\la};q^2)_m}\over
{(q^2,q^2t^{-2};q^2)_m}}
\, {}_4\vp_3\left( {{q^{-2\bar\la}, q^{-2\bar\la}s^{-2}, 
q^{-2m}, q^{-2p}t^2}\atop{q^2s^{-2},q^{-2m-2\la}t^2, 
q^{-2m-2\la}}};q^2, q^{2+2\bar\la-2\la}\right)
$$
and the ${}_4\vp_3$-series is balanced if $q^{2\la}=q^{2\bar\la}$.
Since this is the case for the complementary series and
the strange series, we apply once more Sears's transformation
\cite{\GaspR, (2.10.4)} with
$a$ and $d$ specialised to $q^{-2\bar\la}$ and $q^2s^{-2}$ 
to rewrite this sum as \eqtag{8117}
for $\la$ satisfying $q^{2\la}=q^{2\bar\la}$. Observe that
$q^{2\la}=q^{2\bar\la}$ makes the square root
of $q$-shifted factorials in \eqtag{8120}
equal to $1$. This proves the
result for the complementary series and the strange series.
\qed\enddemo

We next employ the following
$q$-integral representation for the
Askey-Wilson polynomials of \eqtag{815}; 
$$
\aligned
p_m(x;a,b,c,d|q) &= \bigl(A(x;a,b,c;d|q)\bigr)^{-1} 
{{(bc;q)_m}\over{(ad;q)_m}} \\ &\times
\int_{qx/d}^{q/xd} 
{{(dux,du/x,abcdu/q;q)_\infty}\over{(adu/q,bdu/q,cdu/q;q)_\infty}}
{{(q/u;q)_m}\over{(abcdu/q;q)_m}} \Bigl( {{adu}\over q}\Bigr)^m
\, d_qu, \\
A(x;a,b,c;d|q) &= {{q(1-q)}\over{d(x-x^{-1})}}
{{(x^2,x^{-2},q,ab,ac,bc;q)_\infty}\over
{(ax,a/x,bx,b/x,cx,c/x;q)_\infty}},
\endaligned
\tag\Beq{8130}
$$
where the $q$-integral is defined by, cf. \eqtag{451}, 
$$
\int_a^b f(x)\, d_q x=  \int_0^b f(x)\, d_q x - 
\int_0^a f(x)\, d_q x, \quad
\int_0^c f(x)\, d_q x = (1-q)c\sum_{n=0}^\infty f(cq^k)\, q^k,
$$
cf. \S 4.2. 
The $q$-integral representation is in Exercise~7.34 
of \cite{\GaspR}. The proof consists of rewriting the
$q$-integral into the form \cite{\GaspR, (2.10.19)}, which
can be done in such a way that the very-well-poised 
${}_8\vp_7$-series is terminating. The
terminating very-well-poised 
${}_8\vp_7$-series can then be rewritten as a terminating
balanced ${}_4\vp_3$-series by Watson's transformation
formula \cite{\GaspR, (2.5.1)}, which can be recognised as
an Askey-Wilson polynomial in the form \eqtag{815}. 

Using \eqtag{8130} and Lemma \thtag{8110} we have for
$n\leq 0$ 
$$
\multline
\langle v^\bullet_t,e_n\rangle \langle e_n,v^\bullet_s\rangle
= (ts)^n q^{-2n}
A \int_{q^{2+2\la}}^{q^{-2\la}} 
{{(uq^{2+2\la}, uq^{-2\la},q^2s^{-2}t^{-2}u;q^2)_\infty}\over
{(s^{-2}u,ut^{-2}, u;q^2)_\infty}} \\ \times \,
{}_2\vp_1\left( {{q^2s^{-2}t^{-2}, q^2/u}\atop{q^2s^{-2}t^{-2}u}};
q^2, uq^{-2n}\right)\, d_{q^2}u, 
\endmultline
\tag\Beq{8140}
$$
where $A^{-1}=A(q^{1+2\la};qs^{-2},q,qt^{-2};q|q^2)$. 
Interchanging $q$-integrating
and summation is justified, since all sums are absolutely
convergent for $n\leq 0$ because $|q^{2+2\la}|\leq q$ and
$|q^{-2\la}|\leq q$. 
We can next use \eqtag{8140} for the expression 
to extend the left hand side to the case 
$n>0$ by using the analytic continuation of the 
${}_2\vp_1$-series in the $q$-integral of \eqtag{8140}. 
The analytic continuation of the ${}_2\vp_1$-series is
given by \cite{\GaspR, (4.3.2)}, a formula we already used 
for the $c$-function expansion of \eqtag{7140}. 
In this particular case the second term vanishes 
for $n\in\Z$ and the
factor in front of the remaining ${}_2\vp_1$-series can
be simplified using the theta product identity \eqtag{237}.
This gives 
$$
{}_2\vp_1\left( {{q^2s^{-2}t^{-2}, q^2/u}\atop{q^2s^{-2}t^{-2}u}};
q^2, uq^{-2n}\right) = 
(q^2s^{-2}t^{-2})^n 
{}_2\vp_1\left( {{q^2s^{-2}t^{-2}, q^2/u}\atop{q^2s^{-2}t^{-2}u}};
q^2, uq^{2n}\right). 
\tag\Beq{8135}
$$
For another way to see this, rewrite the
left hand side using Heine's transformation 
\cite{\GaspR, (1.4.5)} to recognise the $q$-Bessel coefficient
behaviour. Next \eqtag{238} provides the requested relation. 
 
\proclaim{Lemma \Bthname{8150}} For $z$ in the annulus
$|q/st|<|z|<|st/q|$ we have 
$$
\multline
\sum_{n=-\infty}^0 (ts)^n q^{-2n} (qz)^n
{}_2\vp_1\left( {{q^2s^{-2}t^{-2}, q^2/u}\atop{q^2s^{-2}t^{-2}u}};
q^2, uq^{-2n}\right) + \\
\sum_{n=1}^\infty (ts)^{-n} (qz)^n
{}_2\vp_1\left( {{q^2s^{-2}t^{-2}, q^2/u}\atop{q^2s^{-2}t^{-2}u}};
q^2, uq^{2n}\right) \\ = 
{{(1-q^2/s^2t^2)}\over{(1-q/zst)
(1-qz/st)}}  {{(q^4/s^2t^2, q^2, qu/zst,
quz/st;q^2)_\infty}\over{ (q^3/zst, 
q^3z/st, q^2u/s^2t^2, u;q^2)_\infty}}. 
\endmultline
$$
\endproclaim

\demo{Proof} We can write the left hand side as
$$
\sum_{j=0}^\infty {{(q^2s^{-2}t^{-2}, q^2/u;q^2)_j}\over
{(q^2, q^2s^{-2}t^{-2}u;q^2)_j}}u^j 
\Bigl( \sum_{n=-\infty}^0 (ts)^n q^{-2n-2nj} (qz)^n +
\sum_{n=1}^\infty (ts)^{-n} (qz)^n q^{2nj}\Bigr).
$$
Both sums in parantheses are geometric sums and absolutely
convergent for $z$ in the annulus. The sums in 
parantheses equal
$$
\multline
{{(1-q^{2+4j}/s^2t^2)}\over{(1-q^{2j+1}/zst)
(1-q^{2j+1}z/st)}} = \\
{{(1-q^2/s^2t^2)}\over{(1-q/zst)
(1-qz/st)}} 
{{(1-q^{2+4j}/s^2t^2)\, (q/zst,
qz/st;q^2)_j}\over{ (1-q^2/s^2t^2)
\, (q^3/zst,q^3z/st;q^2)_j}}.
\endmultline
$$
Plugging this back gives the left hand side as a very-well-poised
${}_6\vp_5$-series;
$$
\multline 
{{(1-q^2/s^2t^2)}\over{(1-q/zst)
(1-qz/st)}} \, 
{}_6W_5(q^2/s^2t^2; qz/st, q/zst, q^2/u;
q^2, u) 
= \\ {{(1-q^2/s^2t^2)}\over{(1-q/zst)
(1-qz/st)}}  {{(q^4/s^2t^2, q^2, qu/zst,
quz/st;q^2)_\infty}\over{ (q^3/zst, 
q^3z/st, q^2u/s^2t^2, u;q^2)_\infty}},
\endmultline
$$
where we used the summation formula \cite{\GaspR, (2.7.1)}.
\qed\enddemo

Combining \eqtag{8140} and Lemma \thtag{8150} we see that for
$z$ in the annulus as in Lemma \thtag{8150} we have
$$
\multline
\sum_{n=-\infty}^\infty  
\langle v^\bullet_t,e_n\rangle \langle e_n,v^\bullet_s\rangle
q^nz^n = A {{(1-q^2/s^2t^2)}\over{(1-q/zst)
(1-qz/st)}} {{(q^4/s^2t^2, q^2;q^2)_\infty}\over{ (q^3/zst, 
q^3z/st;q^2)_\infty}} \\ \times
\int_{q^{2+2\la}}^{q^{-2\la}} 
{{(uq^{2+2\la}, uq^{-2\la}, qu/zst,
quz/st;q^2)_\infty}\over
{(s^{-2}u,ut^{-2}, u, u;q^2)_\infty}} \, d_{q^2}u.
\endmultline
$$
The $q$-integral is of the same type as used for the
Askey-Wilson polynomial, and it can be explicitly evaluated
in terms of a very-well poised ${}_8\vp_7$-series by
\cite{\GaspR, (2.10.19)}. The $q$-integral is equal to 
$$
\multline
q^{-2\la}(1-q^2) {{(q^2,q^{-4\la},
q^{2+4\la},q^2/t^2,q^2/s^2,q^2/s^2t^2,
q^{1-2\la}/zst, q^{1-2\la}z/st;q^2)_\infty}\over{(q^{2+2\la}/s^2,
q^{2+2\la}/t^2,q^{2+2\la}, q^{-2\la}/s^2, q^{-2\la}/t^2,
q^{-2\la},q^{-2\la}, q^{2-2\la}/s^2t^2;q^2)_\infty}} \\
\times\, {}_8W_7(q^{-2\la}/s^2t^2;q^{-2\la}/s^2, q^{-2\la}/t^2,
q^{-2\la}, qz/st,q/zst;q^2,q^{2+2\la}).
\endmultline
$$
Plugging this back in and using the value for
$A$ proves Lemma \thtag{690} for $z$ in the annulus 
$|q/st|<|z|<|st/q|$. 
\qed

%%%%%%%%%%%%%%%%%%%%%%%%%%%%%%%%%%%%%%%%%%%%%%%%%%%%%%%%%%%%%%%%%%%%
%R E F E R E N C E S%
%%%%%%%%%%%%%%%%%%%%%%%%%%%%%%%%%%%%%%%%%%%%%%%%%%%%%%%%%%%%%%%%%%%%

\enddocument